\documentclass[12pt]{smfart}

\usepackage{smfenum}
\usepackage[francais]{babel}
\usepackage{smfthm}
\usepackage{amssymb}
\usepackage{eucal}
\usepackage{mathrsfs}

\newcommand{\C}{\mathbf{C}}
\newcommand{\R}{\mathbf{R}}
\newcommand{\N}{\mathbf{N}}
\newcommand{\Q}{{\mathbf{Q}}_p}
\newcommand{\Z}{{\mathbf{Z}}_p}

\newcommand{\Qp}{\Q}
\newcommand{\Cp}{\mathbf{C}_p}
\newcommand{\Zp}{\Z}

\newcommand{\ZZ}{\mathbf{Z}}
\newcommand{\NN}{\N}
\newcommand{\RR}{\R}
\newcommand{\OO}{\mathcal{O}}
\renewcommand{\O}{\OO_L}

\newcommand{\Qpbar}{\overline{\mathbf{Q}}_p}

\newcommand{\BK}{\mathrm{B}(\Z)}

\newcommand{\G}{\mathrm{GL}_2(\Qp)}
\newcommand{\B}{\mathrm{B}(\Qp)}
\newcommand{\K}{\mathrm{GL}_2(\Zp)}
\newcommand{\g}{\mathrm{Gal}(\Qpbar/\Qp)}
\newcommand{\gal}{\g}
\newcommand{\galinf}{\mathrm{Gal}(\Qpbar/F_\infty)}

\renewcommand{\=}{\overset{\text{d\'ef}}{=}}
\newcommand{\val}{\mathrm{val}}
\newcommand{\eps}{\varepsilon}
\renewcommand{\projlim}{\varprojlim}

\newcommand{\Fil}{\mathrm{Fil}}
\renewcommand{\geq}{\geqslant}
\renewcommand{\leq}{\leqslant} 
\newcommand{\pr}{\mathrm{pr}}
\newcommand{\nrm}{|\cdot|}

\newcommand{\bcris}{\mathbf{B}_{\rm cris}} 
\newcommand{\bcontp}{\widetilde{\mathbf{B}}^+_{\rm rig}} \newcommand{\btrigplus}{\widetilde{\mathbf{B}}^+_{\rm rig}} 
\newcommand{\btrig}{\widetilde{\mathbf{B}}^\dagger_{\rm rig}} 
\newcommand{\bdR}{\mathbf{B}_{\rm dR}}  
\newcommand{\bfont}{\mathbf{B}}
\newcommand{\bdag}{\mathbf{B}^\dagger}
\newcommand{\btdag}{\widetilde{\mathbf{B}}^\dagger}
\newcommand{\bplus}{\mathbf{B}^+}

\newcommand{\calE}{\mathscr{E}}
\newcommand{\calO}{\mathscr{O}_{\mathscr{E}}}
\newcommand{\calR}{\mathscr{R}}

\newcommand{\dcris}{\mathrm{D}_{\mathrm{cris}}}
\newcommand{\ddR}{\mathrm{D}_{\mathrm{dR}}}
\newcommand{\dplus}{\mathrm{D}^+}
\newcommand{\dsen}{\mathrm{D}_{\mathrm{sen}}}
\newcommand{\dfont}{\mathrm{D}}
\renewcommand{\ddag}{\mathrm{D}^\dagger}
\newcommand{\drig}{\mathrm{D}^\dagger_{\mathrm{rig}}}
\newcommand{\dsharp}{\mathrm{D}^\sharp}
\newcommand{\mfil}{\mathrm{M}}
\newcommand{\nwach}{\mathrm{N}}
\newcommand{\breuil}{\mathrm{B}}
\newcommand{\st}{\mathrm{St}}

\author[L. Berger]{Laurent Berger}
\address{CNRS \& IH\'ES \\
Le Bois-Marie\\
35 route de Chartres\\
91440 Bures-sur-Yvette \\ 
France}
\email{laurent.berger@ihes.fr}
\urladdr{www.ihes.fr/\~{}lberger/}

\author[C. Breuil]{Christophe Breuil}
\address{CNRS \& IH\'ES\\
Le Bois-Marie\\
35 route de Chartres\\
91440 Bures-sur-Yvette \\ 
France}
\email{breuil@ihes.fr}
\urladdr{www.ihes.fr/\~{}breuil/}

\date{janvier 2006}

\title[Repr\'esentations potentiellement cristallines de $\mathrm{GL}_2(\Qp)$]{Sur quelques repr\'esentations potentiellement cristallines de $\mathrm{GL}_2(\Qp)$}  

\alttitle{On some potentially crystalline representations of $\mathrm{GL}_2(\Qp)$}

\subjclass{11F}

\NumberTheoremsIn{subsection}

\begin{document}

\begin{abstract}
On associe aux repr\'esentations $p$-adiques irr\'eductibles de $\gal$ de dimension $2$ devenant cristallines sur une extension ab\'elienne de $\Q$ des espaces de Banach $p$-adiques $\breuil(V)$ munis d'une action lin\'eaire continue unitaire de $\G$. Lorsque $V$ est de plus $\varphi$-semi-simple, on utilise le $(\varphi,\Gamma)$-module et le module de Wach de $V$ pour montrer que la repr\'esentation $\breuil(V)$ est non nulle, topologiquement irr\'eductible et admissible.
\end{abstract}

\begin{altabstract}
To each $2$-dimensional irreducible $p$-adic representation of $\gal$ which becomes crystalline over an abelian extension of $\Q$, we associate a Banach space $\breuil(V)$ endowed with a linear continuous unitary action of $\G$. When $V$ is moreover $\varphi$-semi-simple, we use the $(\varphi,\Gamma)$-module and the Wach module associated to $V$ to show that the representation $\breuil(V)$ is nonzero, topologically irreducible and admissible.
\end{altabstract} 

\maketitle

\setcounter{tocdepth}{2}
\tableofcontents

\setlength{\baselineskip}{17pt}

\section{Introduction}

\subsection{Introduction}\label{intro}

Soit $p$ un nombre premier et $n\in \N$. Dans la recherche d'une
correspondance \'eventuelle entre (certaines) repr\'esentations
$p$-adiques $V$ de $\g$ de dimension $n$ et (certaines)
repr\'esentations $p$-adiques $\breuil(V)$ de ${\rm GL}_n(\Q)$, un des
cas importants \`a regarder est certainement celui o\`u $n=2$, $V$ est absolument irr\'eductible, devient cristalline sur une extension ab\'elienne de $\Q$ et est $\varphi$-semi-simple. Cette
derni\`ere condition signifie que le Frobenius $\varphi$ sur le
$\varphi$-module
filtr\'e $\dcris(V)$ associ\'e par Fontaine \`a $V$ est semi-simple.
La repr\'esentation $V$ a des poids de Hodge-Tate distincts $i_1< i_2$ et, si
l'on note ${\rm Alg}(V)$ la repr\'esentation alg\'ebrique de $\G$ de
plus haut poids $(i_1,i_2-1)$ et ${\rm Lisse}(V)$ la repr\'esentation
lisse irr\'eductible de $\G$ associ\'ee par la correspondance locale de
Hecke \`a la repr\'esentation de Weil d\'eduite de $V$ par \cite{F5}, la repr\'esentation
$\breuil(V)$ est simplement le compl\'et\'e $p$-adique de la
repr\'esentation localement alg\'ebrique ${\rm Alg}(V)\otimes {\rm
Lisse}(V)$ par rapport \`a un r\'eseau stable par $\G$ et de type fini
sous l'action de $\G$. Notons que ${\rm
Lisse}(V)$ est toujours ici une repr\'esentation de la s\'erie principale. Ainsi, $\breuil(V)$ est un espace de Banach $p$-adique muni d'une action continue unitaire de $\G$ (i.e. laissant une norme invariante). Notons que, lorsque ${\rm Lisse}(V)$ est de dimension $1$, cette d\'efinition de $\breuil(V)$ doit \^etre modifi\'ee.

Le probl\`eme est qu'il n'est pas du tout \'evident qu'un tel r\'eseau
existe, ou, de mani\`ere \'equivalente, que $\breuil(V)$ soit non nul. Dans \cite[th\'eor\`eme 1.3]{Br1}, la non nullit\'e de $\breuil(V)$ est
d\'emontr\'ee lorsque $V$ est cristalline et $i_2-i_1<2p$ (essentiellement), et dans
\cite[th\'eor\`eme 1.3.3]{Br2},
son admissibilit\'e (au sens de \cite{ST3}) et son irr\'eductibilit\'e
topologique (avec une condition suppl\'ementaire pour cette
derni\`ere). La m\'ethode repose sur le calcul de la r\'eduction d'une
boule unit\'e de $\breuil(V)$ modulo l'id\'eal maximal des coefficients.

Lorsque $V$ est de dimension $2$, absolument irr\'eductible mais cette
fois semi-stable non cristalline, $\breuil(V)$ est d\'efini dans \cite{Br2}
et \cite{Br3} et des conjectures analogues (non nullit\'e,
admissibilit\'e, etc.) formul\'ees et tr\`es partiellement
d\'emontr\'ees. P. Colmez dans \cite{Co2} a vu que la th\'eorie des
$(\varphi,\Gamma)$-modules de Fontaine permettait de d\'emontrer
\'el\'egamment ces conjectures dans le cas semi-stable en construisant
un mod\`ele de la restriction au Borel sup\'erieur de la
repr\'esentation duale $\breuil(V)^*$ \`a partir du
$(\varphi,\Gamma)$-module de $V$.

Il \'etait donc naturel de regarder si un tel mod\`ele existait aussi
dans le cas des repr\'esenta\-tions $V$ potentiellement cristallines ci-dessus et s'il permettait de d\'emontrer les conjectures de non nullit\'e, d'irr\'eductibilit\'e et d'admissibilit\'e. La r\'eponse est affirmative et fait l'objet du pr\'esent
article. Comme dans \cite{Co2}, on d\'emontre donc un isomorphisme
Borel-\'equivariant entre $\breuil(V)^*$ et $(\varprojlim_{\psi}
\dfont(V))^{\rm b}$ (th\'eor\`eme \ref{phigammabanach}) o\`u $\dfont(V)$ est le
$(\varphi,\Gamma)$-module associ\'e \`a la repr\'esentation potentiellement cristalline $V$ et o\`u la limite projective consiste en les suites
$\psi$-compatibles born\'ees
d'\'el\'e\-ments de $\dfont(V)$. Pour montrer cet isomorphisme, il est
n\'ecessaire d'une part d'\'etendre au cas potentiellement cristallin consid\'er\'e la th\'eorie des modules de Wach de \cite{Be1} et \cite{W96}, d'autre part de passer par une description interm\'ediaire de
$\breuil(V)$ comme espace de
fonctions continues sur $\Q$ d'un certain type (th\'eor\`eme
\ref{complete}). Les r\'esultats c\^ot\'e $(\varphi,\Gamma)$-modules
permettent alors de d\'eduire le r\'esultat principal de cet article
(\S\ref{resul}) :

\begin{enonce*}{Th\'eor\`eme}
Si $V$ est une repr\'esentation $p$-adique absolument irr\'e\-ductible de dimension $2$ de $\gal$, qui devient cristalline sur une extension ab\'elienne de $\Q$ et qui est $\varphi$-semi-simple, alors $\breuil(V)$ est non nul,
topologiquement irr\'eductible et admissible.
\end{enonce*}

On obtient aussi deux autres corollaires, l'un concernant tous les
r\'eseaux possibles stables par $\G$ dans ${\rm Alg}(V)\otimes {\rm
Lisse}(V)$ (corollaire \ref{reseaurigolo}), l'autre concernant les
vecteurs localement analytiques dans $\breuil(V)$ (corollaires \ref{anal} et \ref{anal2}).

D'autres isomorphismes $\breuil(V)^*\simeq (\varprojlim_{\psi}\dfont(V))^{\rm b}$ sont d\'emontr\'es lorsque $V$ est trianguline non de Rham dans \cite{Co3}, mais les m\'ethodes d'analyse $p$-adique c\^ot\'e $\G$ y sont sensiblement diff\'erentes. Ici, on utilise de mani\`ere essentielle l'existence d'un entrelacement entre deux fa\c cons d'\'ecrire la s\'erie principale ${\rm
Lisse}(V)$ (correspondant essentiellement aux deux fa\c cons d'ordonner les caract\`eres que l'on induit), entrelacement qui {\og passe \`a la compl\'etion $p$-adique \fg} et permet de d\'efinir un {\it entrelacement $p$-adique} entre deux fa\c cons d'\'ecrire le Banach $\breuil(V)$. On s'aper\c coit alors que cet entrelacement $p$-adique a une interpr\'etation en th\'eorie de Hodge $p$-adique : si deux distributions de $\breuil(V)^*$ se correspondent par cet entrelacement (on ne distingue pas ici les deux mani\`eres d'\'ecrire $\breuil(V)$), alors les deux \'el\'ements qu'on leur associe dans $(\varprojlim_{\psi}\dfont(V))^{\rm b}$ sont reli\'es par une condition \'equivalente \`a la donn\'ee de la {\it filtration de Hodge} sur $\dcris(V)$ (voir \S\ref{deuxlemmes} et lemme \ref{versfin}).

Lorsque $V$ n'est pas $\varphi$-semi-simple, signalons que $\breuil(V)^*$ et $(\varprojlim_{\psi}\dfont(V))^{\rm b}$ sont toujours d\'efinis mais que l'on ignore s'ils sont naturellement isomorphes (l'entrelacement ci-dessus est dans ce cas l'identit\'e).

Une premi\`ere version de cet article (octobre 2004), dont une version pr\'eliminaire avait fait l'objet d'un cours au C.M.S. de Hangzhou en ao\^ut 2004 (\cite{BB}), ne traitait que le cas des repr\'esenta\-tions cristallines.

\subsection{Notations}\label{nota}

On fixe $\Qpbar$ une cl\^oture alg\'ebrique de $\Qp$, on note {\og $\val$ \fg} la valuation sur $\Qpbar$ telle que $\val(p) \= 1$, $|\cdot |$ la norme $p$-adique $|x| \= p^{-\val(x)}$ et $\Cp$ le compl\'et\'e de $\Qpbar$ pour $|\cdot|$. On normalise l'isomorphisme de la th\'eorie du corps de classes local en envoyant les uniformisantes sur les Frobenius g\'eom\'etriques. On note $\eps$ le caract\`ere cyclotomique $p$-adique de $\gal$ vu aussi comme caract\`ere de $\Qp^{\times}$. En particulier, $\eps(p)=1$ et $\eps_{\mid \Zp^{\times}} : \Zp^{\times} \to \Zp^{\times}$
est l'identit\'e. On note ${\rm nr}(x)$ le caract\`ere non ramifi\'e
de $\Qp^{\times}$ envoyant $p$ sur $x$. On d\'esigne par $L$ une
extension finie de $\Qp$, $\OO_L$ son anneau d'entiers, $k_L$ son corps
r\'esiduel et $\pi_L$ une uniformisante. On note $V$ une repr\'esen\-tation $p$-adique de $\gal$, c'est-\`a-dire un $L$-espace vectoriel de dimension finie muni d'une action lin\'eaire et continue de $\gal$, et $T$ un $\OO_L$-r\'eseau de $V$ stable par $\gal$. Enfin, $\B$ d\'esigne les matrices triangulaires sup\'erieures dans $\G$ et $\BK$ le sous-groupe de ces matrices qui sont dans $\K$. 

On note $F_n$ l'extension finie de $\Qp$ dans $\Qpbar$ engendr\'ee par les racines $p^n$-i\`emes de l'unit\'e et $F_\infty \= \cup_{n \geq 0} F_n$. On fixe dans tout cet article le choix d'une suite compatible $(\zeta_{p^n})_{n \geq 0}$ de racines primitives $p^n$-i\`emes de l'unit\'e. Le groupe de Galois $\Gamma \= \mathrm{Gal}(F_\infty / \Qp)$ est isomorphe \`a $\Zp^\times$ via le caract\`ere cyclotomique et si $n \geq 1$, alors $\Gamma_n \= \eps^{-1}(1+p^n \Zp)$ s'identifie au groupe de Galois $\mathrm{Gal}(F_\infty / F_n)$. On note $L_n = L \otimes_{\Qp} F_n$, ce qui fait que $L_n$ est un produit de corps et un $L_n[\Gamma]$-module simple. Si $\eta : \Gamma \to \OO_L^\times$ est un caract\`ere fini \`a valeurs dans $L$, on note $G(\eta)$ la somme de Gauss associ\'ee \`a $\eta$ : si $\eta = 1$, alors $G(\eta)=1$ et si $\eta$ est de conducteur $n = n(\eta) \geq 1$, alors $G(\eta) \= \sum_{\gamma \in \Gamma / \Gamma_n} \eta^{-1}(\gamma) \otimes \gamma(\zeta_{p^n}) \in L_n$. On v\'erifie facilement les propri\'et\'es suivantes des sommes de Gauss : 
\begin{itemize}
\item[(i)] si $g \in \gal$, alors $g(G(\eta)) = \eta(g) G(\eta)$ ($g$ agissant lin\'eairement sur $L$);
\item[(ii)] on a $G(\eta) \cdot G(\eta^{-1}) = p^{n(\eta)} \eta(-1)$ et en particulier $G(\eta) \in L_n^\times$.
\end{itemize}

Tous les espaces de Banach $B$ de ce texte sont $p$-adiques et tels que $\|B\|\subseteq |L|$. On appelle $\G$-Banach unitaire un espace de Banach $B$ muni d'une action \`a gauche $L$-lin\'eaire de $\G$ telle que les applications $\G \to B$, $g\mapsto gv$ sont continues pour tout $v\in B$ et telle que, pour un choix de norme $\|\cdot\|$ sur $B$, on a $\|gv\|=\|v\|$ pour tout $g\in \G$ et tout $v\in B$. Un $\G$-Banach unitaire est dit admissible (suivant \cite{ST3}) si le Banach dual est de type fini sur $L\otimes_{\O}\O[[\K]]$ o\`u $\O[[\K]]\=\varprojlim\O[\K/H]$, la limite projective \'etant prise sur les sous-groupes de congruences principaux $H$ de $\K$. 

\section{Repr\'esentations $p$-adiques}

\subsection{Quelques anneaux de s\'eries formelles} \label{series}

Le but de ce paragraphe est d'introduire certains anneaux de
s\'eries formelles ($\calO$, $\calE$, $\calE^+$, $\calR^+$, $\calE^\dagger$ et $\calR$), ainsi que
certaines des structures dont ils sont munis et dont nous avons besoin
dans la suite de cet article. On commence par rappeler succintement les d\'efinitions de ces divers anneaux.

\begin{itemize}
\item[(i)] On note $\calO$ l'anneau form\'e des
s\'eries $\sum_{i \in \ZZ} a_i X^i$ telles que $a_i \in \OO_L$ et
$a_{-i} \to 0$ quand $i \to +\infty$. C'est un anneau
local de corps r\'esiduel $k_L((X))$.

\item[(ii)] On note $\calE \= \calO[1/p]$, c'est un
corps local de dimension $2$.

\item[(iii)] On note $\calE^+ \= L \otimes_{\OO_L} \OO_L[[X]]$.

\item[(iv)] On note $\calR^+$ l'anneau des 
s\'eries formelles $f(X) \in
L[[X]]$ qui convergent sur le disque unit\'e, ce qui fait que
$\calE^+$ s'identifie au sous-anneau 
de $\calR^+$ constitu\'e des
s\'eries formelles \`a coefficients born\'es.

\item[(v)] On note $\calE^\dagger$ le sous-corps de $\calE$ constitu\'e des 
s\'eries $f(X) = \sum_{i \in \ZZ} a_i X^i$ telles qu'il existe $\rho < 1$ pour lequel $|a_{-i}| \rho^{-i} \to 0$ quand $i \to +\infty$.

\item[(vi)] On note $\calR$ l'anneau form\'e des
s\'eries $\sum_{i \in \ZZ} a_i X^i$ telles que $a_i \in L$, telles qu'il existe $\rho < 1$ pour lequel $|a_{-i}| \rho^{-i} \to 0$ quand $i \to +\infty$, et telles que pour tout $\sigma < 1$, on ait $|a_i| \sigma^i \to 0$ quand $i \to +\infty$. Le corps $\calE^\dagger$ s'identifie alors au sous-anneau de $\calR$ constitu\'e des s\'eries formelles \`a coefficients born\'es.
\end{itemize}

Le corps $\calE$ est muni de la norme de Gauss $\|\cdot\|_{\rm Gauss}$
d\'efinie par $\|f(X)\|_{\rm Gauss} \= \sup_{i \in \ZZ} |a_i|$ si $f(X)=\sum_{i  \in \ZZ} a_i X^i$. L'anneau des entiers de $\calE$ pour cette
norme est $\calO$, et la norme de Gauss induit sur $\calO$ la
topologie $\pi_L$-adique. L'application naturelle $\calO \to
k_L((X))$ est alors continue si l'on donne \`a $\calO$ la topologie
$\pi_L$-adique et \`a $k_L((X))$ la topologie discr\`ete. 

On peut d\'efinir une topologie moins fine sur $\calO$, la topologie
faible, d\'efinie par le fait que les $\{\pi_L^i \calO + X^j
\OO_L[[X]]\}_{i,j \geq 0}$ forment une base de voisinages
de z\'ero, et la topologie faible sur $\calE = \cup_{k \geq 0}
\pi_L^{-k} \calO$ qui est la topologie de la limite inductive. Cette
topologie induit la topologie $(\pi_L,X)$-adique sur $\OO_L[[X]]$, et 
l'application naturelle $\calO \to
k_L((X))$ est alors continue si l'on donne \`a $\calO$ la topologie
faible et \`a $k_L((X))$ la topologie $X$-adique. 

Si $\rho < 1$, alors on peut d\'efinir une norme $\|\cdot\|_{D(0,\rho)}$
sur $\calR^+$ par la formule :
\begin{equation}\label{normerho}
\|f(X)\|_{D(0,\rho)} \= \sup_{\substack{z \in \C_p \\ |z| \leq \rho}}
|f(z)|  =  \sup_{i \geq 0} |a_i| \rho^i,  
\end{equation}
si $f(X)=\sum_{i \geq 0} a_i X^i$. L'ensemble des normes $\{
\|\cdot\|_{D(0,\rho)} \}_{0 \leq \rho < 1}$ d\'efinit une topologie sur
$\calR^+$ qui en fait un espace de Fr\'echet.

\begin{defi}\label{dordrer}
Si $f(X) \in \calR^+$ et $r \in \RR_{\geq 0}$, on dit que $f(X)$ est d'ordre $r$ (il serait plus correct de dire {\og d'ordre $\leq r$ \fg}) si pour un $\rho$ tel que $0 < \rho < 1$, la suite $\{ p^{-nr} \|f(X)\|_{D(0,\rho^{1/p^n})} \}_{n \geq 0}$ est
born\'ee.
\end{defi}

Il est facile de voir que si c'est vrai pour un choix de $0
< \rho < 1$, alors c'est vrai pour tout choix de $0 < \rho < 1$. 
Un exemple de s\'erie d'ordre $1$ est donn\'e par $f(X)=\log(1+X)$.

Nous allons maintenant rappeler les formules qui d\'efinissent l'action de $\Gamma$, le frobenius $\varphi$ et l'op\'erateur $\psi$ sur ces anneaux. Soit $R$ l'un des anneaux $\calO$, $\calE$, $\calE^+$, $\calR^+$, $\calE^\dagger$ ou $\calR$.

Commen\c{c}ons par l'action de $\Gamma$; le caract\`ere cyclotomique $\varepsilon : \Gamma \to \Zp^\times$ est un isomorphisme. L'anneau $R$ est alors muni d'une action de $\Gamma$, telle que si $\gamma \in \Gamma$, alors $\gamma$ agit par un morphisme de $L$-alg\`ebres, et $\gamma(X) \= (1+X)^{\varepsilon(\gamma)}-1$. On v\'erifie facilement que $\Gamma$ agit par des isom\'etries, pour toutes les normes et topologies d\'efinies ci-dessus. Les id\'eaux de $\calR^+$ ou de $\calE^+$ qui sont stables sous l'action de $\Gamma$ sont d'une forme tr\`es particuli\`ere. Pour $n \geq 1$, on note : \[ Q_n(X) = \frac{(1+X)^{p^n}-1 }{(1+X)^{p^{n-1}}-1} \] ce qui fait par exemple que $Q_1(X) = ((1+X)^p-1)/ X$ et que $Q_n(X) = \varphi^{n-1} (Q_1(X))$. Le polyn\^ome $Q_n(X)$ est le polyn\^ome minimal sur $\Qp$ de $\zeta_{p^n}-1$ et l'id\'eal de $L[X]$ qu'il engendre est donc stable sous l'action de $\Gamma$.

\begin{lemm}\label{gamstable}
Si $I$ est un id\'eal principal de $\calR^+$, qui est stable par $\Gamma$, alors il existe une suite d'entiers $\{j_n\}_{n \geq 0}$ telle que $I$ est engendr\'e par un \'el\'ement de la forme $X^{j_0}\prod_{n=1}^{+\infty}(Q_n(X)/p)^{j_n}$.

Si $I$ est un id\'eal de $\calE^+$, qui est stable par $\Gamma$, alors il existe une suite finie d'entiers $j_0,j_1,\hdots,j_m$ telle que $I$ est engendr\'e par un \'el\'ement de la forme $X^{j_0}\prod_{n=1}^m Q_n(X)^{j_n}$.
\end{lemm}

Pour $L=\Qp$ ce lemme fait l'objet de \cite[lemme I.3.2]{Be1} et la d\'emonstration s'adapte sans probl\`eme. Signalons tout de m\^eme que les polyn\^omes $Q_n(X)$ ne sont plus n\'ecessairement irr\'eductibles dans $L[X]$, et que leurs diviseurs \'eventuels engendrent des id\'eaux de $L[X]$ stables par des sous-groupes ouverts de $\Gamma$.

L'anneau $R$ est aussi muni d'un morphisme de Frobenius $\varphi$, qui est
lui-aussi un morphisme de $L$-alg\`ebres, tel que $\varphi(X) \=
(1+X)^p-1$. Cette application est continue pour toutes les topologies
ci-dessus et commute \`a l'action de $\Gamma$. 

Passons maintenant \`a l'op\'erateur $\psi$. L'anneau $R$ est un $\varphi(R)$-module libre de rang $p$, dont une base est donn\'ee par $\{(1+X)^i\}_{0 \leq i \leq p-1}$. Si $y \in R$, on peut donc \'ecrire $y=\sum_{i=0}^{p-1} (1+X)^i \varphi(y_i)$.

\begin{defi}\label{defpsianno}
Si $R$ est l'un des anneaux $\calO$, $\calE$, $\calE^+$, $\calR^+$, $\calE^\dagger$ ou $\calR$, alors on d\'efinit un op\'erateur $\psi : R \to R$ par la formule $\psi(y) \= y_0$ si $y=\sum_{i=0}^{p-1} (1+X)^i \varphi(y_i)$.  
\end{defi}

Cet op\'erateur v\'erifie alors $\psi(\varphi(x)y) = x \psi(y)$ et commute \`a l'action de $\Gamma$. Il ne commute pas \`a $\varphi$ et n'est pas $\OO_L[[X]]$-lin\'eaire. 

\subsection{Repr\'esentations $p$-adiques et $(\varphi,\Gamma)$-modules}  \label{pgmod}

Rappelons qu'une repr\'esen\-tation $L$-lin\'eaire de $\gal$ est un $L$-espace
vectoriel $V$ de dimension finie muni d'une action lin\'eaire et continue de $\gal$. Il est assez difficile de d\'ecrire un tel objet, et nous allons voir
dans ce paragraphe qu'une certaine classe de {\og
$(\varphi,\Gamma)$-modules \fg} permet d'en donner une
description explicite.

Un $\varphi$-module sur $\calO$ est un $\calO$-module de type fini
$\dfont$ muni d'un morphisme $\varphi$-semi-lin\'eaire $\varphi:\dfont
\to \dfont$. On \'ecrit $\varphi^*(\dfont)$ pour le $\calO$-module
engendr\'e par $\varphi(\dfont)$ dans $\dfont$ et on dit que $\dfont$ est \'etale si $\dfont=\varphi^*(\dfont)$. 
Un $\varphi$-module sur $\calE$ est un $\calE$-espace vectoriel de
dimension finie $\dfont$ muni d'un morphsime $\varphi$-semi-lin\'eaire
$\varphi:\dfont \to \dfont$. On dit que $\dfont$ est \'etale si $\dfont$ a un $\calO$-r\'eseau stable par $\varphi$ et \'etale.
Un $(\varphi,\Gamma)$-module est un $\varphi$-module muni d'une
action continue de $\Gamma$ par des morphismes semi-lin\'eaires 
(par rapport \`a l'action de $\Gamma$ sur les coefficients) 
et commutant \`a $\varphi$.

Si $\dfont$ est un $(\varphi,\Gamma)$-module \'etale sur $\calE$, et si
$\bfont = \widehat{\calE^{\rm nr}}$ est l'anneau construit par Fontaine dans
\cite[\S A1.2]{F90}, alors $V(\dfont) \= (\bfont \otimes_{\calE}
\dfont)^{\varphi=1}$ est un $L$-espace vectoriel de dimension
$\dim_{\calE}(\dfont)$ muni d'une action lin\'eaire et continue de $\gal$:
c'est une repr\'esen\-tation $L$-lin\'eaire de $\gal$. On a alors le
r\'esultat suivant (cf. \cite[\S A3.4]{F90}) :

\begin{theo}\label{fontequiv}
Le foncteur $\dfont \mapsto V(\dfont)$ est une \'equivalence de cat\'egories :
\begin{itemize}
\item[(i)] de la cat\'egorie des $(\varphi,\Gamma)$-modules \'etales sur
$\calE$ vers la cat\'egorie des repr\'esen\-tations $L$-lin\'eaires de
$\gal$;
\item[(ii)] de la cat\'egorie des $(\varphi,\Gamma)$-modules \'etales sur
$\calO$ vers la cat\'egorie des $\OO_L$-repr\'esenta\-tions de
$\gal$.
\end{itemize}
\end{theo}

L'inverse de ce foncteur est not\'e $V \mapsto \dfont(V)$, et on peut montrer que $\dfont(V) = (\bfont \otimes_{\Qp} V)^{\galinf}$. Les
$(\varphi,\Gamma)$-modules \'etales nous donnent donc un moyen
commode de travailler avec les repr\'esen\-tations $L$-lin\'eaires. En
th\'eorie, les $(\varphi,\Gamma)$-modules permettent de retrouver
tous les objets que l'on peut associer aux repr\'esen\-tations
$L$-lin\'eaires (et en pratique, c'est souvent le cas). Dans ce
paragraphe et les suivants, nous allons rappeler quelques unes de ces
constructions. 

Si $R$ est l'un des anneaux $\calE^\dagger$ ou $\calR$, alors on d\'efinit de la m\^eme mani\`ere la notion de $\varphi$-module et de $(\varphi,\Gamma)$-module sur $R$. 

Si $\ddag$ est un $\varphi$-module sur $\calE^\dagger$, alors $\dfont \= \calE \otimes_{\calE^\dagger} \ddag$ est un $\varphi$-module sur $\calE$ et on dit que $\ddag$ est \'etale si $\dfont$ l'est. Le r\'esultat ci-dessous (dont l'analogue pour des $\varphi$-modules sans action de $\Gamma$ est faux) est d\^u \`a Cherbonnier et Colmez (cf. \cite[\S III.5]{CC98}).

\begin{theo}\label{ccres}
Le foncteur $\ddag \mapsto \calE \otimes_{\calE^\dagger} \ddag$, de la cat\'egorie des $(\varphi,\Gamma)$-modules \'etales sur $\calE^\dagger$ vers la cat\'egorie des $(\varphi,\Gamma)$-modules \'etales sur $\calE$, est une \'equivalence de cat\'egories.
\end{theo}

En particulier, si $V$ est une repr\'esen\-tation $p$-adique de $\gal$, alors on peut lui associer, gr\^ace au th\'eor\`eme ci-dessus, un $(\varphi,\Gamma)$-module sur $\calE^\dagger$ not\'e $\ddag(V)$. Il existe d'ailleurs un anneau $\bdag \subset \bfont$ tel que $\ddag(V) = (\bdag \otimes_{\Qp} V)^{\galinf}$.

Si $\ddag$ est un $\varphi$-module sur $\calE^\dagger$, alors $\drig \= \calR \otimes_{\calE^\dagger} \ddag$ est un $\varphi$-module sur $\calR$ et $\ddag$ est \'etale si et seulement si $\drig$ est {\og pur de pente nulle \fg} au sens de \cite{KK04}. Nous donnons pour r\'ef\'erence le r\'esultat ci-dessous (cf. \cite{KK04}) qui ne nous servira pas dans le reste de cet article :

\begin{theo}\label{kkres}
Le foncteur $\ddag \mapsto \calR \otimes_{\calE^\dagger} \ddag$, de la cat\'egorie des $\varphi$-modules \'etales sur $\calE^\dagger$ vers la cat\'egorie des $\varphi$-modules de pente nulle sur $\calR$, est une \'equivalence de cat\'egories.
\end{theo}

Pour terminer, donnons d\`es \`a pr\'esent la d\'efinition ci-dessous, qui joue un r\^ole important pour les repr\'esen\-tations que nous consid\'erons dans cet article.

\begin{defi}\label{defhf}
Si $V$ est une repr\'esen\-tation $L$-lin\'eaire de
$\gal$, on dit que $V$ est de hauteur finie s'il existe un sous-$\calE^+$-module $\dplus$ de $\dfont$ qui est stable par $\varphi$ et $\Gamma$ et tel que l'application naturelle $\calE \otimes_{\calE^+} \dplus \to \dfont$ est un isomorphisme. 
\end{defi}

On peut montrer (cf. \cite[\S B1.8]{F90}) que si $\bplus =A_S^+[1/p]$ est l'anneau construit dans \cite[\S B1.8]{F90}, et que l'on pose $\dplus(V) = (\bplus \otimes_{\Qp} V)^{\galinf}$, alors tout $\calE^+$-module $\dplus$ satisfaisant les conditions de la d\'efinition ci-dessus est contenu dans $\dplus(V)$ et donc que $V$ est de hauteur finie si et seulement si $\dfont(V)$ a une base constitu\'ee d'\'el\'ements de $\dplus(V)$. Les repr\'esen\-tations de hauteur finie sont \'etudi\'ees dans \cite[\S B]{F90} et dans \cite{Co4}.

Passons maintenant \`a la d\'efinition de l'op\'erateur $\psi$ sur les $(\varphi,\Gamma)$-modules. Si $\dfont$ est un $\varphi$-module \'etale (sur $\calE^\dagger$ ou sur $\calE$) et si $y \in \dfont$, alors on peut \'ecrire $y=\sum_{i=0}^{p-1} (1+X)^i \varphi(y_i)$ o\`u les $y_i \in \dfont$ sont bien d\'etermin\'es.

\begin{defi}\label{defipsi}
Si $\dfont$ est un $\varphi$-module \'etale sur $\calE^\dagger$ ou sur $\calE$, alors on d\'efinit un op\'erateur $\psi : \dfont \to
\dfont$ par la formule $\psi(y) \= y_0$ si $y=\sum_{i=0}^{p-1} (1+X)^i
\varphi(y_i)$.
\end{defi}

Cet op\'erateur v\'erifie alors $\psi(\varphi(x)y) = x \psi(y)$ et
$\psi(x \varphi(y)) = \psi(x) y$ si $x \in \calE^\dagger$ (ou $\calE$) et $y \in \dfont$, et commute \`a l'action de $\Gamma$ si $\dfont$ est un $(\varphi,\Gamma)$-module. C'est cet op\'erateur qui permet de faire le lien entre les repr\'esen\-tations $L$-lin\'eaires de $\gal$ et les repr\'esen\-tations de $\G$ (ce qui est l'objet de cet article), et aussi le lien entre les $(\varphi,\Gamma)$-modules et la cohomologie d'Iwasawa des repr\'esen\-tations de $\gal$, ce que nous rappelons ci-dessous.

Si $V$ est une repr\'esen\-tation $L$-lin\'eaire, alors les invariants de $\dfont(V)$ sous l'action de $\psi$ jouent un r\^ole important. Nous commen\c{c}ons par rappeler ci-dessous deux propri\'et\'es importantes de ce module.

\begin{prop}\label{tadvpeo}
On a $\dfont(V)^{\psi=1} = \ddag(V)^{\psi=1}$ et le $\calE$-espace vectoriel $\dfont(V)$ a une base constitu\'ee d'\'el\'ements de $\dfont(V)^{\psi=1}$.
\end{prop}

\begin{proof}
Le fait que $\dfont(V)^{\psi=1} = \ddag(V)^{\psi=1}$ est un r\'esultat de Cherbonnier dont on trouvera une d\'emonstration dans \cite[\S III.3]{CC99} et la deuxi\`eme assertion est d\'emontr\'ee dans \cite[\S I.7]{CC99}.
\end{proof}

\begin{coro}\label{dpsinozero}
Si $V$ est une repr\'esen\-tation $L$-lin\'eaire de $\gal$, alors $\dfont(V)^{\psi=1} \neq 0$.
\end{coro}

Rappelons bri\`evement la d\'efinition de la cohomologie d'Iwasawa. Si $V$ est une repr\'esen\-tation $L$-lin\'eaire de $\gal$, si $T$ est un $\OO_L$-r\'eseau $\g$-stable de $V$, et si $i \geq 0$, alors on d\'efinit :
\[ H^i_{\rm Iw}(\Qp,T) \= \varprojlim_{{\rm cor}_{F_{n+1}/F_n}} H^i({\rm  Gal}(\Qpbar/F_n),T). \]  

Le groupe $H^i_{\rm Iw}(\Qp,V) \= \Qp \otimes_{\Zp} H^i_{\rm Iw}(\Qp,T)$
ne d\'epend alors pas du choix de $T$. Les $\Qp \otimes_{\Zp} \Zp[[\Gamma]]$-modules $H^i_{\rm Iw}(\Qp,V)$ ont \'et\'e \'etudi\'es en d\'etail par Perrin-Riou, qui a notamment montr\'e (cf. \cite[\S 3.2]{BP94}) :

\begin{prop}\label{bpriw}
Si $V$ est une repr\'esen\-tation $p$-adique de $\gal$, alors on a $H^i_{\rm
Iw}(\Qp,V)=0$ si $i \neq 1,2$. De plus : 
\begin{itemize}
\item[(i)] le sous-module de $\Qp \otimes_{\Zp} \Zp[[\Gamma]]$-torsion de $H^1_{\rm Iw}(\Qp,V)$ est isomorphe \`a $V^{{\rm Gal}(\Qpbar / F_\infty)}$ et $H^1_{\rm Iw}(\Qp,V)/\{{\rm torsion}\}$ est un $\Qp \otimes_{\Zp} \Zp[[\Gamma]]$-module libre de rang $\dim_{\Qp}(V)$;
\item[(ii)] $H^2_{\rm Iw}(\Qp,V)=(V^*(1)^{{\rm Gal}(\Qpbar /
  F_\infty)})^*$. 
\end{itemize}
\end{prop}

Les espaces $\dfont(V)^{\psi=1}$ et $\dfont(V)/(1-\psi)$ sont eux aussi des $\Qp \otimes_{\Zp} \Zp[[\Gamma]]$-modules et on a le r\'esultat
suivant (cf. \cite[\S II.1]{CC99}) : 

\begin{prop}\label{cciw}
On a des isomorphismes canoniques de $\Qp \otimes_{\Zp} \Zp[[\Gamma]]$-modules $H^1_{\rm Iw}(\Qp,V) \simeq \dfont(V)^{\psi=1}$ et $H^2_{\rm Iw}(\Qp,V) \simeq \dfont(V)/(1-\psi)$. 
\end{prop}

Terminons ce paragraphe avec un lemme technique concernant l'action de $\psi$ sur les $(\varphi,\Gamma)$-modules.

\begin{lemm}\label{noheart}
Si $T$ est une $\OO_L$-repr\'esen\-tation sans torsion ou une $k_L$-repr\'esen\-tation et si $M \subset \dfont(T)$ est un $\OO_L[[X]]$-module de type fini et stable par $\psi$ tel que $\cap_{j \geq 0} \psi^j(M) = 0$, alors $M=0$.
\end{lemm}

\begin{proof}
Montrons tout d'abord le lemme pour les $k_L$-repr\'esen\-tations. Le module $M$ n'est pas n\'ecessairement stable par $\varphi$, mais il existe $r \geq 0$ tel que $\varphi(M) \subset X^{-r} M$ et alors $\varphi(X^r M) \subset X^r M$ ce qui fait que $X^r M \subset \psi(X^r M)$ et donc que :
\[ X^r M \subset \cap_{j \geq 0} \psi^j(X^r M) \subset \cap_{j \geq 0} \psi^j(M) = 0, \] 
ce qui fait que $X^r M=0$ et donc que $M=0$.

Passons \`a pr\'esent au cas des $\OO_L$-repr\'esen\-tations. Dans ce cas, en consid\'erant l'image de $M$ dans $\dfont(T/\pi_L T)$, et en utilisant le lemme pour la $k_L$-repr\'esen\-tation $T/\pi_L T$, on voit qu'un $M$ qui satisfait les hypoth\`eses du lemme est inclus dans $\pi_L \dfont(T)$ et en it\'erant ce proc\'ed\'e, on trouve que $M \subset \cap_{k \geq 0} \pi_L^k \dfont(T) = 0$. 
\end{proof}

\subsection{Topologie faible et treillis}\label{topotr}

Nous allons maintenant nous int\'eresser \`a la topologie des $(\varphi,\Gamma)$-modules sur $\calO$ et $\calE$. Si $\dfont$ est un
$\calO$-module libre de rang $d$, alors le choix d'une base de $\dfont$ donne un isomorphisme $\dfont \simeq \calO^d$ et on peut munir $\dfont$ de la topologie faible induite par cet isomorphisme. Un petit calcul montre qu'une
application $\calO$-lin\'eaire $\calO^d \to \calO^d$ est
n\'ecessairement continue et donc que la topologie d\'efinie sur $\dfont$
par $\dfont \simeq \calO^d$ ne d\'epend pas du choix d'une base de $\dfont$. 

\begin{lemm}\label{bornetreillis}
Si $P$ est une partie d'un $\calO$-module libre $\dfont$, 
et $M(P)$ est le $\OO_L[[X]]$-module engendr\'e par $P$, alors 
les propri\'et\'es suivantes sont \'equivalentes:
\begin{itemize}
\item[(i)] $P$ est born\'ee pour la topologie faible;
\item[(ii)] $M(P)$ est born\'e pour la topologie faible;
\item[(iii)] pour tout $j \geq 1$, l'image de $M(P)$ dans 
$\dfont/\pi_L^j \dfont$ est un $\OO_L[[X]]$-module de type fini. 
\end{itemize}
\end{lemm}

\begin{proof}
Choisissons une base de $\dfont$, et notons $\dplus$ le $\OO_L[[X]]$-module engendr\'e par cette base. Ainsi, la topologie faible sur $\dfont$ est d\'efinie par le fait que les $\{\pi_L^i \dfont + X^j \dplus \}_{i,j \geq 0}$ forment une base de voisinages de z\'ero. En particulier, une partie
$P$ de $\dfont$ est born\'ee si et seulement si pour tout $k \geq 0$, il
existe $n(k,P) \in \ZZ$ 
tel que $P \subset \pi_L^k \dfont + X^{n(k,P)} \dplus$. Comme le
$\OO_L[[X]]$-module engendr\'e par $\pi_L^k \dfont + X^{n(k,P)} \dplus$ est $\pi_L^k \dfont + X^{n(k,P)} \dplus$ lui-m\^eme, on voit que les propri\'et\'es (i) et (ii) sont \'equivalentes. Il reste donc \`a montrer que les $\OO_L[[X]]$-modules born\'es sont ceux qui satisfont (iii), c'est-\`a-dire qu'un $\OO_L[[X]]$-module $M$ est born\'e si et seulement si pour tout $j \geq 1$, l'image de $M$ dans $\dfont/\pi_L^j \dfont$ est un $\OO_L[[X]]$-module de type fini. Si $M$ est born\'e, alors par d\'efinition, pour tout $j
\geq 1$, il existe $n(j,M)$ tel que $M \subset \pi_L^j \dfont + X^{n(j,M)}
\dplus$ ce qui fait que l'image de $M$ dans $\dfont/\pi_L^j \dfont$ est contenue
dans celle de $X^{n(j,M)} \dplus$ et est de type fini puisque $\OO_L[[X]]$
est un anneau noetherien. R\'eciproquement, si $M$ satisfait (iii), 
alors pour tout $j \geq 1$ l'image de $M$ dans $\dfont/\pi_L^j \dfont$ est un
$\OO_L[[X]]$-module de type fini et est donc contenue dans $X^{n(j,M)}
\dplus$ pour un $n(j,M) \in \ZZ$, ce qui fait que 
$M \subset \pi_L^j \dfont + X^{n(j,M)} \dplus$ et donc que $M$ est born\'e pour la topologie faible.
\end{proof}

\begin{defi}\label{deftrei}
Un treillis de $\dfont$ est un $\OO_L[[X]]$-module born\'e $M \subset \dfont$
tel que l'image de $M$ dans $\dfont/\pi_L \dfont$ en est un $k_L[[X]]$-r\'eseau. Un treillis d'un $\calE$-module $\dfont$ est un treillis d'un
$\calO$-r\'eseau de $\dfont$.
\end{defi}

Les treillis font l'objet d'une \'etude d\'etaill\'ee dans \cite[\S
4]{Co2}. Rappelons le r\'esultat de base sur les treillis 
stables par $\psi$ d'un $\varphi$-module $\dfont$ (c'est la proposition
4.29 de \cite{Co2}):  

\begin{prop}\label{col429}
Si $\dfont$ est un $\varphi$-module \'etale sur $\calO$, il existe un
unique treillis $\dsharp$ de $\dfont$ v\'erifiant les propri\'et\'es
suivantes :
\begin{itemize}
\item[(i)] quels que soient $x \in \dfont$ et $k \in \NN$, il existe $n(x,k) \in
  \NN$ tel que $\psi^n(x) \in \dsharp + p^k \dfont$ si $n \geq n(x,k)$;
\item[(ii)] l'op\'erateur $\psi$ induit une surjection de $\dsharp$ sur
  lui-m\^eme.
\end{itemize}
De plus :
\begin{itemize}
\item[(iii)] si $N$ est un sous-$\OO_L[[X]]$ module born\'e 
de $\dfont$ et $k \in \NN$, il existe
$n(N,k)$ tel que $\psi^n(N) \subset \dsharp + p^k \dfont$ si $n \geq
n(N,k)$;
\item[(iv)] si $N$ est un treillis de $\dfont$ stable par $\psi$ tel que
  $\psi$ induise une surjection de $N$ sur lui-m\^eme, alors $N
  \subset \dsharp$ et $\dsharp / N$ est annul\'e par $X$. 
\end{itemize}
\end{prop}

Si $V$ est une repr\'esen\-tation $L$-lin\'eaire de $\gal$, on note $\dsharp(V)$ le treillis associ\'e \`a $\dfont(V)$ dans la proposition ci-dessus. Le lemme ci-dessous pr\'ecise le (iv) de cette proposition.

\begin{lemm}\label{sco429}
Si $V$ est irr\'eductible et de dimension $\geq 2$, et si $N$ est un treillis de $\dfont(V)$ stable par $\psi$ et tel que $\psi$ induise une surjection de $N$ sur lui-m\^eme, alors $N = \dsharp(V)$.
\end{lemm}

\begin{proof}
Si $V$ est irr\'eductible et de dimension $\geq 2$, alors $V^{\mathrm{Gal}(\Qpbar/\Qp^{\mathrm{ab}})} = 0$ et le lemme suit alors du (iii) de la remarque 5.5 de \cite{Co2} et de la d\'emonstration de la proposition 4.47 de \cite{Co2}. 
\end{proof}

Le reste de ce paragraphe est consacr\'e au d\'ebut de l'\'etude des limites projectives d\'efinies ci-dessous.

\begin{defi}\label{defilimproj}
\begin{itemize}
\item[(i)] On note $(\projlim_{\psi} \dfont(V))^{\rm b}$ le $L$-espace vectoriel des suites $(x_n)_{n \geq 0}$ d'\'el\'e\-ments de $\dfont(V)$ telles que l'ensemble $\{x_n\}_{n \geq 0}$ est born\'e (d'o\`u le {\og b \fg}) pour 
la topologie faible et telles que $\psi(x_{n+1})=x_n$.
\item[(ii)] On note $(\projlim_{\psi} \dsharp(V))^{\rm b}$ le $L$-espace vectoriel des suites $(x_n)_{n \geq 0}$ d'\'el\'e\-ments de $\dsharp(V)$ telles que l'ensemble $\{x_n\}_{n \geq 0}$ est born\'e pour la topologie faible et telles que $\psi(x_{n+1})=x_n$.
\item[(iii)] Si $T$ est un $\OO_L$-r\'eseau de $V$ stable par $\gal$, on note
$(\projlim_{\psi}\dsharp(T))^{\rm b}$ le $\OO_L$-module des suites
$\psi$-compatibles et born\'ees pour la topologie faible d'\'el\'ements de $\dsharp(T)$.
\end{itemize}
Ces trois espaces sont munis de la topologie de la limite projective.
\end{defi}

\begin{prop}\label{ddiese}
L'injection $\dsharp(V) \hookrightarrow \dfont(V)$
induit un isomorphisme topologique $(\projlim_{\psi}\dsharp(V))^{\rm b} \to (\projlim_{\psi}\dfont(V))^{\rm b}$ et si $T$ est un $\OO_L$-r\'eseau de $V$ stable par $\gal$, alors $L \otimes_{\OO_L} (\projlim_{\psi} \dsharp(T))^{\mathrm{b}} = (\projlim_{\psi}\dsharp(V))^{\rm b}$.
\end{prop}

\begin{proof}
Il est clair que l'application $(\projlim_{\psi}\dsharp(V))^{\rm b} \to (\projlim_{\psi}\dfont(V))^{\rm b}$ est injective. Fixons un $\OO_L$-r\'eseau $T$ de $V$ stable par $\gal$. Si $x=(x_n)_{n \geq 0} \in (\projlim_{\psi}\dfont(V))^{\rm b}$, alors par d\'efinition l'ensemble $P \= \{x_n\}_{n \geq 0}$ est born\'e pour la topologie faible et par le lemme \ref{bornetreillis}, le $\OO_L[[X]]$-module $M(P)$ engendr\'e par $P$ est born\'e (pour la topologie faible). Quitte \`a multiplier $x$ par une
puissance de $\pi_L$, on peut d'ailleurs supposer que $x_m \in \dfont(T)$
pour tout $m \geq 0$. Si $k,m \geq 0$, alors
la proposition \ref{col429} appliqu\'ee \`a $N=M(P)$ montre que 
$x_m=\psi^n(x_{m+n}) \in \dsharp(T) + p^k \dfont(T)$ si $n$ est assez
grand. Comme c'est vrai pour tout $k$, on en d\'eduit que $x_m \in
\dsharp(T)$ pour tout $m$ et donc que l'application 
$L \otimes_{\OO_L} (\projlim_{\psi}\dsharp(T))^{\mathrm{b}} \to   
(\projlim_{\psi}\dfont(V))^{\mathrm{b}}$ est une bijection. C'est un
hom\'eomorphisme car la topologie de $\dsharp(T)$ est la topologie
induite par la topologie faible de $\dfont(T)$ via l'inclusion $\dsharp(T)
\subset \dfont(T)$. On en d\'eduit les deux points de la proposition.
\end{proof}

Rappelons que $\psi:\dsharp(T)\to \dsharp(T)$ est
surjective, ce qui fait que les applications de transition dans 
$\projlim_{\psi}\dsharp(T)$ le sont.
Le lemme ci-dessous et son corollaire
seront utilis\'es dans le paragraphe \S\ref{laouiw}. Nous verrons en effet plus bas que, \emph{pour les repr\'esen\-tations que nous consid\'erons}, on a $(\projlim_{\psi}\dsharp(T))^{\mathrm{b}} = \projlim_{\psi}\dsharp(T)$.

\begin{lemm}\label{psietddiese}
L'application naturelle $(\projlim_{\psi}\dsharp(T)) / (1-\psi)
\to \dsharp(T) / (1-\psi)$ est un isomorphisme.
\end{lemm}

\begin{proof}
Cette application est \'evidemment surjective, et nous allons montrer
qu'elle est injective, c'est-\`a-dire que si  $x=(x_n)_{n \geq 0} \in
\projlim_{\psi}\dsharp(T)$, avec $x_0 \in
(1-\psi)\dsharp(T)$, alors $x \in (1-\psi)
\projlim_{\psi}\dsharp(T)$. Soit $y_0 \in \dsharp(T)$ tel que
$(1-\psi)y_0 = x_0$. Pour tout $m \geq 0$, il existe $y^0_m \in
\dsharp(T)$ tel que $\psi^m(y_m^0) = y_0$ et on a alors
$(1-\psi)y_m^0 - x_m \in \dsharp(T)^{\psi^m =0}$. L'op\'erateur
$1-\psi$ est 
bijectif sur $\dsharp(T)^{\psi^m =0}$ 
(un inverse \'etant donn\'e par $1+\psi+\psi^2+\cdots+\psi^{m-1}$)
et il existe
donc $y_m \in \dsharp(T)$ tel que $(1-\psi)y_m = x_m$. Pour tout $k
\geq 0$, soit $z_k \= (z_{k,n})_{n \geq 0} \in \projlim_{\psi}
\dsharp(T)$ tel que $z_{k,k} = y_k$.
Comme $\projlim_{\psi}\dsharp(T)$ est
un espace topologique compact, la suite $\{z_k\}_{k \geq 0}$ a une
valeur d'adh\'erence $z$ et comme $(1-\psi)z_k \to x$ quand $k
\to \infty$ (par la continuit\'e de $\psi$, voir la
proposition \ref{actgcont} ci-dessous), 
on voit que $(1-\psi)z = x$ et donc que $x \in (1-\psi)
\projlim_{\psi}\dsharp(T)$.
\end{proof}

\begin{coro}\label{h2ddiese}
On a un isomorphisme de $\OO_L[[\Gamma]]$-modules
$(\projlim_{\psi}\dsharp(T)) / (1-\psi) \simeq H^2_{\rm Iw}(\Q,T)$. 
\end{coro}

\begin{proof}
La proposition 4.43 de \cite{Co2} affirme que $\dsharp(T) / (1-\psi) =
\dfont(T)/(1-\psi)$, et  la remarque II.3.2 de \cite{CC99} (cf. la proposition \ref{cciw} ci-dessus) nous dit que $\dfont(T)/(1-\psi) = H^2_{\rm Iw}(\Qp, T)$ ce qui fait que l'on a une succession d'isomorphismes : 
\[ (\projlim_{\psi}\dsharp(T)) / (1-\psi) \simeq \dsharp(T) /
(1-\psi) \dsharp(T) \simeq  \dfont(T)/(1-\psi)\dfont(T) \simeq  H^2_{\rm
  Iw}(\Qp,T). \]
\end{proof}

\subsection{Th\'eorie de Hodge $p$-adique}\label{cristallin}

Le but de ce paragraphe est de rappeler la classification de certaines repr\'esen\-tations $p$-adiques, les repr\'esen\-tations apc (qui sont d\'efinies ci-dessous), en terme de certains $(\varphi,\gal)$-modules filtr\'es.

Afin de classifier certaines repr\'esen\-tations $L$-lin\'eaires du groupe
$\gal$, Fontaine a introduit (entre autres) les anneaux $\bcris$ et
$\bdR$ (cf. \cite{F3}). Ces anneaux v\'erifient les propri\'et\'es suivantes:
\begin{itemize}
\item[(i)] L'anneau $\bcris$ est une $\Qp$-alg\`ebre munie d'une action
  de $\gal$, telle que $\bcris^{\gal}=\Qp$ et d'un Frobenius $\varphi$ qui
  commute \`a l'action de $\gal$;
\item[(ii)] le corps $\bdR$ est le corps des fractions d'un anneau 
  complet de valuation discr\`ete
  $\bdR^+$ (dont le corps r\'esiduel est
  $\Cp$), et il est donc muni de la filtration d\'efinie par les
  puissances de l'id\'eal maximal. Il est aussi muni d'une action
  continue de $\gal$, telle que $\bdR^{\gal}=\Qp$ et la filtration 
  est stable sous l'action de $\gal$;
\item[(iii)] on a une inclusion naturelle $\bcris \subset \bdR$ et une
  suite exacte (dite {\og fondamentale \fg}) : 
  \[ 0 \to 
  \Qp \to \bcris^{\varphi=1}
  \to \bdR/\bdR^+ \to 0; \]
\item[(iv)] il existe $t \in \bcris$ tel que $\varphi(t)=pt$ et $t$ est un
  g\'en\'erateur de l'id\'eal maximal de $\bdR^+$. 
  Le choix d'un tel $t$, qui est d\'etermin\'e par le choix d'une suite compatible  $(\zeta_{p^n})_{n \geq 0}$ de racines primitives $p^n$-i\`emes de l'unit\'e, d\'etermine une application injective $\calR^+ \to L \otimes_{\Qp} \bcris$, donn\'ee par la formule $f(X) \mapsto f(\exp(t)-1)$, ce qui fait par exemple que $t$ est l'image de $\log(1+X)$, et cette injection commute \`a $\gal$ et \`a $\varphi$; 
\item[(v)] si $m \geq 0$, alors on a une application injective $\iota_m \= \varphi^{-m} : L \otimes_{\Qp} \bcris \to  L \otimes_{\Qp} \bdR$, donn\'ee sur $\calR^+$ par la formule $f(X) \mapsto f(\zeta_{p^m} \exp(t/p^m) - 1)$, et la filtration de $\calR^+$ induite par cette application est donn\'ee par la filtration {\og ordre d'annulation en $\zeta_{p^m}-1$ \fg}.
\end{itemize}

\'Etant donn\'ee une repr\'esen\-tation $L$-lin\'eaire $V$ de $\gal$, on pose :
\begin{equation}\label{defdcrisddr}
\dcris(V) \= (\bcris \otimes_{\Qp} V)^{\gal}
\quad\text{et}\quad\ddR(V) \= (\bdR
\otimes_{\Qp} V)^{\gal}.
\end{equation} 

En g\'en\'eral, ces $\Qp$-espaces vectoriels sont de
dimensions inf\'erieures ou \'egales \`a $\dim_{\Qp}(V)$ et on dit que $V$ est
cristalline (resp. de de Rham) si $\dim_{\Qp} \dcris(V)$
(resp. $\dim_{\Qp} \ddR(V)$) est \'egale \`a $\dim_{\Qp}(V)$. 
Comme $\bcris \subset \bdR$, une repr\'esen\-tation cristalline
est n\'ecessairement de de Rham. 

Le Frobenius $\varphi$ de $\bcris$ commute \`a l'action de
$\gal$ et la filtration de $\bdR$ est stable par $\gal$, 
ce qui fait que $\dcris(V)$
est un $\varphi$-module et que $\ddR(V)$ est un $\Qp$-espace vectoriel
filtr\'e. Si $V$ est une repr\'esen\-tation cristalline, 
alors l'application naturelle de 
$\dcris(V)$ dans $\ddR(V)$ est un isomorphisme 
et $\dcris(V)$ est donc un $\varphi$-module filtr\'e. 
Si $V$ est $L$-lin\'eaire, alors $\dcris(V)$ et $\ddR(V)$ sont
naturellement des $L$-espaces vectoriels, et $\varphi:\dcris(V)
\to \dcris(V)$ ainsi que la filtration sur $\ddR(V)$ sont 
$L$-lin\'eaires.

Si $V$ est une repr\'esen\-tation de de Rham, les poids de
Hodge-Tate de $V$ sont par d\'efinition les entiers $h$ tels que ${\rm Fil}^{-h} \ddR(V) \neq {\rm Fil}^{-h+1} \ddR(V)$, compt\'es avec la
multiplicit\'e $\dim_L {\rm Fil}^{-h}
\ddR(V) / {\rm Fil}^{-h+1} \ddR(V)$.

Si $D$ est un $\varphi$-module de dimension $1$, on d\'efinit
$t_N(D) \= {\rm val}(\alpha)$ o\`u $\alpha$ 
est la matrice de $\varphi$ dans une
base de $D$ et si $D$ est un espace vectoriel filtr\'e de dimension
$1$, on d\'efinit $t_H(D)$ comme \'etant le plus grand $h \in
\ZZ$ tel que ${\rm Fil}^h D \neq 0$. Si $D$ est un $\varphi$-module de
dimension $\geq 1$, on d\'efinit $t_N(D) \= t_N(\det D)$ et $t_H(D) \=
t_H(\det D)$, ce qui fait que $t_H(D)$ est aussi l'oppos\'e de la
somme des poids de Hodge-Tate de $V$, compt\'es avec multiplicit\'es. 

Si $D$ est un $\varphi$-module filtr\'e, on dit que $D$ est
admissible si $t_N(D) = t_H(D)$ et si $t_N(D') - t_H(D') \geq
0$ pour tout sous-$\varphi$-module $D' \subset D$. 
Le fait que pour tout $h \geq 1$, on a ${\rm Fil}^{h+1}
\bcris^{\varphi=p^h} = 0$
permet de montrer (cf. \cite[\S 5.4]{F4}) que si $V$ est une repr\'esen\-tation
cristalline, alors $\dcris(V)$ est admissible. 

On peut aussi d\'efinir la notion de $(\varphi,\gal)$-modules filtr\'es $L$-lin\'eaires admissibles (cf. \cite[\S 4]{F4}), et la notion de repr\'esen\-tation potentiellement cristalline (cf. \cite[\S 5]{F4}); on a alors le th\'eor\`eme de Colmez-Fontaine (voir \cite[th\'eor\`eme A]{CF} et aussi \cite{Be2}) :

\begin{theo}\label{CFthm}
Le foncteur $\dcris(\cdot)$ est une \'equivalence de cat\'egories
de la cat\'ego\-rie des repr\'esen\-tations $L$-lin\'eaires
potentiellement cristallines de $\gal$ dans la cat\'egorie des $(\varphi,\gal)$-modules filtr\'es $L$-lin\'eaires admissibles.
\end{theo}

Passons maintenant \`a la classe de repr\'esen\-tations qui nous int\'eressent.

\begin{defi}\label{defapc}
Si $V$ est une repr\'esen\-tation $p$-adique de $\gal$, alors on dit que $V$ est apc ({\og ab\'eliennement potentiellement cristalline \fg}) s'il existe $n \geq 0$ tel que la restriction de $V$ \`a $\mathrm{Gal}(\Qpbar/F_n)$ est cristalline. 
\end{defi}

Si $V$ est apc, alors on d\'efinit $n(V)$ comme \'etant le plus petit entier $n \geq 1$ tel que la restriction de $V$ \`a $\mathrm{Gal}(\Qpbar/F_n)$ est cristalline. On a donc $n(V)=1$ plut\^ot que $n(V)=0$ si $V$ est cristalline, ceci pour des raisons techniques li\'ees aux modules de Wach (cf. ci-dessous). Remarquons que bien entendu, une repr\'esen\-tation apc est de de Rham. Si $V$ est une repr\'esen\-tation apc, alors on pose $\dcris(V) \= (\bcris \otimes_{\Qp} V)^{\mathrm{Gal}(\Qpbar/F_n)}$ et cette d\'efinition ne d\'epend pas de $n \geq n(V)$; $\dcris(V)$ est alors un $L$-espace vectoriel muni d'un frobenius $L$-lin\'eaire et d'une action de $\Gamma$ qui commute \`a $\varphi$ et qui est triviale sur $\Gamma_{n(V)}$. Si $F_{n(V)} \subset K \subset F_\infty$, alors $K \otimes_{\Qp} \dcris(V) = K \otimes_{\Qp} \ddR(V)$. 

\begin{defi}\label{defmv}
On pose $m(V) \= \inf_\chi n(V(\chi))$ o\`u $\chi$ parcourt l'ensemble des caract\`eres d'ordre fini de $\Gamma$. C'est donc le {\og conducteur essentiel \fg} de $V$. 
\end{defi}

Les repr\'esen\-tations qui nous int\'eressent dans la suite de cet article sont les repr\'esen\-tations $L$-lin\'eaires $V$ de $\gal$ qui sont apc, absolument irr\'eductibles, de dimension $2$ et dont les poids de Hodge-Tate sont $0$ et $k-1$ avec $k \geq 2$. Par le th\'eor\`eme \ref{CFthm}, pour se donner une telle repr\'esen\-tation, il suffit de se donner un $(\varphi,\gal)$-module filtr\'e $L$-lin\'eaire admissible v\'erifiant certaines conditions. 

Soient $\alpha$ et $\beta$ deux caract\`eres localement constants $\alpha, \beta : \Qp^\times \to L^\times$, qui v\'erifient $-(k-1) < \val(\alpha(p)) \leq \val(\beta(p)) < 0$ et $\val(\alpha(p)) + \val(\beta(p)) = -(k-1)$ et qui sont triviaux sur $1+p^n \Zp$ pour un $n \geq 1$. On pose $\alpha_p = \alpha(p)^{-1}$ et $\beta_p = \beta(p)^{-1}$, ce qui fait que $\alpha_p$ et $\beta_p$ appartiennent \`a l'id\'eal maximal de $\OO_L$. 

\begin{defi}\label{defdab}
On note $D(\alpha,\beta)$ le $(\varphi,\gal)$-module filtr\'e d\'efini par $D(\alpha,\beta) = L \cdot e_\alpha \oplus L \cdot e_\beta$ o\`u :

\begin{itemize}
\item[(i)] Si $\alpha \neq \beta$, alors :
\[  \begin{cases}  \varphi(e_\alpha) & = \alpha_p^{-1} e_\alpha \\  
 \varphi(e_\beta) & = \beta_p^{-1} e_\beta  \end{cases} 
 \qquad\text{et si $g \in \Gamma$, alors :}\qquad
 \begin{cases}  g(e_\alpha) & = \alpha(\chi(g)) e_\alpha \\  
 g(e_\beta) & = \beta(\chi(g)) e_\beta  \end{cases} 
 \]
et
\[ \mathrm{Fil}^i (L_n \otimes_L D(\alpha,\beta)) =
 \begin{cases} 
 L_n \otimes_L D(\alpha,\beta) & \text{si $i \leq -(k-1)$;} \\
 L_n \cdot (e_\alpha + G(\beta \alpha^{-1})  \cdot e_\beta) & \text{si $-(k-2) \leq i \leq 0$;} \\
 0 & \text{si $i \geq 1$.} \\
 \end{cases}
  \]
\item[(ii)] Si $\alpha = \beta$, alors :
\[  \begin{cases} \varphi(e_\alpha) & = \alpha_p^{-1} e_\alpha \\  
 \varphi(e_\beta) & =  \beta_p^{-1} (e_\beta - e_\alpha) \end{cases} 
 \qquad\text{et si $g \in \Gamma$, alors :}\qquad
 \begin{cases}  g(e_\alpha) & = \alpha(\chi(g)) e_\alpha \\  
 g(e_\beta) & = \beta(\chi(g)) e_\beta  \end{cases} \]
et
\[ \mathrm{Fil}^i (L_n \otimes_L D(\alpha,\beta)) =
 \begin{cases} 
 L_n \otimes_L D(\alpha,\beta) & \text{si $i \leq -(k-1)$;} \\
 L_n \cdot e_\beta & \text{si $-(k-2) \leq i \leq 0$;} \\
 0 & \text{si $i \geq 1$.} \\
 \end{cases}
  \]
\end{itemize}
\end{defi}

\begin{prop}\label{yatouapc}
Si $V$ est une repr\'esen\-tation $L$-lin\'eaire apc de $\gal$, absolument irr\'eductible, de dimension $2$ et dont les poids de Hodge-Tate sont $0$ et $k-1$ avec $k \geq 2$, alors il existe deux caract\`eres $\alpha$ et $\beta$ comme ci-dessus tels que $\dcris(V) = D(\alpha, \beta)$.

R\'eciproquement, si $\alpha$ et $\beta$ sont deux tels caract\`eres, alors il existe une repr\'esen\-tation $L$-lin\'eaire apc $V$ de $\gal$, absolument irr\'eductible, telle que $\dcris(V) = D(\alpha, \beta)$.
\end{prop}

Voir \cite[\S 5.5]{Co3} pour une d\'emonstration. Terminons ce paragraphe en remarquant que si $V$ d\'enote la repr\'esen\-tation associ\'ee \`a $D(\alpha, \beta)$, alors $n(V) = \sup(n(\alpha),n(\beta),1)$ tandis que $m(V) = \sup (n(\alpha^{-1}\beta),1)$.

\subsection{Th\'eorie de Hodge $p$-adique et $(\varphi,\Gamma)$-modules}\label{thpg}

Dans le paragraphe \S\ref{pgmod}, nous avons rappel\'e quelques points de la th\'eorie des $(\varphi,\Gamma)$-modules et dans le paragraphe \S\ref{cristallin}, nous avons rappel\'e quelques points de la th\'eorie de Hodge $p$-adique. Ces deux th\'eories ne sont pas ind\'ependantes, et dans ce paragraphe nous allons rappeler le r\'esultat essentiel (le th\'eor\`eme \ref{lbcris} ci-dessous) permettant de passer de l'une \`a l'autre. Avant de faire cela, rappelons que dans \cite{Ber1}, on a construit un anneau $\btrig$ qui a la propri\'et\'e que $\bdag \subset \btrig$ et que $\bcris \subset \btrig[1/t]$, ce qui fait que si $V$ est une repr\'esen\-tation $L$-lin\'eaire de $\gal$, alors $\ddag(V) \subset \btrig[1/t] \otimes_{\Qp} V$ et $\dcris(V) \subset  \btrig[1/t] \otimes_{\Qp} V$. L'un des r\'esultats principaux de \cite{Ber1} (cf. \cite[th\'eor\`eme 0.2]{Ber1}) est alors le suivant.

\begin{theo}\label{lbcris}
Si $V$ est une repr\'esen\-tation apc, alors :
\begin{itemize}
\item[(i)] $\ddag(V) \subset  \calR[1/t] \otimes_L \dcris(V)$ et de plus : 
\[ \calR[1/t] \otimes_{\calE^\dagger} \ddag(V) = \calR[1/t] \otimes_L \dcris(V); \]
\item[(ii)] si l'on suppose de plus que les poids de Hodge-Tate de $V$ sont $\geq 0$, alors $\ddag(V) \subset  \calR \otimes_L \dcris(V)$.
\end{itemize}
\end{theo}

Comme on l'a rappel\'e dans la proposition \ref{tadvpeo}, $\ddag(V)$ a une base form\'ee d'\'el\'ements de $\ddag(V)^{\psi=1}$. L'avantage de ces \'el\'ements est que, comme ils sont fix\'es par $\psi$, ils tendent \`a avoir peu de d\'enominateurs, comme le montre la proposition ci-dessous.

\begin{prop}\label{regupsi}
Si $V$ est une repr\'esen\-tation apc, telle que les pentes de $\varphi$ sur $\dcris(V)$ sont $<0$, et si $y \in \calR \otimes_L \dcris(V)$ v\'erifie $\psi(y)=y$, alors $y \in \calR^+ \otimes_L \dcris(V)$.
\end{prop}

\begin{proof} 
On suppose que $L$ contient les valeurs propres de $\varphi$ sur $\dcris(V)$ et fixe une base $e_1,\hdots,e_d$ de $\dcris(V)$ dans laquelle la matrice $(p_{i,j})$ de $\varphi$ est triangulaire sup\'erieure. Si l'on \'ecrit $y = \sum_{i=1}^d y_i \otimes \varphi(e_i)$, alors l'\'equation $\psi(y)=y$ devient : \[ \psi(y_k) = p_{k,k} y_k + \sum_{j > k} p_{k,j} y_j, \] 
pour $k=1,\hdots,d$. Comme on a suppos\'e que les $p_{k,k}$ sont de valuations $<0$, la proposition 1.10 de \cite{Co2} montre que si $\psi(y_k) - p_{k,k} y_k \in \calR^+$, alors $y_k \in \calR^+$, ce qui permet de conclure par r\'ecurrence descendante sur $k$. 
\end{proof}

Le th\'eor\`eme ci-dessous est une g\'en\'eralisation de \cite[th\'eor\`eme 1]{Co4}. 

\begin{theo}\label{apchf}
Si $V$ est une repr\'esen\-tation apc, alors $V$ est de hauteur finie.
\end{theo}

\begin{proof}
Il s'agit de montrer que le $\calE$-espace vectoriel $\dfont(V)$ admet une base dont les \'el\'ements appartiennent \`a $\dplus(V)$. 
Comme les s\'eries $\varphi^n(X)$ sont inversibles dans $\calE$, il suffit de montrer que le $\calE$-espace vectoriel $\dfont(V)$ admet une base dont les \'el\'ements appartiennent \`a $\dplus(V)[\{\varphi^n(1/X)\}_{n \geq 0}]$.
Le corollaire I.7.6 de \cite{CC99} montre que $\dfont(V)$ admet une base dont les \'el\'ements appartiennent \`a $\dfont(V)^{\psi=1}$ et il suffit par suite de voir que si $y \in \dfont(V)^{\psi=1}$, alors $y \in \dplus(V)[\{\varphi^n(1/X)\}_{n \geq 0}]$. Quitte \`a tordre $V$, on suppose que ses poids de Hodge-Tate sont $>0$.

Si $y \in \dfont(V)^{\psi=1}$, alors $y \in \ddag(V)$ par la proposition \ref{tadvpeo} et donc par le th\'eor\`eme \ref{lbcris}, $y \in \calR \otimes_L \dcris(V)$. La proposition \ref{regupsi} nous dit alors que $y \in \calR^+ \otimes_L \dcris(V)$. On en d\'eduit que d'une part $y \in \bdag \otimes_{\Qp} V$ et d'autre part que $y \in \bcontp[1/t] \otimes_{\Qp} V$. Comme $\bdag \cap \bcontp[1/t] = \bplus[\{\varphi^n(1/X)\}_{n \geq 0}]$, on en d\'eduit que $y \in \dplus(V)[\{\varphi^n(1/X)\}_{n \geq 0}]$.
\end{proof}

Si $V$ est une repr\'esen\-tation de hauteur finie, alors l'application $\iota_m = \varphi^{-m} : \bplus \to \bdR^+$ \'evoqu\'ee au paragraphe pr\'ec\'edent nous donne une application $\iota_m : \dplus(V) \to \bdR^+ \otimes_{\Qp} V$ qui se prolonge d'ailleurs en $\iota_m : \calR^+[1/t] \otimes_{\calE^+} \dplus(V) \to \bdR \otimes_{\Qp} V$.

\section{Repr\'esentations apc}\label{wach}

La th\'eorie des modules de Wach, qui est d\'evelopp\'ee pour les repr\'esen\-tations cristallines dans \cite{Be1}, s'\'etend aux repr\'esen\-tations apc et l'objet de cette partie est de donner des d\'emonstrations des r\'esultats que l'on obtient. Remarquons que l'hypoth\`ese {\og apc \fg} est optimale en un certain sens puisque Wach a montr\'e qu'une repr\'esen\-tation de de Rham est de hauteur finie si et seulement si elle est apc (cf. \cite[\S A.5]{W96}). La plupart des r\'esultats de ce chapitre, ainsi que leurs d\'emonstrations, sont similaires \`a ceux de \cite{Be1}.

\subsection{Modules de Wach}\label{constwa}

On dit qu'une repr\'esen\-tation de Hodge-Tate est positive si ses poids de Hodge-Tate sont $\leq 0$ (terminologie un peu malheureuse). Rappelons que l'on a montr\'e au paragraphe \S\ref{thpg} qu'une repr\'esen\-tation apc est de hauteur finie.

\begin{theo}\label{exunmod}
Si $V$ est une repr\'esen\-tation apc positive, alors il existe un unique sous-$\calE^+$-module $\nwach(V)$ de $\dplus(V)$ qui satisfait les conditions suivantes :
\begin{itemize}
\item[(i)] On a $\dfont(V)  = \calE \otimes_{\calE^+} \nwach(V)$;
\item[(ii)] l'action de $\Gamma$ pr\'eserve $\nwach(V)$ et est finie sur $\nwach(V) / X \nwach(V)$;
\item[(iii)] il existe $h \geq 0$ tel que $X^h \dplus(V) \subset \nwach(V)$.
\end{itemize}
Le module $\nwach(V)$ est alors stable par $\varphi$.
\end{theo}

\begin{proof}
Le point de d\'epart de la d\'emonstration est un r\'esultat de Wach (cf. \cite[p. 380]{W96}), qui affirme que si $V$ est de hauteur finie, alors elle est positive et apc si et seulement s'il existe $N  \subset \dplus(V)$ satisfaisant les conditions (i) et (ii) ci-dessus. Il nous reste donc \`a montrer que l'on peut imposer la condition (iii) et qu'alors $N$ est uniquement d\'etermin\'e. 

Soit donc $N$ v\'erifiant (i) et (ii) et soit $I_N$ l'id\'eal de $\calE^+$ constitu\'e de l'ensemble des $\lambda \in \calE^+$ tels que $\lambda \cdot \dplus(V) \subset N$. Cet id\'eal est non nul (par le (i)) et stable par $\Gamma$ (par le (ii)) ce qui fait que, par le lemme \ref{gamstable}, il existe $\alpha_0,\hdots,\alpha_n$ tels que $I_N$ est engendr\'e par $X^{\alpha_0} Q_1(X)^{\alpha_1} \cdots Q_n(X)^{\alpha_n}$. Si l'on pose $\nwach(V) = \dplus(V) \cap N[\{Q_i(X)^{-1}\}_{i \geq 1}]$, alors $\nwach(V)$ est un $\calE^+$-module libre de rang $d$ (puisque $\calE^+$ est principal et que l'on a $N \subset \nwach(V) \subset \dplus(V)$)  qui est toujours stable par $\Gamma$. L'application naturelle $N / X N \to \nwach(V) / X \nwach(V)$ est injective et donc bijective (puisque $N/XN$ et $\nwach(V)/X\nwach(V)$ sont deux $L$-espaces vectoriels de dimension $d$). Enfin, si $y \in \dplus(V)$, alors $X^{\alpha_0} Q_1(X)^{\alpha_1} \cdots Q_n(X)^{\alpha_n} y \in N$ et donc $X^{\alpha_0} y \in \nwach(V)$ ce qui fait que (iii) est v\'erifi\'ee avec $h=\alpha_0$.

Ceci montre l'existence de $\nwach(V)$, et il reste \`a en montrer l'unicit\'e. Ceci montrera d'ailleurs que $\nwach(V)$ est stable par $\varphi$, puisque $\nwach(V)$ et $\nwach(V)+\varphi^* \nwach(V)$ satisfont tous les deux les conditions (i), (ii) et (iii). Supposons donc que l'on ait deux $\calE^+$-modules $N_1$ et $N_2$ satisfaisant les trois conditions du th\'eor\`eme. 
La condition (iii) implique qu'il existe $r \geq 0$ tel que $X^r N_1 \subset N_2$. Supposons que $r \geq 1$. On a dans ce cas une suite exacte $0 \to X^r N_1 \to N_2 \to N_2 / X^r N_1 \to 0$ et en appliquant le lemme du serpent \`a la multiplication par $X$ dans cette suite, on obtient :
\[ 0 \to (N_2 / X^r N_1)[X] \to X^r N_1  / X^{r+1} N_1 \to N_2 / X N_2 \to  N_2 / (X N_2 + X^r N_1) \to 0. \]
Par la condition (ii), un sous-groupe ouvert du groupe $\Gamma$ agit par $\chi^r$ sur $X^r N_1  / X^{r+1} N_1$ et par $1$ sur $N_2 / X N_2$ ce qui fait que, si $r \geq 1$, alors l'application $X^r N_1  / X^{r+1} N_1 \to N_2 / X N_2$ est nulle et donc que l'application $N_2 / X N_2 \to  N_2 / (X N_2 + X^r N_1)$ est un isomorphisme. On en conclut que si $r \geq 1$, alors en fait $X^r N_1 \subset X N_2$ et donc que $X^{r-1} N_1 \subset N_2$; en it\'erant ce proc\'ed\'e, on se ram\`ene \`a $r=0$ ce qui fait que $N_1 \subset N_2$. Par sym\'etrie, on en conclut aussi que $N_2 \subset N_1$ et donc finalement que $N_1 = N_2$, ce qui termine la d\'emonstration du th\'eor\`eme.
\end{proof}

Si $V$ est une repr\'esen\-tation apc telle que $V(1)$ est positive, alors on voit que $\nwach(V) = X \cdot \nwach(V(1))$, ce qui justifie la d\'efinition suivante.

\begin{defi}\label{wachgen}
Si $V$ est une repr\'esen\-tation apc, et $h \geq 0$ un entier tel que $V(-h)$ est positive, alors on pose $\nwach(V)=X^{-h} \nwach(V(-h))$.
\end{defi}

Comme $\nwach(V) \subset \bplus \otimes_{\Qp} V$, on peut composer les applications $\iota_n : \bplus \to \bdR^+$ et $\theta : \bdR^+ \to \Cp$ pour obtenir une application $\theta \circ \iota_n : \nwach(V) \to \Cp \otimes_{\Qp} V$, qui est nulle sur $Q_n(X) \nwach(V)$.

Comme $V$ est une repr\'esen\-tation de Hodge-Tate de $\gal$, on a une d\'ecomposi\-tion $(\Cp\otimes_{\Qp} V)^{\galinf} \simeq \oplus_{j=1}^d \widehat{L}_{\infty} (-h_j)$, o\`u les $h_j$ sont les oppos\'es des poids de Hodge-Tate de $V$, et on peut montrer (cf. \cite[theorem 3]{Sn80} par exemple) que la r\'eunion $\dsen(V) \= (\Cp \otimes_{\Qp} V)^{\galinf}_{\mathrm{fini}}$ des sous-$F_{\infty}$-espaces vectoriels de
dimension finie stables par $\Gamma$ de $(\Cp \otimes_{\Qp} V)^{\galinf}$ est \'egale \`a $\oplus_{j=1}^d L_{\infty} (-h_j)$.

Le (ii) du th\'eor\`eme \ref{exunmod} nous renseigne sur le module de Wach {\og \'evalu\'e en $0$ \fg} et le lemme suivant nous renseigne sur le module de Wach {\og \'evalu\'e en $\zeta_{p^n}-1$ pour $n \geq 1$ \fg}.

\begin{lemm}\label{injthion}
Si $n \geq 1$, alors l'application $\theta \circ \iota_n : \nwach(V) / Q_n(X) \nwach(V) \to \Cp \otimes_{\Qp} V$ est injective et son image est incluse dans $\dsen(V)$.
\end{lemm}

\begin{proof}
Commen\c{c}ons par montrer que l'application d\'eduite de $\theta \circ \iota_n$ est injective. Si $y \in \nwach(V)$ est tel que $\theta \circ \iota_n(y) =0$, alors $y$ est divisible par $Q_n(X)$ dans $\bplus \otimes_{\Qp} V$ et donc dans $\dplus(V)$. La d\'efinition de $\nwach(V)$ donn\'ee dans la d\'emonstration du th\'eor\`eme \ref{exunmod} montre alors que $y$ est divisible par $Q_n(X)$ dans $\nwach(V)$ et donc que $y \in Q_n(X) \nwach(V)$ ce qui fait que l'application $\theta \circ \iota_n : \nwach(V) / Q_n(X) \nwach(V) \to \Cp \otimes_{\Qp} V$ est injective. 

L'image de cette application est un sous-$F_n$-espace vectoriel de dimension finie et stable par $\Gamma$ de $(\Cp \otimes_{\Qp} V)^{\galinf}$, ce qui fait qu'il est inclus dans $\dsen(V)$.
\end{proof}

On en d\'eduit une application $L_\infty \otimes_{L_n} \nwach(V) / Q_n(X) \nwach(V) \to \dsen(V)$. Nous verrons dans la suite de cet article que si $n \geq m(V)$, alors cette application est une bijection.

Pour terminer ce paragraphe, donnons le lemme technique ci-dessous, qui nous servira dans la suite.

\begin{lemm}\label{egvwd}
Si $h \geq 0$ et si $y \in \calR^+ \otimes_{\calE^+} \nwach(V)$ sont tels que $\gamma(y)=\chi(\gamma)^h y$ pour tout $\gamma$ dans un sous-groupe ouvert de $\Gamma$, alors $y \in X^h \calR^+ \otimes_{\calE^+} \nwach(V)$. 
\end{lemm}

\begin{proof}
Si $\gamma$ appartient de plus \`a un sous-groupe ouvert de $\Gamma$ qui agit trivialement sur $\nwach(V)/X \nwach(V)$, alors on voit que d'une part $\gamma - \chi(\gamma)^j$ envoie $X^j \calR^+ \otimes_{\calE^+} \nwach(V)$ dans $X^{j+1} \calR^+ \otimes_{\calE^+} \nwach(V)$ pour tout $j \geq 0$, et donc que :
\[ (\gamma-1) \cdot (\gamma - \chi(\gamma)) \cdots (\gamma - \chi(\gamma)^{h-1})  \calR^+ \otimes_{\calE^+} \nwach(V) \subset X^h \calR^+ \otimes_{\calE^+} \nwach(V) \]
et d'autre part que :
\[ y = \frac{(\gamma-1) \cdot (\gamma - \chi(\gamma)) \cdots (\gamma - \chi(\gamma)^{h-1})y}{(\chi(\gamma)^h-1) \cdot (\chi(\gamma)^h - \chi(\gamma)) \cdots (\chi(\gamma)^h - \chi(\gamma)^{h-1})}, \]
d'o\`u l'on d\'eduit le lemme.
\end{proof}

\subsection{De $\nwach(V)$ \`a $\dcris(V)$}\label{mwmf}

Le th\'eor\`eme \ref{lbcris} combin\'e au th\'eor\`eme \ref{exunmod} nous montre que si $V$ est une repr\'esen\-tation apc, alors $\dcris(V) \subset \calR[1/t] \otimes_{\calE^+} \nwach(V)$. On peut en fait pr\'eciser ce r\'esultat, c'est l'objet du th\'eor\`eme ci-dessous.

\begin{theo}\label{dcrinwac}
Si $V$ est une repr\'esen\-tation apc positive, alors $\dcris(V) = (\calR^+ \otimes_{\calE^+} \nwach(V))^{\Gamma_n}$ pour $n$ assez grand.
\end{theo}

\begin{proof}
Le d\'ebut de la d\'emonstration est le m\^eme que celui de \cite[prop II.2.1]{Be1} : on choisit une base de $\nwach(V)$, et on appelle $P$ et $G$ les matrices de $\varphi$ et d'un \'el\'ement $\gamma \in \Gamma$ (diff\'erent de $1$ et qui agit trivialement sur $\nwach(V)/X \nwach(V)$) dans cette base, ce qui fait que $P \varphi(G) = G \gamma(P)$. En particulier, on a $\det P \varphi(\det G) = \det G \gamma(\det P)$ et comme $\det G$ est une unit\'e, le lemme \ref{gamstable} montre qu'il existe $\alpha_0,\hdots,\alpha_n$ et une unit\'e $U$ tels que $\det P = U X^{\alpha_0} Q_1(X)^{\alpha_1} \cdots Q_n(X)^{\alpha_n}$. On a par ailleurs $\alpha_0 = 0$ car $\det P / \gamma(\det P) = \det G / \varphi(\det G)$ et le membre de gauche est \'egal \`a $\chi(\gamma)^{\alpha_0}$ modulo $X$ alors que le membre de droite vaut $1$ modulo $X$. 

Par le th\'eor\`eme \ref{lbcris}, on a $\dcris(V) \subset \calR \otimes_{\calE^+} \nwach(V)$; soit $M \in \mathrm{M}(d,\calR)$ la matrice d'une base de $\dcris(V)$ dans la base de $\nwach(V)$ que l'on a choisie pr\'ecedemment. Si $P_0$ d\'enote la matrice de $\varphi$ sur $\dcris(V)$, si $Q$ d\'enote la transpos\'ee de la matrice des cofacteurs de $P$, et si $N=X^{\alpha_1} \varphi(X)^{\alpha_2} \cdots \varphi^{n-1}(X)^{\alpha_n} M$, alors un petit calcul montre que $\varphi(N) = U^{-1} QNP_0$ et le corollaire I.4.2 de \cite{Be1} montre qu'alors $N \in \mathrm{M}(d,\calR^+)$. Il existe donc $n,h \geq 0$
tels que $\dcris(V) \subset \varphi^n(X^{-h}) \calR^+ \otimes_{\calE^+} \nwach(V)$.

On en d\'eduit que si $y \in \dcris(V)$, alors $t^h y \in \calR^+ \otimes_{\calE^+} \nwach(V)$ et le lemme \ref{egvwd} nous dit alors que $(t/X)^h y \in \calR^+ \otimes_{\calE^+} \nwach(V)$, ce qui fait finalement que l'on a $(\varphi^n(X)/X)^h \dcris(V) \subset \calR^+ \otimes_{\calE^+} \nwach(V)$. En appliquant $\varphi^n$ aux deux membres, on trouve que $(\varphi^{2n}(X)/\varphi^n(X))^h \dcris(V) \subset \calR^+ \otimes_{\calE^+} \nwach(V)$ et comme $\varphi^n(X)/X$ et $\varphi^{2n}(X)/\varphi^n(X)$ sont premiers entre eux, on en d\'eduit que $\dcris(V) \subset \calR^+ \otimes_{\calE^+} \nwach(V)$.
\end{proof}

Par le th\'eor\`eme \ref{dcrinwac}, on a une inclusion $\calR^+ \otimes_L \dcris(V) \subset \calR^+ \otimes_{\calE^+} \nwach(V)$; la th\'eorie des diviseurs \'el\'ementaires marche pour l'anneau $\calR^+$ (cf. \cite[\S 4.2]{Be1}), et nous notons $\delta_1,\hdots,\delta_d$ les id\'eaux (principaux) ainsi d\'etermin\'es, ordonn\'es de telle mani\`ere que $\delta_1 \mid \cdots \mid \delta_d$. Comme $\calR^+ \otimes_L \dcris(V)$ et $\calR^+ \otimes_{\calE^+} \nwach(V)$ sont munis d'une action de $\Gamma$, les id\'eaux $\delta_i$ sont stables sous cette action et par le lemme \ref{gamstable}, il existe des entiers $\{\beta_{n,i}\}_{n \geq 0, 1 \leq i \leq d}$ tels que :
\[ \delta_i =  X^{\beta_{0,i}} \cdot \prod_{n=1}^\infty \left( \frac{Q_n(X)}{p}\right)^{\beta_{n,i}} \cdot \calR^+. \]
Le calcul des $\beta_{n,i}$ est un point important de la th\'eorie des modules de Wach des repr\'esen\-tations apc.

\begin{theo}\label{calcdivelem}
Soit $V$ une repr\'esen\-tation apc positive, soient $h_1,\hdots,h_d$ les oppos\'es des poids de Hodge-Tate de $V$, rang\'es dans l'ordre croissant, et soit $t_W^n(V) \= \beta_{n,1} + \cdots + \beta_{n,d}$ pour $n \geq 1$, o\`u les $\beta_{n,i}$ sont d\'efinis ci-dessus.

On a alors :
\begin{itemize}
\item[(i)] $\beta_{0,i} = 0$ pour tout $i$;
\item[(ii)] si $n \geq 1$, alors $\beta_{n,i} \in \{h_1,\hdots,h_d\}$ pour tout $i$;
\item[(iii)] si $n \geq m(V)$, alors en fait $\beta_{n,i} = h_i$ pour tout $i$; 
\item[(iv)] on a $t_W^1(V) \leq t_W^2(V) \leq \cdots \leq t_W^{n(V)} (V) = t_H(V) \= h_1 + \cdots + h_d$. 
\end{itemize}
\end{theo}

La d\'emonstration de ce th\'eor\`eme n\'ecessite quelques lemmes pr\'eparatoires. 

Comme $\calR^+ \otimes_L \dcris(V) \subset \calR^+ \otimes_{\calE^+} \nwach(V)$, on peut {\og localiser et compl\'eter \fg} via l'application $\iota_n : \calR^+ \to L_n[[t]]$ et on en d\'eduit que pour tout $n \geq 0$ : 
\[ L_n[[t]] \otimes_L \dcris(V) \subset L_n[[t]] \otimes^{\iota_n}_{\calE^+} \nwach(V). \]
Les diviseurs \'el\'ementaires de cette inclusion sont alors les id\'eaux : $(t^{\beta_{n,1}}),\hdots,(t^{\beta_{n,d}})$.

\begin{prop}\label{gambet}
On a $\nwach(V)/Q_n(X)\nwach(V) = \oplus_{i=1}^d L_n(-\beta_{n,i})$ en tant que repr\'esen\-tations de $\Gamma_{\max(n,n(V))}$.
\end{prop}

\begin{proof}
Etant donn\'e que $L_n((t)) \otimes_L \dcris(V) = L_n((t)) \otimes^{\iota_n}_{\calE^+} \nwach(V)$, le $L_n$-module $\nwach(V)/Q_n(X) = L_n[[t]] \otimes^{\iota_n}_{\calE^+} \nwach(V) / t$ est un sous-quotient stable par $\Gamma$ de $L_n((t)) \otimes_L \dcris(V)$ et comme $L_n$ est un $L_n[\Gamma]$-module simple, il existe des entiers $\alpha_1 \leq \cdots \leq \alpha_d$ tels qu'il est isomorphe \`a $\oplus_{i=1}^d L_n(-\alpha_i)$ en tant que repr\'esen\-tations de $\Gamma_{\max(n,n(V))}$. Il s'agit donc de montrer que $\alpha_i = \beta_{n,i}$.

Par la th\'eorie des diviseurs \'el\'ementaires, il existe une base de $L_n[[t]] \otimes^{\iota_n}_{\calE^+} \nwach(V)$ telle que les \'el\'ements de cette base, multipli\'es par les $t^{\beta_{n,i}}$, forment une base de $L_n[[t]] \otimes_L \dcris(V)$. Si $G$ est la matrice d'un \'el\'ement $\gamma \in \Gamma_{\max(n,n(V))}$ dans cette base de $L_n[[t]] \otimes^{\iota_n}_{\calE^+} \nwach(V)$ et si $\Lambda = \mathrm{diag} (t^{\beta_{n,1}}, \hdots, t^{\beta_{n,d}})$, alors la matrice de $\gamma$ dans la base construite ci-dessus de $L_n[[t]] \otimes_L \dcris(V)$ est $\Lambda^{-1} G \gamma(\Lambda) = \Lambda^{-1} G \Lambda \cdot \Lambda^{-1} \gamma(\Lambda)$. Le groupe $\Gamma_{\max(n,n(V))}$ agit trivialement sur $\dcris(V)$ ce qui fait que 
$\Lambda^{-1} G \gamma(\Lambda) \in \mathrm{Id} + t \mathrm{M}_d(L_n[[t]])$. Comme on a par ailleurs $\gamma(\Lambda)^{-1} \Lambda = \mathrm{diag}(\chi(\gamma)^{-\beta_{n,1}}, \hdots, \chi(\gamma)^{-\beta_{n,d}})$, on en d\'eduit que  :
\[ \Lambda^{-1} G \Lambda = \Lambda^{-1} G \gamma(\Lambda) \cdot \gamma(\Lambda)^{-1} \Lambda \in \mathrm{diag}(\chi(\gamma)^{-\beta_{n,1}}, \hdots, \chi(\gamma)^{-\beta_{n,d}}) + t \mathrm{M}_d(L_n[[t]]). \]

Le polyn\^ome caract\'eristique de $G$ est le m\^eme que celui de $\Lambda^{-1} G \Lambda$, et en r\'eduisant ce polyn\^ome modulo $t$ et en comparant les valeurs propres de l'action de $\gamma$, on en d\'eduit que $\alpha_i = \beta_{n,i}$.
\end{proof}

\begin{lemm}\label{htwtsee}
Soient $m$ et $n$ deux entiers $\geq 1$, $W_n$ et $W_\infty$ un $L_n$-module et un $L_\infty$-module tous deux libres et munis d'actions semi-lin\'eaires de $\Gamma$ et tels que, en tant que repr\'esen\-tations de $\Gamma_m$, ils soient isomorphes \`a $\oplus_{i=1}^d L_n(a_i)$ et  \`a $\oplus_{i=1}^d L_\infty(b_i)$ pour des entiers $a_1 \leq \cdots \leq a_d$ et $b_1 \leq \cdots \leq b_d$. 

S'il existe une injection $\Gamma$-\'equivariante de $W_n$ dans $W_\infty$, alors $a_i \in \{  b_1, \hdots, b_d \}$ pour tout $i$, et si en plus
$m \leq n$, alors en fait $a_i = b_i$ pour tout $i$.
\end{lemm}

\begin{proof}
Soit $\nabla$ l'op\'erateur $F_\infty$-lin\'eaire provenant de l'action de l'alg\`ebre de Lie de $\Gamma$. On a $W_n = \oplus_{i=1}^d W_n^{\nabla = a_i}$ et $W_\infty = \oplus_{i=1}^d W_\infty^{\nabla = b_i}$ et en comparant les valeurs propres de $\nabla$, on voit que $a_i \in \{  b_1, \hdots, b_d \}$ pour tout $i$. 

Pour montrer que $a_i = b_i$ pour tout $i$ si $m \leq n$, il suffit voir que sous cette hypoth\`ese, l'application naturelle $F_\infty \otimes_{F_n} W_n \to W_\infty$ est bijective. Pour des raisons de dimensions, il suffit de voir qu'elle est injective. Si ce n'\'etait pas le cas, il existerait un \'el\'ement $\sum \lambda_j \otimes w_j \in F_\infty \otimes_{F_n} W_n$ dont l'image dans $W_\infty$ est nulle, avec $w_j$ appartenant \`a l'un des $L_n(a_i)$, et de longueur minimale. En faisant agir $\Gamma_n$ et en utilisant l'hypoth\`ese de minimalit\'e de la longueur de la relation, on voit que l'on a forc\'ement $\lambda_j / \lambda_1 \in F_n$, ce qui fait que l'application $F_\infty \otimes_{F_n} W_n \to W_\infty$ est bien bijective.
\end{proof}

\begin{proof}[D\'emonstration du th\'eor\`eme \ref{calcdivelem}]
Commen\c{c}ons par remarquer que quitte \`a tordre $V$ par un caract\`ere d'ordre fini (ce qui ne change pas la valeur des $\beta_{n,i}$ ou des $h_i$), on peut supposer que $m(V)=n(V)$.

Le fait que $\beta_{0,i} = 0$ r\'esulte du fait que $\nwach(V)/X$ est isomorphe \`a $\oplus_{i=1}^d L(-\beta_{0,i})$ en tant que repr\'esen\-tations de $\Gamma_{n(V)}$ (par un raisonnement analogue \`a celui de la proposition \ref{gambet}) d'une part et du fait qu'un sous-groupe ouvert de $\Gamma$ agit trivialement sur $\nwach(V)/X$ par d\'efinition d'autre part. 

Le fait que pour tous $n,i \geq 1$, $\beta_{n,i}$ est l'un des $h_j$, et que $\beta_{n,i} = h_i$ si $n \geq n(V)$ r\'esulte de la proposition \ref{gambet}, et de la r\'eunion des lemmes \ref{injthion} et \ref{htwtsee}.

On en d\'eduit imm\'ediatement que $t_W^n(V) = \sum_{i=1}^d h_i$ pour tout $n \geq n(V)$. Si $\delta = \delta_1 \times \cdots \times \delta_d$, alors un calcul simple, qui utilise le fait que $\varphi^*(\calR^+ \otimes_L \dcris(V)) =  \calR^+ \otimes_L \dcris(V)$, montre que l'id\'eal de $\calR^+$ engendr\'e par le d\'eterminant de $\varphi$ sur $\nwach(V)$ est engendr\'e par $\delta/\varphi(\delta)$, c'est-\`a-dire qu'il est engendr\'e par :
\[ Q_1(X)^{t_W^1(V)} \cdot Q_2(X)^{t_W^2(V)-t_W^1(V)} \cdots Q_{n(V)}(X)^{t_W^{n(V)}(V)-t_W^{n(V)-1}(V)}. \]
Comme $\nwach(V)$ est stable par $\varphi$, on en d\'eduit que $t_W^1(V) \leq t_W^2(V) \leq \cdots \leq t_W^{n(V)} (V) $.
\end{proof}

\begin{rema}
Dans le cas o\`u $m(V)=1$, on retrouve les r\'esultats de \cite{Be1}. Dans les cas o\`u $m(V) > 1$, nous verrons des exemples (cf. la proposition \ref{caspartic}) qui montrent que le d\'eterminant de $\varphi$ sur $\nwach(V)$ peut effectivement \^etre divisible par $Q_i(X)$ avec $i > 1$.
\end{rema}

Les trois r\'esultats ci-dessous sont des corollaires imm\'ediats du th\'eor\`eme \ref{calcdivelem} et de sa d\'emonstration.

\begin{coro}\label{inclphi}
Si $V$ est une repr\'esen\-tation apc positive, alors :
\[ Q_1(X)^{t_W^1(V)} \cdot Q_2(X)^{t_W^2(V)-t_W^1(V)} \cdots Q_{n(V)}(X)^{t_W^{n(V)}(V)-t_W^{n(V)-1}(V)} \cdot \nwach(V) \subset \varphi^*(\nwach(V)). \]
\end{coro}

\begin{coro}\label{inclrev}
Si $V$ est une repr\'esen\-tation apc dont les poids de Hodge-Tate sont dans l'intervalle $[-h;0]$, alors $\calR^+ \otimes_{\calE^+} \nwach(V) \subset (t/X)^{-h} \calR^+ \otimes_L \dcris(V)$. 
\end{coro}

\begin{coro}\label{isocomp}
Si $V$ est une repr\'esen\-tation apc positive, alors : 
\begin{multline*} 
X^{t_W^1(V)} \cdot \varphi(X)^{t_W^2(V)-t_W^1(V)} \cdots \varphi^{n(V)-1}(X)^{t_W^{n(V)}(V)-t_W^{n(V)-1}(V)} \cdot \bplus \otimes_{\Qp} V \\ 
\subset (L \otimes_{\Qp} \bplus) \otimes_{\calE^+} \nwach(V). 
\end{multline*}
En particulier, si les poids de Hodge-Tate de $V$ sont dans un intervalle $[-h,0]$, alors $\varphi^{n(V)-1}(X^{dh}) \bplus \otimes_{\Qp} V \subset (L \otimes_{\Qp} \bplus) \otimes_{\calE^+} \nwach(V)$.
\end{coro}

\begin{prop}\label{wdsame}
Si $\lambda : \dcris(V) \to \nwach(V) / X \nwach(V)$ d\'enote l'application d\'eduite de l'inclusion de $\dcris(V)$ dans $\calR^+ \otimes_{\calE^+} \nwach(V)$ et de la projection de $\calR^+ \otimes_{\calE^+} \nwach(V)$ sur $\nwach(V) / X \nwach(V)$, alors $\lambda$ est un isomorphisme qui commute \`a $\varphi$ et \`a l'action de $\Gamma$.
\end{prop}

\begin{proof}
Il est clair que $\lambda$ commute \`a $\varphi$ et \`a l'action de $\Gamma$, et il suffit donc de montrer que c'est un isomorphisme. Pour des raisons de dimension, il suffit de montrer que $\lambda$ est injective, c'est-\`a-dire que $\dcris(V) \cap  X \calR^+ \otimes_{\calE^+} \nwach(V) = 0$. Nous allons montrer par r\'ecurrence sur $j \geq 1$ que $\dcris(V) \cap  X \calR^+ \otimes_{\calE^+} \nwach(V) \subset X^j \calR^+ \otimes_{\calE^+} \nwach(V)$. Pour $j=1$, c'est trivialement vrai. 

Choisissons $\gamma \in \Gamma$ diff\'erent de $1$ mais agissant trivialement sur $\dcris(V)$ et sur $\nwach(V) / X$. L'application $\gamma-\chi(\gamma)^j$ est alors bijective sur $\dcris(V)$, et envoie $X^j \calR^+ \otimes_{\calE^+} \nwach(V)$ dans $X^{j+1} \calR^+ \otimes_{\calE^+} \nwach(V)$ ce qui fait que $\dcris(V) \cap  X^j \calR^+ \otimes_{\calE^+} \nwach(V) \subset X^{j+1} \calR^+ \otimes_{\calE^+} \nwach(V)$ et la proposition est d\'emontr\'ee.
\end{proof}

\subsection{Une autre construction de $\nwach(V)$}\label{altwach}

Dans ce paragraphe, on suppose que $V$ est une repr\'esen\-tation apc dont les poids de Hodge-Tate sont dans l'intervalle $[-h;0]$.
Par le corollaire \ref{inclrev}, on a $\nwach(V) \subset t^{-h} \calR^+ \otimes_L \dcris(V)$ et l'objet de ce paragraphe est de montrer comment l'on peut r\'ecup\'erer $\nwach(V)$ comme sous-module de $t^{-h} \calR^+ \otimes_L \dcris(V)$.

\begin{lemm}\label{surfil}
Si $m \geq 0$, alors l'image de $L_m[[t]] \otimes^{\iota_m}_{\calE^+} \nwach(V)$ dans $L_m((t)) \otimes_L \dcris(V)$ est incluse dans $\Fil^0 ( t^{-h} L_m[[t]] \otimes_L \dcris(V) )$ et si $m \geq m(V)$, alors l'application :
\[ L_m[[t]] \otimes^{\iota_m}_{\calE^+} \nwach(V) \to \Fil^0 ( t^{-h} L_m[[t]] \otimes_L \dcris(V) ) \]
est un isomorphisme.
\end{lemm}

\begin{proof}
L'inclusion suit du fait que $\nwach(V) \subset \bplus \otimes_{\Qp} V$ ce qui fait que l'image de $\nwach(V)$ par $\iota_m$ est bien dans le $\Fil^0$ de  $L_m((t)) \otimes_L \dcris(V)$, et l'hypoth\`ese sur les poids de $V$ implique que $\Fil^0 ( L_m((t)) \otimes_L \dcris(V) ) = \Fil^0 ( t^{-h} L_m[[t]] \otimes_L \dcris(V) )$.

Supposons maintenant que $m \geq m(V)$; en consid\'erant (par exemple) une base de $L_m \otimes_L \dcris(V)$ adapt\'ee \`a la filtration, on voit que le d\'eterminant de l'inclusion de $\Fil^0 ( t^{-h} L_m[[t]] \otimes_L \dcris(V) )$ dans $t^{-h} L_m[[t]] \otimes_L \dcris(V)$ est le m\^eme que celui de l'inclusion de $L_m[[t]] \otimes^{\iota_m}_{\calE^+} \nwach(V)$ dans $t^{-h} L_m[[t]] \otimes_L \dcris(V)$ (c'est-\`a-dire $t^{dh-\sum_{i=1}^d h_i}$) ce qui fait que l'application du lemme est bien un isomorphisme dans ce cas-l\`a.
\end{proof}

\begin{defi}\label{defmfil}
On appelle $\mfil(V)$ l'ensemble des $f \in t^{-h} \calR^+ \otimes_L \dcris(V)$ tels que $\iota_m(f) \in \Fil^0 ( t^{-h} L_m[[t]] \otimes_L \dcris(V) )$ pour tout $m \geq m(V)$.
\end{defi}

\begin{prop}\label{wacrob}
Le $\calR^+$-module $\mfil(V)$ est libre de rang $d$, stable par $\psi$, et il v\'erifie :
\[ X^{-h} \calR^+ \otimes_{\calE^+} \nwach(V) 
\subset \mfil(V) \subset \varphi^{m(V)-1}(X)^{-h} \cdot \calR^+ \otimes_{\calE^+} \nwach(V). \] 
\end{prop}

\begin{proof}
Le fait que $X^{-h} \calR^+ \otimes_{\calE^+} \nwach(V) 
\subset \mfil(V)$ suit imm\'ediatement du corollaire \ref{inclrev} et de la premi\`ere partie du lemme \ref{surfil}. Comme les applications $\iota_m$ sont continues, $\mfil(V)$ est un sous-module ferm\'e de $t^{-h} \calR^+ \otimes_L \dcris(V)$ et par le th\'eor\`eme de Forster (cf. par exemple \cite[th\'eor\`eme 4.10]{Ber1}), il est libre de rang $\leq d$ et comme il contient $X^{-h} \calR^+ \otimes_{\calE^+} \nwach(V)$, son rang est en fait \'egal \`a $d$.

Montrons \`a pr\'esent que $\mfil(V)$ est stable par $\psi$. Il est clair que $t^{-h} \calR^+ \otimes_L \dcris(V)$ est stable par $\psi$, et pour conclure, nous allons montrer que si $\iota_{m+1}(f) \in \Fil^0 ( t^{-h} L_{m+1}[[t]] \otimes_L \dcris(V) )$, alors $\iota_m(\psi(f)) \in \Fil^0 ( t^{-h} L_m[[t]] \otimes_L \dcris(V) )$. Pour cela, rappelons que si $g(X) \in \calR^+$, alors :
\[ (\varphi \circ \psi)(g(X)) = \frac{1}{p} \sum_{\eta^p=1} g(\eta (1+X)-1) \] et donc que :
\[ \iota_m(\psi(f)) = \frac{1}{p} \sum_{\eta^p=1} \iota_{m+1}(f)(\eta (1+X)-1). \]

Si $\gamma \in \Gamma$ est tel que $\gamma(\zeta_{p^{m+1}}) = \eta \zeta_{p^{m+1}}$, et si $z(t) \= f^{\varphi^{-(m+1)}} (\zeta_{p^{m+1}} \exp(t/p^{m+1})-1)$, alors $f^{\varphi^{-(m+1)}} (\eta \zeta_{p^{m+1}} \exp(t/p^{m+1})-1) = \gamma(z(\chi(\gamma)^{-1} t))$ et comme $\Fil^0 ( t^{-h} L_{m+1}[[t]] \otimes_L \dcris(V) )$ est un $L_{m+1}[[t]]$-module qui a une base fix\'ee par $\Gamma$, on en d\'eduit que si $\iota_{m+1}(f) \in \Fil^0 ( t^{-h} L_{m+1}[[t]] \otimes_L \dcris(V) )$, alors $\iota_m(\psi(f)) \in \Fil^0 ( t^{-h} L_m[[t]] \otimes_L \dcris(V) )$ et donc que $\mfil(V)$ est stable par $\psi$. 

Montrons enfin que $\mfil(V) \subset \varphi^{m(V)-1}(X)^{-h} \cdot \calR^+ \otimes_{\calE^+} \nwach(V)$; pour cela, fixons une base $n_1,\hdots,n_d$ de $\nwach(V)$. Si $f \in \mfil(V)$, alors $f \in t^{-h} 
\calR^+ \otimes_{\calE^+} \nwach(V)$, et l'on peut donc \'ecrire $f = \sum_{i=1}^d f_i n_i$ avec $f_i \in t^{-h} \calR^+$. La deuxi\`eme inclusion de la proposition revient \`a dire que pour tout $m \geq m(V)$, la fonction $f_i$ n'a pas de p\^ole en $\zeta_{p^m}-1$, c'est-\`a-dire que $\iota_m(f_i) \in L_m[[t]]$. Par le lemme \ref{surfil}, la famille $\iota_m(n_1), \hdots, \iota_m(n_d)$ forme une base de $\Fil^0 ( t^{-h} L_m[[t]] \otimes_L \dcris(V) )$. Si $f \in \mfil(V)$, alors par d\'efinition $\iota_m(f) \in \Fil^0 ( t^{-h} L_m[[t]] \otimes_L \dcris(V) )$, et donc $\iota_m(f_i) \in L_m[[t]]$.
\end{proof}

\begin{rema}\label{warbtr}
Si $m(V)=1$, alors on trouve que $\calR^+ \otimes_{\calE^+} \nwach(V)$ s'identifie \`a l'ensemble des $f \in (t/X)^{-h} \calR^+ \otimes_L \dcris(V)$ tels que $\iota_m(f) \in \Fil^0 ( t^{-h} L_m[[t]] \otimes_L \dcris(V) )$ pour tout $m \geq 1$.
\end{rema}

Les calculs pr\'ec\'edents montrent comment r\'ecup\'erer $\calR^+ \otimes_{\calE^+} \nwach(V)$, et la deuxi\`eme \'etape est d'extraire $\nwach(V)$ de $\calR^+ \otimes_{\calE^+} \nwach(V)$. Pour cela, nous avons besoin de la notion d'ordre (de croissance). 

Rappelons (cf. la d\'efinition \ref{dordrer}) qu'une s\'erie $f \in \calR^+$ est d'ordre $r$ si et seulement si quel que soit $\rho$ tel que $0 < \rho < 1$, la suite $\{ p^{-nr} \| f(X) \|_{D(0, \rho^{1/p^n})}   \}_{n \geq 0}$ est born\'ee. On pose alors $\| f \|_{\rho,r} = \sup_{n \geq 0} p^{-nr} \| f(X) \|_{D(0, \rho^{1/p^n})}$. Si $f(X)=\sum_{j=0}^\infty a_n X^n \in \calR^+$, alors on peut montrer (cf. le lemme \ref{ordre} ci-dessous pour un \'enonc\'e pr\'ecis) que $f$ est d'ordre $r$ si et seulement si la suite $\{ (n+1)^{-r} |a_n| \}_{n \geq 0}$ est born\'ee. Si c'est le cas, alors on pose $\| f \|_r = \sup_{n \geq 0} (n+1)^{-r} |a_n|$, et les normes $\|\cdot\|_r$ et $\|\cdot\|_{\rho,r}$ sont \'equivalentes.

La notion d'ordre est \'etendue aux \'el\'ements de $\btrigplus[1/t]$ (voir \cite[\S V.3]{Be2}) ce qui permet de parler de l'ordre des \'el\'ements de $\calR^+[1/t] \otimes_L \dcris(V)$ et de $\calR^+[1/t] \otimes_{\calE^+} \nwach(V)$ (ces deux espaces \'etant bien-s\^ur les m\^emes). Dans $\btrigplus$, on peut utiliser le frobenius pour donner une troisi\`eme d\'efinition de l'ordre : si $f \in \btrigplus$, alors $f$ est d'ordre $r$ si et seulement si quel que soit $\rho$ tel que $0 < \rho < 1$, la suite $\{ p^{-nr} \| \varphi^{-n} (f) \|_{D(0, \rho)} \}_{n \geq 0}$ est born\'ee. On en d\'eduit en particulier le r\'esultat suivant (rappelons que la matrice d'un endomorphisme est dite \^etre sous forme de Jordan si elle est triangulaire sup\'erieure et adapt\'ee \`a la d\'ecomposition en espaces caract\'eristiques).

\begin{lemm}\label{ordcph}
Si l'on se donne une base $\{ e_1, \hdots, e_d \}$ de $\dcris(V)$ dans laquelle la matrice de $\varphi$ est sous forme de Jordan, la valeur propre correspondant \`a $e_i$ \'etant not\'ee $\alpha_i^{-1}$, alors $y = \sum_{i=1}^d y_i \otimes e_i \in \calR^+ \otimes_L \dcris(V)$ est d'ordre $r$ si et seulement si $y_i$ est d'ordre $r+\val(\alpha_i)$.
\end{lemm}

\begin{proof}
Comme la base $\{ e_1, \hdots, e_d \}$ respecte la d\'ecomposition de $\dcris(V)$ en espaces caract\'eristiques, on peut supposer qu'il n'y a qu'une seule valeur propre $\alpha$, et en tordant l'action du frobenius, on peut supposer par ailleurs qur $r=0$ et que $\alpha=1$. Si l'on note $P$ la matrice de $\varphi$, l'ensemble $\{ P^j \}_{j \in \ZZ}$ est alors born\'e (pour la topologie $p$-adique) et on en d\'eduit que $\{ \varphi^{-n}(y) \}_{n \geq 0 }$ est born\'e si et seulement si pour tout $1 \leq i \leq d$, $\{ \varphi^{-n}(y_i) \}_{n \geq 0 }$ est born\'e (pour la norme $\| \cdot \|_{D(0, \rho)}$), c'est-\`a-dire si et seulement si pour tout $1 \leq i \leq d$, la s\'erie $y_i$ est d'ordre $0$.
\end{proof}

\begin{lemm}\label{old335}
Si $g(X) \in \calE^+$, alors l'ensemble des \'el\'ements de $g(X)^{-1} \calR^+ \otimes_{\calE^+} \nwach(V)$ qui sont d'ordre $0$ est $g(X)^{-1} \nwach(V)$.
\end{lemm}

\begin{proof}
On a $g(X)^{-1} \calR^+ \otimes_{\calE^+} \nwach(V) \subset \btrig \otimes_{\Qp} V$ et l'ensemble des \'el\'ements d'ordre $0$ de $\btrig \otimes_{\Qp} V$ est $\btdag \otimes_{\Qp} V$. Le lemme r\'esulte alors de ce que : 
\begin{multline*} 
\btdag \otimes_{\Qp} V \cap g(X)^{-1} \calR^+ \otimes_{\calE^+} \nwach(V) \\
= \calE^\dagger \otimes_{\calE^+} \nwach(V) \cap g(X)^{-1} \calR^+ \otimes_{\calE^+} \nwach(V) 
= g(X)^{-1} \nwach(V), 
\end{multline*}
puisque $g(X)^{-1} \calR^+ \cap \calE^\dagger = g(X)^{-1} \calE^+$.
\end{proof}

\begin{theo}\label{old337}
On se donne une repr\'esen\-tation apc $V$ dont les poids de Hodge-Tate sont dans l'intervalle $[-h;0]$, ainsi qu'une base $\{ e_1, \hdots, e_d \}$ de $\dcris(V)$ dans laquelle la matrice de $\varphi$ est sous forme de Jordan, la valeur propre correspondant \`a $e_i$ \'etant not\'ee $\alpha_i^{-1}$. 
Si l'on appelle $\mfil^0(V)$ l'ensemble des $f = \sum_{i=1}^d t^{-h} f_i \otimes e_i \in t^{-h} \calR^+ \otimes_L \dcris(V)$ tels que :
\begin{itemize}
\item[(i)] $f_i$ est d'ordre $h+\val(\alpha_i)$;
\item[(ii)] $\iota_m(f) \in \Fil^0 ( t^{-h} L_m[[t]] \otimes_L \dcris(V) )$ pour tout $m \geq m(V)$,
\end{itemize}
alors $\mfil^0(V)$ est un $\calE^+$-module libre de rang $d$ stable par $\psi$, et $X^{-h} \nwach(V) \subset \mfil^0(V) \subset \varphi^{m(V)-1}(X)^{-h} \cdot \nwach(V)$.
\end{theo}

\begin{proof}
Par le lemme \ref{ordcph}, l'ensemble $\mfil^0(V)$ est le sous-ensemble de $\mfil(V)$ constitu\'e des \'el\'ements qui sont d'ordre nul. Le fait que  $X^{-h} \nwach(V) 
\subset \mfil^0(V) \subset \varphi^{m(V)-1}(X)^{-h} \cdot \nwach(V)$ suit alors de la proposition \ref{wacrob} et du lemme \ref{old335}. On en d\'eduit imm\'ediatement que $\mfil^0(V)$ est un $\calE^+$-module libre de rang $d$. Comme $\mfil(V)$ est stable par $\psi$, le fait que $\mfil^0(V)$ est stable par $\psi$ suit du fait que si une s\'erie $f \in \calR^+$ est d'ordre $r$, alors $\psi(f)$ est elle-aussi d'ordre $r$ (cf. par exemple \cite[proposition 1.15]{Co2}).
\end{proof}

Si $r$ est un r\'eel, on note $\calR^+_r$ l'anneau des s\'eries d'ordre $r$; cet anneau est un espace de Banach pour la norme $\|\cdot\|_r$. De la d\'efinition de $\mfil^0(V)$ donn\'ee ci-dessus, on d\'eduit une application $\lambda : \mfil^0(V) \to \prod_{i=1}^d \calR^+_{h+\val(\alpha_i)}$, qui \`a $f = \sum_{i=1}^d t^{-h} f_i \otimes e_i$ associe $(f_1,\hdots,f_d)$. 

\begin{prop}\label{comtop}
L'application $\lambda$ d\'efinie ci-dessus est continue et son image est ferm\'ee, si l'on munit $\mfil^0(V)$ de la topologie $p$-adique et $\prod_{i=1}^d \calR^+_{h+\val(\alpha_i)}$ de la topologie provenant de la norme $\| (f_1,\hdots,f_d) \| \= 
\sup \| f_i \|_{h+\val(\alpha_i)}$.
\end{prop}

\begin{proof}
Si $m \geq m(V)$, alors $\Fil^0 (t^{-h} L_m [[t]] \otimes_L \dcris(V))$ est de codimension finie dans $t^{-h} L_m[[t]] \otimes_L \dcris(V)$ et le noyau de l'application :
\[ \prod_{i=1}^d \calR^+_{h+\val(\alpha_i)} \to t^{-h} L_m [[t]] \otimes_L \dcris(V) / \Fil^0 (t^{-h} L_m [[t]] \otimes_L \dcris(V)), \] qui \`a $(f_1,\hdots,f_d)$ associe l'image de $\sum_{i=1}^d \varphi^{-m}(t^{-h} f_i \otimes  e_i)$ est donc un sous-espace ferm\'e $E_m$ de $\prod_{i=1}^d
\calR^+_{h+\val(\alpha_i)}$ ce qui fait que l'image de $\lambda$, $E \= \cap_{m \geq 1} E_m$ est elle aussi ferm\'ee dans $\prod_{i=1}^d \calR^+_{h+\val(\alpha_i)}$. L'application naturelle $\mfil^0(V) \to E$ est un isomorphisme par d\'efinition;   par le th\'eor\`eme de l'image ouverte, pour montrer que $\mfil^0(V) \to E$ est un isomorphisme topologique, il suffit de montrer que l'application
$\mfil^0(V) \to  \prod_{i=1}^d \calR^+_{h+\val(\alpha_i)}$ est continue,
c'est-\`a-dire qu'elle est born\'ee. Pour cela, il suffit d'observer qu'une boule unit\'e de $\mfil^0(V)$ (pour une norme induisant la topologie $p$-adique) est un $\OO_L[[X]]$-module de type fini.
\end{proof}

\begin{prop}\label{shetwa}
Si $V$ est irr\'eductible et de dimension $\geq 2$, alors il existe $\ell \geq 0$ tel que $\dsharp(V) = \psi^\ell(\nwach(V))$.
\end{prop}

\begin{proof}
Le fait que $\varphi(\nwach(V)) \subset \nwach(V)$ implique que $\nwach(V) \subset \psi(\nwach(V))$. Par ailleurs, le fait que $\nwach(V) \subset \mfil^0(V)$ et que $\mfil^0(V)$ est stable par $\psi$ implique que pour tout $k \geq 0$, on a $\psi^k(\nwach(V)) \subset \mfil^0(V)$. La suite des $\psi^k(\nwach(V))$ est donc une suite croissante de sous-$\calE^+$-modules d'un $\calE^+$-module de type fini, et comme $\calE^+$ est noetherien, il existe $\ell \geq 0$ tel que $\psi^\ell(\nwach(V)) = \psi^{\ell+k}(\nwach(V))$ pour tout $k \geq 0$. Le $\calE^+$-module $\psi^\ell(\nwach(V))$ est un treillis stable par $\psi$ et sur lequel $\psi$ est surjectif, et le lemme \ref{sco429} nous dit alors que $\dsharp(V) = \psi^\ell(\nwach(V))$. 
\end{proof}

\begin{coro}\label{lpslm}
Si $V$ est irr\'eductible et de dimension $\geq 2$, alors $\dsharp(V) \subset \mfil^0(V) \subset \varphi^{m(V)-1}(X)^{-h} \cdot \dsharp(V)$.
\end{coro}

\begin{proof}
Cela r\'esulte directement de la proposition \ref{shetwa} ci-dessus et du th\'eor\`eme \ref{old337} puisque $\nwach(V) \subset \dsharp(V)$.
\end{proof}

\begin{lemm}
Si $V$ est irr\'eductible et de dimension $\geq 2$, alors la topologie faible sur $\dsharp(V)$ est induite par la topologie $(p,X)$-adique sur $\mfil^0(V)$ via l'inclusion $\dsharp(V) \subset \mfil^0(V)$. 

Une suite d'\'el\'ements de $\dsharp(V)$ est born\'ee pour la topologie faible si et seulement si elle est born\'ee pour la topologie $p$-adique dans $\mfil^0(V)$.
\end{lemm}

\begin{proof}
Ces deux points suivent du fait que $\dsharp(V)$ est un $\calE^+$-module de type fini.
\end{proof}

\begin{coro}\label{debilcor}
Si $T$ est un $\OO_L$-r\'eseau de $V$ qui est stable par $\gal$, alors $(\projlim_\psi \dsharp(T))^{\mathrm{b}} = \projlim_\psi \dsharp(T)$ et donc $(\projlim_\psi \dsharp(V))^{\mathrm{b}} = L \otimes_{\OO_L} \projlim_\psi \dsharp(T)$.
\end{coro}

\begin{lemm}\label{samelim}
L'application naturelle $(\projlim_\psi \dsharp(V))^{\mathrm{b}} \to (\projlim_\psi \mfil^0(V))^{\mathrm{b}}$ est un isomorphisme. 
\end{lemm}

\begin{proof}
Cette application est bien s\^ur injective, et il faut v\'erifier qu'elle est surjective. Si $T$ est un $\OO_L$-r\'eseau de $V$ stable par $\gal$, et si l'on pose $\mfil^0(T) = \dfont(T) \cap \mfil^0(V)$, alors l'intersection $\cap_{n \geq 0} \psi^n( \mfil^0(T) )$ est un treillis de $\dfont(T)$ (puisqu'il contient $\dsharp(T)$) stable par $\psi$ et sur lequel $\psi$ est surjectif. Par la proposition \ref{col429}, on en d\'eduit que $\cap_{n \geq 0} \psi^n( \mfil^0(T) ) = \dsharp(T)$ et donc que si $(v_n)_{n \geq 0} \in \projlim_\psi \mfil^0(T)$, alors $v_n \in \dsharp(T)$ pour tout $n$, d'o\`u le lemme.
\end{proof}

\subsection{Repr\'esentations apc de dimension $2$ et repr\'esen\-tations du Borel}

On reprend maintenant les notations de la fin du paragraphe \S\ref{cristallin}, c'est-\`a-dire que l'on se donne une repr\'esen\-tation apc $V$ qui est absolument irr\'eductible, de dimension $2$ et dont les poids de Hodge-Tate sont $0$ et $k-1$ pour un entier $k \geq 2$. En particulier, si l'on pose $h=k-1$, alors les r\'esultats des paragraphes \S\S3.1-3.3 s'appliquent \`a $V(-h)$.

La proposition ci-dessous ne nous servira pas, mais elle a l'inter\^et de fournir des exemples de modules de Wach pour lesquels le d\'eterminant de $\varphi$ est divisble par $Q_i(X)$ avec $i \geq 2$ (contrairement \`a ce qui se passe si $m(V)=1$). Rappelons que pour les repr\'esen\-tations qui nous int\'eressent, $m(V)$ est le plus petit entier $m \geq 1$ tel que $G(\beta \alpha^{-1}) \in F_m$.

\begin{prop}\label{caspartic}
Dans les notations du th\'eor\`eme \ref{calcdivelem}, on a $t_W^i(V(-h)) = 0$ pour $i \leq m(V)-1$ et $t_W^i(V(-h)) = h$ pour $i \geq m(V)$. En particulier, le d\'eterminant de $\varphi$ sur $\nwach(V(-h))$ est $Q_{m(V)}^h$.  
\end{prop}

\begin{proof}
Observons que comme les oppos\'es des poids de Hodge-Tate de $V(-h)$ sont $0$ et $h$, on a forc\'ement $t_W^n(V(-h)) \in \{0; h\}$. Si $m \leq m(V)-1$, alors par d\'efinition $G(\beta \alpha^{-1}) \notin F_m$ ce qui fait que $\Fil^0(t^{-h} L_m[[t]] \otimes_L \dcris(V(-h))) = L_m[[t]] \otimes_L \dcris(V(-h))$ et le lemme \ref{surfil} montre alors que les $\beta_{m,i}$ du th\'eor\`eme \ref{calcdivelem} sont tous les deux nuls et donc que $t_W^m(V(-h)) = 0$ dans ce cas. Enfin, le (iii) du th\'eor\`eme \ref{calcdivelem} implique que $t_W^i(V(-h)) = h$ pour $i \geq m(V)$.
\end{proof}

Revenons \`a pr\'esent \`a $\dsharp(V)$. Par le corollaire \ref{lpslm}, on a une inclusion $\dsharp(V(-h)) \subset \mfil^0(V(-h))$ et on en d\'eduit des applications (bien s\^ur toujours injectives) : 
\[ (\projlim_\psi \dsharp(V(-h)))^{\mathrm{b}} \to \projlim_\psi \mfil^0(V(-h)) \to \projlim_\psi t^{-h} \calR^+ \otimes_L \dcris(V(-h)). \] 
Remarquons que $\dsharp(V(-h)) = \dsharp(V)(-h)$ et que $\dcris(V(-h)) = t^h \dcris(V)(-h)$, ce qui fait que l'on a une injection $(\projlim_\psi \dsharp(V))^{\mathrm{b}} \to \projlim_\psi \calR^+ \otimes_L \dcris(V)$.

Nous en arrivons maintenant au r\'esultat principal de ce chapitre.

\begin{theo}\label{limproj}
Supposons $\alpha \neq \beta$.
Si $( w_{\alpha,n} )_{n \geq 0}$ et $( w_{\beta,n} )_{n \geq 0}$ sont deux suites d'\'el\'ements de $\calR^+$, alors la suite : 
\[ ( w_{\alpha,n} \otimes e_\alpha + w_{\beta,n} \otimes e_\beta )_{n \geq 0} \in \projlim \calR^+ \otimes_L \dcris(V) \]
appartient \`a $(\projlim_\psi \dsharp(V))^{\mathrm{b}}$ si et seulement si :
\begin{itemize}
\item[(i)] pour tout $n \geq 0$, $w_{\alpha,n}$ est d'ordre $\val(\alpha_p)$ et $w_{\beta,n}$ est d'ordre $\val(\beta_p)$, et les deux suites $( \| w_{\alpha,n} \|_{\val(\alpha_p)})_{n \geq 0}$ et $( \| w_{\beta,n}\|_{\val(\beta_p)} )_{n \geq 0}$ sont born\'ees;
\item[(ii)] pour tout $m \geq m(V)$, et pour tout $n \geq 0$, on a $\varphi^{-m}(w_{\alpha,n} \otimes e_\alpha + w_{\beta,n} \otimes e_\beta) \in \Fil^0( L_m[[t]] \otimes_L \dcris(V))$, c'est-\`a-dire $\varphi^{-m}(w_{\alpha,n})\alpha_p^m - G(\alpha \beta^{-1}) \varphi^{-m}(w_{\beta,n})\beta_p^m \in t^{k-1} L_m[[t]]$;
\item[(iii)] pour tout $n \geq 0$, on a $\psi(w_{\alpha,n+1})=\alpha_p^{-1} w_{\alpha,n}$ et $\psi(w_{\beta,n+1})=\beta_p^{-1} w_{\beta,n}$.
\end{itemize}
\end{theo}

\begin{proof}
Posons $w_n = w_{\alpha,n} \otimes e_\alpha + w_{\beta,n} \otimes e_\beta \in t^{-h} \calR^+ \otimes_L \dcris(V(-h))$. 

\begin{itemize}
\item La condition (iii) est \'equivalente \`a dire que pour tout $n \geq 0$, on a $\psi(w_{n+1})=w_n$. 

\item Le fait que $w_{\alpha,n}$ est d'ordre $\val(\alpha_p)$ et $w_{\beta,n}$ est d'ordre $\val(\beta_p)$ avec la condition (ii) sont \'equivalents, par le th\'eor\`eme \ref{old337}, au fait que $w_n \in \mfil^0(V(-h))$.

\item Le fait que les deux suites $\{ \| w_{\alpha,n} \|_{\val(\alpha_p)} \}_{n \geq 0}$ et $\{ \| w_{\beta,n}\|_{\val(\beta_p)} \}_{n \geq 0}$ sont born\'ees est alors \'equivalent, par la proposition \ref{comtop}, au fait que la suite $(w_n)_{n \geq 0}$ est born\'ee pour la topologie $p$-adique dans $\mfil^0(V(-h))$.
\end{itemize}

On voit donc que les conditions (i), (ii) et (iii) sont satisfaites si et seulement si $(w_n)_{n \geq 0} \in (\projlim_\psi \mfil^0(V(-h)))^{\mathrm{b}}$. Le th\'eor\`eme r\'esulte alors du lemme \ref{samelim}, qui nous dit que l'application $(\projlim_\psi \dsharp(V(-h)))^{\mathrm{b}} \to  (\projlim_\psi \mfil^0(V(-h)))^{\mathrm{b}}$ est un isomorphisme (et de l'identification de $\dsharp(V)$ \`a $\dsharp(V(-h))$).
\end{proof}

Le th\'eor\`eme ci-dessus nous fournit une description de $(\projlim_\psi \dsharp(V))^{\mathrm{b}}$ en termes de fonctions analytiques v\'erifiant certaines conditions. Nous reviendrons l\`a-dessus plus loin dans l'article.

Pour le moment, terminons ce paragraphe en d\'efinissant une action du groupe :
\[ \B \= \left\{ \begin{pmatrix} \Qp^\times & \Qp \\ 0 & \Qp^\times \end{pmatrix} \right\} \subset \G \]
sur $(\projlim_\psi \dsharp(V))^{\mathrm{b}}$. Nous allons aussi montrer que cette action est continue, topologiquement irr\'eductible, et que le lemme de Schur est v\'erifi\'e.

On fixe un caract\`ere lisse $\chi$ de $\Qp^\times$ et on munit $(\projlim_\psi \dsharp(V))^{\mathrm{b}}$ d'une action de $\B$ comme suit. Tout \'el\'ement $g \in \B$ peut s'\'ecrire comme produit :
\[ g = \begin{pmatrix} x & 0 \\ 0 & x \end{pmatrix} \cdot
\begin{pmatrix} 1 &  0 \\ 0 & p^j  \end{pmatrix} \cdot 
\begin{pmatrix} 1 &  0 \\ 0 & a  \end{pmatrix} \cdot 
\begin{pmatrix} 1 &  z \\ 0 & 1 \end{pmatrix}, \]
o\`u $x \in \Qp^\times$, $j\in\ZZ$, $a \in \Zp^{\times}$ et $z
\in \Qp$. 

\begin{defi}\label{acgt}
Si $v=(v_i)_{i \geq 0} \in (\projlim_\psi \dsharp(V))^{\mathrm{b}}$, alors on pose pour $i\geq 0$ :
\begin{align*}
\left( \begin{pmatrix} x & 0 \\ 0 & x \end{pmatrix} \cdot v \right)_i
& = \chi^{-1}(x) v_i;  \\
\left( \begin{pmatrix} 1 &  0 \\ 0 & p^j  
\end{pmatrix} \cdot v \right)_i & = v_{i-j} = \psi^j(v_i); \\
\left( \begin{pmatrix} 1 &  0 \\ 0 & a  
\end{pmatrix} \cdot v \right)_i & = \gamma_a^{-1}(v_i), 
\text{ o\`u $\gamma_a \in \Gamma$ est tel que $\eps(\gamma_a) = a$;}\\
\left( \begin{pmatrix} 1 &  z \\ 0 & 1  
\end{pmatrix} \cdot v \right)_i & = 
\psi^j((1+X)^{p^{i+j} z} v_{i+j}),\text{ pour $i+j \geq -{\rm val}(z)$.}
\end{align*}
\end{defi}

On laisse le soin au lecteur de v\'erifier que les formules ci-dessus d\'efinissent bien une action du groupe $\B$ sur $(\projlim_\psi \dsharp(V))^{\mathrm{b}}$. On d\'efinit aussi une structure de $\OO_L[[X]]$-module sur $(\projlim_\psi \dsharp(V))^{\mathrm{b}}$ en posant $(1+X)^z \cdot v \= \left(\begin{smallmatrix} 1 & z \\ 0 & 1 \end{smallmatrix}\right) \cdot v$ pour $z
\in \Zp$. 

\begin{prop}\label{actgcont}
L'application $\B \times (\projlim_\psi \dsharp(V))^{\mathrm{b}} \to
(\projlim_\psi \dsharp(V))^{\mathrm{b}}$ est continue. 
\end{prop}

\begin{proof}
On v\'erifie qu'il suffit de montrer que l'application $\psi : \projlim_{\psi} \dsharp(T) \to \projlim_{\psi} \dsharp(T)$ est continue et que l'application $\BK \times \projlim_{\psi} \dsharp(T) \to \projlim_{\psi}\dsharp(T)$ est continue. Si $E$ et $\{X_i\}_{i\in I}$ sont des espaces topologiques
et si pour tout $i$, on se donne une application
continue $f_i : E \times X_i \to E \times X_i$, alors
l'application $E \times \prod_{i \in I} X_i \to \prod_{i \in
I} X_i$ donn\'ee par $(e,(x_i)_i) \mapsto (f_i(e,x_i))_i$ est
continue. En effet, la diagonale $\Delta_E$ de $\prod_{i \in I}
E$ y est ferm\'ee, et l'application  $(e,(x_i)_i) \mapsto (f_i(e,x_i))_i$
est la composition des applications :
\[ E \times \prod_{i \in I} X_i = \Delta_E \times \prod_{i \in I} X_i
\subset \prod_{i \in I} (E \times X_i) \overset{\prod_{i \in I} f_i} 
{\longrightarrow}  \prod_{i \in I} (E \times X_i) 
\overset{\prod_{i \in I} {\rm pr}_i} {\longrightarrow}
\prod_{i \in I} X_i.  \]
On se ram\`ene donc \`a montrer que si chaque $f_i$
est continue, alors $\prod_{i \in I} f_i$ est continue (pour la
topologie produit) ce qui est laiss\'e en exercice facile au lecteur. 

Afin de montrer la proposition, il suffit donc de montrer que l'application $\psi:
\dsharp(T) \to \dsharp(T)$ est continue, et que l'application
$\BK \times \dsharp(T)  \to \dsharp(T)$ est continue. Commen\c{c}ons par montrer que si $V$ est une repr\'esen\-tation apc, alors (a) l'ensemble $\{p^j \dfont(T) + X^k \dsharp(T)\}_{j,k \geq 0}$ est une base de voisinages de z\'ero pour la topologie faible et (b) l'ensemble $\{p^j \dfont(T) +\varphi(X)^k \dsharp(T)\}_{j,k \geq 0}$ est aussi une base de voisinages de z\'ero pour la topologie faible. La proposition \ref{shetwa} implique que $\nwach(T) \subset \dsharp(T) \subset \varphi^{m(V)-1}(X)^{-h} \nwach(T)$, ce qui montre le point (a) puisque $\nwach(T)$ est un $\OO_L[[X]]$-module libre qui engendre $\dfont(T)$ sur $\calO$. Pour montrer le point (b), et comme $p^j \dfont(T) + \varphi(X)^k \dsharp(T) \subset p^j \dfont(T) + X^k \dsharp(T)$, il suffit de montrer (par exemple) que pour tout $k \geq 0$, il existe $\ell$ tel que $p^j \dfont(T) + X^\ell \dsharp(T) \subset p^j \dfont(T) + \varphi(X)^k \dsharp(T)$. Pour cela, remarquons que si $m
\geq j$ est tel que $p^m \geq k$, alors 
$p^j \dfont(T) + X^{p^{m+1}} \dsharp(T) \subset 
p^j \dfont(T) + \varphi(X)^{p^m} \dsharp(T)
\subset p^j \dfont(T) + \varphi(X)^k \dsharp(T)$ et on peut donc prendre
$\ell=p^{m+1}$. Le fait que l'op\'erateur $\psi: \dsharp(T) \to 
\dsharp(T)$ est continu pour la
topologie faible r\'esulte alors du fait que
$\psi\left(p^j \dfont(T) + \varphi(X)^k \dsharp(T)\right) 
=  p^j \dfont(T)+ X^k \dsharp(T)$. 

Montrons maintenant que l'application naturelle $\BK \times \dsharp(T)
\to \dsharp(T)$ est continue (pour la topologie
faible). Comme 
$\BK$ agit par multiplication par des \'el\'ements de $\OO_L[[X]]$ ou
par action de $\Gamma$, et que la topologie faible de $\dsharp(T)$
est d\'efinie par une base de voisinages de $0$ qui sont des
$\OO_L[[X]]$-modules stables par $\Gamma$, 
on voit pour chaque $g \in \BK$
l'action de $g$ sur $\dsharp(T)$ est continue. Il reste donc \`a montrer que si $W$ est un voisinage de z\'ero dans
$\dsharp(T)$, alors il existe un sous-groupe normal ouvert $U$ de
$\BK$ tel que $u(x)-x \in W$ pour tous $(u,x) \in U \times
\dsharp(T)$. Cela r\'esulte des faits suivants :
\begin{itemize}
\item si $a \in \Zp$ et $n \geq 0$, 
alors $(1+X)^{p^n a}-1 \in (p,X)^{n+1}$;
\item il existe une constante $c$ telle que si $\gamma \in \Gamma_{n(V)}$ et $n \geq 1$,  alors $(\gamma^{p^n }-1)\nwach(T) \subset (p,X)^{n-c} \nwach(T)$,
\end{itemize}
que nous laissons en exercices au lecteur.
\end{proof}

La proposition suivante est le {\og lemme de Schur \fg} pour la
repr\'esen\-tation de $\B$ que l'on vient de d\'efinir.

\begin{prop}\label{schur}
Toute application continue $\B$-\'equivariante :
\[ f: (\projlim_\psi \dsharp(V))^{\mathrm{b}}
\longrightarrow (\projlim_\psi \dsharp(V))^{\mathrm{b}} \]
est scalaire, i.e. est la multiplication par un \'el\'ement de $L$.
\end{prop}

\begin{proof}
Par compacit\'e, on peut (quitte \`a multiplier $f$ par une puissance de $p$) supposer que que $f$ envoie $\projlim_\psi \dsharp(T)$ dans $\projlim_\psi \dsharp(T)$. Notons $\pr : \projlim_\psi \dsharp(T) \to \dsharp(T)$ l'application $v = (v_n)_{n \geq 0} \mapsto v_0$. Commen\c{c}ons par montrer que si $v=(v_n)_{n \geq 0}$, alors $\pr \circ f(v)$ ne d\'epend que de $v_0 = \pr(v)$. Soit $K_n$ l'ensemble des $v \in \projlim_\psi \dsharp(T)$ dont les $n$ premiers termes sont nuls, ce qui fait que pour $n \geq 1$, $K_n$ est un sous-$\OO_L[[X]]$-module ferm\'e et stable par $\psi$ et $\Gamma$ de $\projlim_\psi \dsharp(T)$ et que $\psi(K_n)=K_{n+1}$. 
Si l'on pose $M = \pr \circ f (K_1)$, alors $\psi^n(M) = \pr \circ f (K_{n+1})$.
On a $\cap_{n \geq 1}  \pr \circ f (K_n) = \{ 0 \}$; en effet, si l'on peut \'ecrire $y = \pr \circ f (k_n)$ pour tout $n$, alors comme $k_n \to 0$, on a n\'ecessairement $y=0$. On en d\'eduit que $\cap_{n \geq 1} \psi^n(M) = 0$ et donc par le lemme \ref{noheart} que $M=0$. Ceci revient \`a dire que si l'on se donne $v=(v_n)_{n \geq 0}$, et que $v_0 = 0$, alors $f(v)_0 = 0$.

Pour tout $w \in \dsharp(T)$, soit $\widetilde{w}$ un \'el\'ement de 
$\projlim_\psi \dsharp(T)$ tel que $\widetilde{w}_0 = w$. Les calculs pr\'ec\'edents montrent que l'application $h : \dsharp(T) \to \dsharp(T)$ donn\'ee par $h(w) \= \pr \circ f (\widetilde{w})$ est bien d\'efinie, et qu'elle est $\OO_L[[X]]$-lin\'eaire et commute \`a $\psi$ et \`a l'action de $\Gamma$. Par \cite[proposition 4.7]{Co2}, elle s'\'etend en une
application de $(\varphi,\Gamma)$-modules $h: \dfont(V) \to \dfont(V)$ qui est n\'ecessairement scalaire car $V$ est irr\'eductible. Comme $f(v)_n = h(\psi^{-n} v)$ car $f$ commute \`a $\psi^{-n}$, on en d\'eduit que $f$ est aussi scalaire. 
\end{proof}

\begin{prop}\label{actgirred}
L'action de $\B$ sur $(\projlim_\psi \dsharp(V))^{\mathrm{b}}$ est
topologiquement irr\'e\-ducti\-ble. 
\end{prop}

\begin{proof}
C'est la r\'eunion de la d\'emonstration du corollaire 4.59 et du (iii) de la remarque 5.5 de \cite{Co2} (puisque $V$ et donc $V^*$ est absolument irr\'eductible). 
\end{proof}

\section{Repr\'esentations apc irr\'eductibles de $\G$}

Le but de cette partie est de d\'efinir des espaces de Banach $p$-adiques $\breuil(V)$ munis d'une action lin\'eaire continue de $\G$ et d'en commencer l'\'etude. Les repr\'esentations $\breuil(V)$ sont associ\'ees aux repr\'esentations irr\'eductibles $V$ de dimension $2$ de $\g$ qui deviennent cristallines sur une extension ab\'elienne de $\Q$.

\subsection{Fonctions de classe $\mathcal{C}^r$ et distributions d'ordre $r$}\label{analyse}

Le but de ce paragraphe est de rappeler les d\'efinitions et \'enonc\'es (classiques) d'analyse $p$-adique utilis\'es dans la suite. Pour les preuves, nous renvoyons \`a \cite{Co5}.

Si $f:\Z\to L$ est une fonction quelconque, 
on pose pour $n$ entier positif ou nul :
$$a_n(f)\=\sum_{i=0}^n(-1)^i\binom{n}{i}f(n-i).$$

Soit $r$ un nombre r\'eel positif ou nul.

\begin{defi}\label{classecr} 
Une fonction $f : \Z \to L$ est de classe $\mathcal{C}^r$ si
$n^r |a_n(f)| \to 0$ dans $\RR_{\geq 0}$ quand $n \to + \infty$. 
\end{defi} 

Rappelons que $f$ est continue si et seulement si $a_n(f)$ tend $p$-adiquement vers $0$ quand $n$ tend vers l'infini, de sorte que les fonctions de classe $\mathcal{C}^0$ au sens de la d\'efinition \ref{classecr} sont pr\'ecis\'ement les fonctions continues sur $\Z$. Toute fonction de classe $\mathcal{C}^r$ est aussi de classe $\mathcal{C}^s$ pour $0\leq s\leq r$ et est donc en particulier toujours continue. On note $\mathcal{C}^r(\Z,L)$ le $L$-espace vectoriel des fonctions $f:\Z\to L$ de classe $\mathcal{C}^r$.

Toute fonction continue $f : \Z \to L$ s'\'ecrit (d\'eveloppement de Mahler) :
\begin{equation}\label{mahler}
f(z)=\sum_{n=0}^{+\infty}a_n(f)\binom{z}{n} 
\end{equation} 
o\`u $z \in \Z$ et $\binom{z}{0}\=1$, $\binom{z}{n} = z(z-1)\cdots (z-n+1) / n!$ si $n\geq 1$. De plus, $\|f\|_0\={\rm sup}_{z\in\Z}|f(z)|$ co\"\i ncide avec ${\rm sup}_{n\geq 0}|a_n(f)|$. On v\'erifie facilement que l'espace $\mathcal{C}^r(\Z,L)$ est un espace de Banach pour la norme $\|f\|_r\={\rm sup}_{n\geq 0}(n+1)^r|a_n(f)|$.

La terminologie {\og de classe $\mathcal{C}^r$ \fg} provient du fait que, lorsque $r$ est entier strictement positif, les fonctions de $\mathcal{C}^r(\Z,L)$ sont aussi les fonctions sur $\Z$ qui, moralement, admettent $r$ d\'eriv\'ees avec la $r$-i\`eme d\'eriv\'ee continue, ce qui explique la terminologie. Par exemple, les fonctions localement analytiques sur $\Z$ sont de classe $\mathcal{C}^r$ pour tout $r\in \R_{\geq 0}$.

\begin{theo}[Amice-V\'elu, Vishik]\label{amice}
Soit $d$ un entier tel que $r-1<d$. Le sous-$L$-espace vectoriel ${\rm Pol}^d(\Z,L)$ de $\mathcal{C}^r(\Z,L)$ des fonctions $f:\Z\to L$ localement polynomiales de degr\'e (local) au plus $d$ est dense dans $\mathcal{C}^r(\Z,L)$.
\end{theo}

\begin{defi}
On appelle distribution temp\'er\'ee d'ordre $r$ un \'el\'ement de
l'espace $\mathcal{C}^r(\Z,L)^*$, c'est-\`a-dire une forme lin\'eaire continue sur l'espace de Banach des fonctions de classe $\mathcal{C}^r$.
\end{defi}

On dit parfois aussi distribution temp\'er\'ee d'ordre $\leq r$. Nous donnons maintenant deux descriptions des distributions temp\'er\'ees d'ordre $r$. 

Par le th\'eor\`eme \ref{amice}, l'inclusion ${\rm Pol}^d(\Z,L)\subsetneq \mathcal{C}^r(\Z,L)$ induit lorsque $r-1<d$ une injection :
$$\mathcal{C}^r(\Z,L)^*\hookrightarrow {\rm Pol}^d(\Z,L)^*$$
o\`u ${\rm Pol}^d(\Z,L)^*$ est l'espace vectoriel dual de ${\rm Pol}^d(\Z,L)$. Rappelons que si $U$ est un ouvert de $\Zp$ et $\mathbf{1}_U$ sa fonction caract\'eristique, alors on note $\int_U f(z) d \mu(z)$ pour $\mu(\mathbf{1}_U(z)f(z))$.
 
\begin{theo}[Amice-V\'elu, Vishik]\label{amice2}
Soit $\mu\in {\rm Pol}^d(\Z,L)^*$ et supposons que $r-1<d$. Alors $\mu\in \mathcal{C}^r(\Z,L)^*$ si et seulement s'il existe une constante $C_{\mu}\in L$ telle que, $\forall\ a\in \Z$, $\forall\ j\in \{0,\cdots,d\}$ et $\forall\ n\in \NN$ :  
\begin{equation}\label{borneamice}
\int_{a+p^n\Z}(z-a)^jd\mu(z)\in C_{\mu}p^{n(j-r)}\O.
\end{equation}
\end{theo}

Remarquons que le plus petit entier $d$ tel que le th\'eor\`eme \ref{amice2} s'applique est la partie enti\`ere de $r$. Lorsque $\mu$ est d'ordre $r$ et $r-1<d$, on pose :
\begin{equation}\label{normedistrd}
\|\mu\|_{r,d}\={\rm sup}_{a\in\Z}{\rm sup}_{j\in \{0,\cdots,d\}}{\rm sup}_{n\in\NN}\left(p^{n(j-r)}\left|\int_{a+p^n\Z}(z-a)^jd\mu(z)\right|\right).
\end{equation}
On peut montrer que $\|\mu\|_{r,d}$ est une norme sur $\mathcal{C}^r(\Z,L)^*$ qui redonne la topologie d'espace de Banach de $\mathcal{C}^r(\Z,L)^*$ et qui est \'equivalente \`a la norme :
\begin{equation}\label{normedistr}
\|\mu\|_{r}\={\rm sup}_{a\in\Z}{\rm sup}_{j,n\in\NN}\left(p^{n(j-r)}\left|\int_{a+p^n\Z}(z-a)^jd\mu(z)\right|\right).
\end{equation}
En particulier, la majoration (\ref{borneamice}) est \'equivalente \`a la m\^eme majoration pour tout $j\in \N$ (et tout $a\in \Z$, $n\in \N$), quitte peut-\^etre \`a modifier $C_{\mu}$.

Nous aurons besoin du lemme suivant :

\begin{lemm}\label{lemmech}
Soit $r\in \R_{\geq 0}$ et $d$ la partie enti\`ere de $r$. Soit $n\in \N$ et :
$$f(z)\=\sum_{a\in \{0,\hdots,p^n-1\}}{\bf 1}_{a+p^n\Z}(z)\sum_{i=0}^d\lambda_{a,i}(z-a)^i\in {\rm Pol}^d(\Z,L)$$
o\`u $\lambda_{a,i}\in L$. Alors :
\begin{equation}\label{arghh}
\sup_{\mu\in \mathcal{C}^r(\Z,L)^*}\frac{\left|\int_{\Z}f(z)d\mu(z)\right|}{\|\mu\|_{r,d}}=\sup_{a\in \{0,\hdots,p^n-1\}}\sup_{i\in \{0,\hdots,d\}}\left|\lambda_{a,i}\right|p^{n(r-i)}.
\end{equation}
\end{lemm}
\begin{proof}
Par (\ref{normedistrd}), on voit que le r\'eel de gauche est plus petit que celui de droite. Pour $(a,i)\in \{0,\hdots,p^n-1\}\times \{0,\hdots,d\}$, il n'est pas difficile de construire une forme lin\'eaire $\mu_{a,i}\in {\rm Pol}^d(\Z,L)^*$ \`a support dans $a+p^n\Z$ satisfaisant (\ref{borneamice}) telle que $\int_{a+p^n\Z}(z-a)^jd\mu_{a,i}(z)=0$ si $j\ne i$ ($j\in \{0,\hdots, d\}$), $\int_{a+p^n\Z}(z-a)^id\mu_{a,i}(z)=p^{n(i-r)}$ et $\|\mu_{a,i}\|_{r,d}=1$ (les d\'etails sont laiss\'es en exercice au lecteur). En particulier : $$\frac{\left|\int_{a+p^n\Z}f(z)d\mu_{a,i}(z)\right|}{\|\mu_{a,i}\|_{r,d}}=\left|\lambda_{a,i}\right|p^{n(r-i)}$$ 
est inf\'erieur au r\'eel de gauche. Comme cela est vrai pour tout $(a,i)\in \{0,\hdots,p^n-1\}\times \{0,\hdots,d\}$, on en d\'eduit le r\'esultat.
\end{proof}

L'espace vectoriel ${\rm An}(\Z,L)$ des fonctions localement analytiques sur $\Z$ (muni de sa topologie d'espace de type compact, cf. \cite[\S16]{Sc} ou \cite[\S1.4.3]{Co5}), est {\it a fortiori} dense dans $\mathcal{C}^r(\Z,L)$ et on dispose donc aussi d'une injection continue entre duaux continus :
$$\mathcal{C}^r(\Z,L)^*\hookrightarrow {\rm An}(\Z,L)^*.$$
La transform\'ee d'Amice-Mahler :
\begin{equation}\label{Yvette} 
\mu\longmapsto \sum_{n=0}^{+\infty}\mu\left(\binom{z}{n}\right)X^n 
\end{equation} 
induit un isomorphisme topologique entre ${\rm An}(\Z,L)^*$ et
$\calR^+$. Rappelons que (voir \S\ref{series}) : 
$$\calR^+\=\left\{\sum_{n=0}^{+\infty}a_nX^n\mid a_n\in L,\ {\rm lim}_{n \to \infty} |a_n|\rho^n=0\ \forall\ \rho\in [0,1[\right\}$$
et que cet espace est
muni de la topologie d'espace de Fr\'echet induite par la collection
des normes $\|\cdot \|_{D(0,\rho)}={\rm sup}_{n \geq 0}(|a_n|\rho^n)$ pour $0<\rho<1$. Rappelons \'egalement (d\'efinition \ref{dordrer}) qu'un \'el\'ement $w\in \calR^+$ est d'ordre $r$ si, pour un (ou de mani\`ere \'equivalente tous les) $\rho\in ]0,1[$, la suite $\left(p^{-nr}\|w\|_{D(0,\rho^{1/p^n})}\right)_{n \geq 0}$ est born\'ee dans $\R_{\geq 0}$.

\begin{lemm}\label{ordre}
Soit $w=\sum_{n=0}^{+\infty}a_nX^n\in \calR^+$.
\begin{itemize}
\item[(i)] Un \'el\'ement $w=\sum_{n=0}^{+\infty}a_nX^n\in \calR^+$ est
d'ordre $r$ si et seulement si $\{n^{-r}|a_n|\}_{n \geq 1}$ 
est born\'e (dans $\R_{\geq 0}$) lorsque $n$ varie.
\item[(ii)] Les normes ${\rm sup}_{n \geq 0}\left(p^{-nr}\|w\|_{D(0,\rho^{1/p^n})}\right)$ et ${\rm sup}_{n \geq 0} \left((n+1)^{-r}|a_n|\right)$ sont \'equivalentes pour tout $\rho\in ]0,1[$.
\end{itemize}
\end{lemm}
\begin{proof}
Voir \cite[\S2.1]{Co5}.
\end{proof}

Le r\'esultat suivant d\'ecoule imm\'ediatement des d\'efinitions et du lemme \ref{ordre} :  

\begin{prop}\label{amice3}
Soit $\mu\in {\rm An}(\Z,L)^*$. Alors $\mu\in \mathcal{C}^r(\Z,L)^*$ si et seulement l'\'el\'ement $\sum_{n=0}^{+\infty}\mu\left(\binom{z}{n}\right)X^n$ est d'ordre $r$ dans $\calR^+$.
\end{prop}

On peut montrer que $\|\mu\|_{r}'\={\rm sup}_n((n+1)^{-r}\left|\mu\left(\binom{z}{n}\right)\right|)$ est une norme sur $\mathcal{C}^r(\Z,L)^*$ qui redonne la topologie d'espace de Banach de $\mathcal{C}^r(\Z,L)^*$ (\cite[\S2.3.1]{Co5}).

\subsection{D\'efinition de $\breuil(V)$}\label{definition1}

Le but de ce paragraphe est de donner une premi\`ere d\'efinition de $\breuil(V)$ comme espace fonctionnel.

Soit $V$ une repr\'esentation apc absolument irr\'eductible comme au \S\ref{cristallin} (rappelons que $V$ \`a pour poids de Hodge-Tate $(0,k-1)$ avec n\'ecessairement $k\geq 2$ sinon $V$ n'est pas absolument irr\'eductible). Les valeurs propres $\alpha_p^{-1}$ et $\beta_p^{-1}$ du semi-simplifi\'e de $\varphi$ sur $\dcris(V)=D(\alpha_p,\beta_p)$ sont alors telles que ${\rm val}(\alpha_p)>0$, ${\rm val}(\beta_p)>0$ et ${\rm val}(\alpha_p)+{\rm val}(\beta_p)=k-1$. Remarquons que $\varphi$ est semi-simple si et seulement si $\alpha\ne \beta$. Quitte \`a \'echanger $\alpha$ et $\beta$, on suppose dans la suite ${\rm val}(\alpha_p)\geq {\rm val}(\beta_p)$.

Soit $B(\alpha)$ l'espace de Banach suivant. Son $L$-espace vectoriel sous-jacent est form\'e des fonctions $f:\Q\to L$ v\'erifiant les deux conditions :
\begin{itemize}
\item[(i)] $f_{\mid \Z}$ est une fonction de classe $\mathcal{C}^{{\rm val}(\alpha_p)}$;
\item[(ii)] $(\beta\alpha^{-1})(z)|z|^{-1} z^{k-2}f(1/z) _{\mid \Z-\{0\}}$ se prolonge sur $\Z$ en une fonction de classe $\mathcal{C}^{{\rm val}(\alpha_p)}$.
\end{itemize}
Comme espace vectoriel, on a donc :
\begin{equation}\label{f1f2} 
B(\alpha)\simeq \mathcal{C}^{{\rm val}(\alpha_p)}(\Z,L)\oplus 
\mathcal{C}^{{\rm val}(\alpha_p)}(\Z,L),\ f\mapsto f_1\oplus f_2 
\end{equation} 
o\`u, pour $z\in \Z$, $f_1(z)\=f(pz)$ et $f_2(z)\=(\beta\alpha^{-1})(z)|z|^{-1} z^{k-2} f(1/z)$. On munit $B(\alpha)$ de la norme :  
$$\|f\|\={\rm max}\left(\|f_1\|_{{\rm val}(\alpha_p)},\|f_2\|_{{\rm val}(\alpha_p)}\right),$$
qui en fait un espace de Banach en vertu du \S\ref{analyse}. On fait agir $L$-lin\'eairement $\G$ \`a gauche sur les fonctions de $B(\alpha)$ comme suit : 
\begin{equation}\label{action} 
\begin{pmatrix}a&b\\c&d\end{pmatrix} \cdot f(z)
=\alpha(ad-bc)(\beta\alpha^{-1})(-cz+a)|-cz+a|^{-1}(-cz+a)^{k-2}
 f\left(\frac{dz-b}{-cz+a}\right). 
\end{equation} 
Notons que $\left(\begin{smallmatrix}a &0\\0&a\end{smallmatrix}\right)$ agit par la multiplication par $\varepsilon^{k-2}(a){(\alpha\beta)(a)|a |^{-(k-1)}}\in \O^{\times}$.

\begin{lemm}\label{action2}
Si $f\in B(\alpha)$ et $g\in \G$, alors $g\cdot f\in B(\alpha)$ et l'action de $\G$ sur $B(\alpha)$ se fait par automorphismes continus.
\end{lemm}
\begin{proof}
On pose $r\= {\rm val}(\alpha_p)$, $d$ la partie enti\`ere de $r$ (donc $d\leq k-2$) et on munit l'espace de Banach $\mathcal{C}^r(\Z,L)$ de la norme induite par son bidual, c'est-\`a-dire :
$$\|f\|\=\sup_{\mu\in \mathcal{C}^r(\Z,L)^*}\frac{\left|\int_{\Z}f(z)d\mu(z)\right|}{\|\mu\|_{r,d}}$$
qui redonne la topologie d'espace de Banach de $\mathcal{C}^r(\Z,L)$ par \cite[lemme 9.9]{Sc} et les r\'esultats du \S\ref{analyse}. L'assertion du lemme est triviale si $g$ est scalaire. Par la d\'ecomposition de Bruhat et le cas scalaire, on est r\'eduit \`a montrer la stabilit\'e et la continuit\'e pour les matrices $g$ de la forme $\left(\begin{smallmatrix}0 &p\\1&0\end{smallmatrix}\right)$, $\left(\begin{smallmatrix}1 &0\\0&\lambda\end{smallmatrix}\right)$ et $\left(\begin{smallmatrix}1 &\lambda\\0&1\end{smallmatrix}\right)$ (avec $\lambda\in \Q^{\times}$). Le premier cas est \'evident puisqu'il envoie $f=(f_1,f_2)\in B(\alpha)$ sur $(f_2,f_1)\in B(\alpha)$ \`a multiplication pr\`es par des scalaires (cf. (\ref{f1f2})). Quitte \`a conjuguer par $\left(\begin{smallmatrix} 0 &p\\1&0\end{smallmatrix}\right)$ et \`a multiplier par un scalaire convenable, on peut prendre $\lambda\in \Z-\{0\}$ dans le deuxi\`eme et $g$ envoie alors $(f_1(z),f_2(z))\in B(\alpha)$ sur $(f_1(\lambda z),f_2(\lambda z))$ (\`a multiplication pr\`es par des scalaires). Il suffit donc de montrer que l'application $f(\cdot)\mapsto f(\lambda\cdot)$ est bien d\'efinie et continue de $\mathcal{C}^r(\Z,L)$ dans $\mathcal{C}^r(\Z,L)$. Elle est bien d\'efinie par \cite[prop 1.38]{Co5}. Pour la continuit\'e, par le th\'eor\`eme \ref{amice}, il suffit de montrer qu'il existe $c\in |L^{\times}|$ tel que, si $f\in {\rm Pol}^d(\Z,L)\subset \mathcal{C}^r(\Z,L)$, alors $\|f(\lambda\cdot)\|\leq c\|f(\cdot)\|$. En \'ecrivant, pour un $n\in \N$ assez petit et des $\lambda_{n,a,i}\in L$ convenables :
$$f(z)=\sum_{a\in \{0,\hdots,p^n-1\}}{\bf 1}_{a+p^n\Z}(z)\sum_{i=0}^d\lambda_{a,i}(z-a)^i,$$
un calcul donne :
$$f(\lambda z)=\sum_{\substack{a\in \{0,\hdots,p^n-1\} \\
{\rm val}(a)\geq {\rm val}(\lambda)}}
{\bf 1}_{\frac{a}{\lambda}+p^{n-{\rm val}(\lambda)}\Z}(z)\sum_{i=0}^d\lambda_{a,i}\lambda^i\left(z-\frac{a}{\lambda}\right)^i.$$
On en d\'eduit gr\^ace \`a (\ref{normedistrd}) :
\begin{eqnarray*}
\|f(\lambda\cdot)\|&\leq &\sup_{\substack{a\in \{ 0,\hdots,p^n-1 \} \\ {\rm val}(a)\geq {\rm val}(\lambda)}}\sup_{i\in \{0,\hdots,d\}}\left|\lambda_{a,i}\lambda^i\right| p^{(n-{\rm val}(\lambda))(r-i)}\\
&\leq &\sup_{a\in \{0,\hdots,p^n-1\}}\sup_{i\in \{0,\hdots,d\}}\left|\lambda_{a,i}\right| p^{n(r-i)}p^{-r{\rm val}(\lambda)}\\
&\leq & \|f(\cdot)\|
\end{eqnarray*}
o\`u la derni\`ere in\'egalit\'e r\'esulte du lemme \ref{lemmech} et du fait que ${\rm val}(\lambda)\geq 0$. Passons au dernier cas. Quitte \`a conjuguer $g= \left(\begin{smallmatrix}1 &\lambda\\0&1\end{smallmatrix}\right)$ par un \'el\'ement convenable de $\left(\begin{smallmatrix} 1 &0\\0&\Q^{\times}\end{smallmatrix}\right)$, on peut supposer $\lambda=p$ et $g$ envoie alors $(f_1,f_2)$ sur $(f_1(z+1),(1+pz)^{k-2}f_2(z/(1+pz)))$ (\`a multiplication pr\`es par des scalaires). Il suffit donc de montrer que les applications $f(\cdot)\mapsto f(\cdot+1)$ et $f(\cdot)\mapsto (1+pz)^{k-2}f(\cdot /(1+p\cdot ))$ sont bien d\'efinies et continues de $\mathcal{C}^r(\Z,L)$ dans $\mathcal{C}^r(\Z,L)$.  Elles sont bien d\'efinies encore par \cite[\S1.5]{Co5}. La continuit\'e se v\'erifie par un argument analogue au pr\'ec\'edent avec un calcul utilisant (\ref{normedistrd}) et (\ref{normedistr}) dont on laisse les d\'etails au lecteur.
\end{proof}

La repr\'esentation $B(\alpha)$ doit \^etre pens\'ee comme une induite parabolique. Sans donner un sens formel \`a ce qui suit, on a un isomorphisme $\G$-\'equivariant :
$$\left({\rm Ind}_{\B}^{\G}\alpha\otimes d^{k-2}\beta\nrm^{-1} \right)^{\mathcal{C}^{{\rm val}(\alpha_p)}}\simeq B(\alpha)$$
o\`u l'espace de gauche est celui des fonctions $F:\G\to L$ qui sont de classe $\mathcal{C}^{{\rm val}(\alpha_p)}$ (oublions que nous n'avons pas d\'efini de telles fonctions dans ce cadre !) telles que :
\begin{equation}\label{equfonc}
F\left( \begin{pmatrix}a&b\\0&d\end{pmatrix}g\right) = \alpha(a)\beta(d)| d|^{-1} d^{k-2}F(g)
\end{equation}
avec action de $\G$ donn\'ee par $(g\cdot F)(g')\=F(g'g)$. On passe de $F$ \`a une fonction $f\in B(\alpha)$ en posant :
\begin{equation}\label{fonction} 
f(z)\=F\left(\begin{pmatrix}0&1\\-1&z\end{pmatrix}\right).
\end{equation}
On d\'efinit de m\^eme $B(\beta)$ muni d'une action de $\G$ par automorphismes en \'echan\-geant partout $\alpha$ et $\beta$. 

Voici des exemples importants de fonctions dans $B(\alpha)$ :

\begin{lemm}\label{appartient}
Pour $0\leq j<{\rm val}(\alpha_p)$ et $a\in \Q$, les fonctions $z\mapsto z^j$ et les fonctions :
$$z\longmapsto (\beta\alpha^{-1})(z-a)| z-a|^{-1}(z-a)^{k-2-j}$$
(prolong\'ees par $0$ en $a$) sont dans $B(\alpha)$. 
\end{lemm} 
\begin{proof} 
En faisant agir $\left(\begin{smallmatrix}0 &1\\1&0\end{smallmatrix}\right)$ sur $z^j$, il suffit de traiter les deuxi\`emes fonctions. En faisant agir $\left(\begin{smallmatrix}1 &a\\0&1\end{smallmatrix}\right)$, il suffit m\^eme par le lemme \ref{action2} de traiter le cas $a=0$ et comme $z\mapsto z^j$ est clairement de classe $\mathcal{C}^{{\rm val}(\alpha_p)}$ sur $\Z$, il suffit de montrer que $z\mapsto f(z)\=(\beta\alpha^{-1})(z)|z|^{-1} z^{k-2-j}$ est de classe $\mathcal{C}^{{\rm val}(\alpha_p)}$ sur $\Z$. Soit $f_0$ la fonction nulle sur $\Z$ et, pour $n\in \ZZ$, $n>0$, posons $f_n(z)\=(\beta\alpha^{-1})(z)| z |^{-1} z^{k-2-j}$ si ${\rm val}(z)<n$, $f_n(z)\=0$ sinon. La fonction $f_n$ est de classe $\mathcal{C}^{{\rm val}(\alpha_p)}$ sur $\Z$ puisqu'elle est localement polynomiale. Il suffit de montrer que $f_{n+1}-f_n\to 0$ dans $\mathcal{C}^{{\rm val}(\alpha_p)}(\Z,L)$ quand $n\to +\infty$ (car $\sum_{n=0}^{\infty}(f_{n+1}-f_n)=f\in \mathcal{C}^{{\rm val}(\alpha_p)}(\Z,L)$ puisque $\mathcal{C}^{{\rm val}(\alpha_p)}(\Z,L)$ est complet). Par \cite[lemme 9.9]{Sc}, il suffit de v\'erifier que $f_{n+1}-f_n\to 0$ dans $(B(\alpha)^*)^*$, i.e. que :  
$$\sup_{\mu\in \mathcal{C}^{{\rm val}(\alpha_p)}(\Z,L)^*}\frac{\left|\int_{\Z}(f_{n+1}(z)-f_n(z))d\mu(z)\right|}{\|\mu\|_{{\rm val}(\alpha_p)}}\longrightarrow 0\ \text{quand $n\to +\infty$.}$$
Si $S\subset \Z^{\times}$ est un syst\`eme de repr\'esentants des classes de $(\ZZ/p^{m(V)}\ZZ)^{\times}$, alors :
$$\int_{\Z}(f_{n+1}(z)-f_n(z))d\mu(z)=\left(\frac{\alpha_p p}{\beta_p}\right)^{\!n}\sum_{a_i\in S}\beta\alpha^{-1}(a_i)\int_{p^na_i+p^{n+m(V)}\Z}z^{k-2-j}d\mu(z).$$ 
En \'ecrivant $z^{k-2-j}=(z-p^na_i+p^na_i)^{k-2-j}$ et en d\'eveloppant, on obtient :
\begin{multline*}
\left|\int_{\Z}(f_{n+1}(z)-f_n(z))d\mu(z)\right| \\
\leq  p^{-n(2{\rm val}(\alpha_p)-k+2)}\cdot
\sup_{\substack{a_i\in S \\ 0\leq j'\leq k-2-j}}p^{-nj'}\left|\int_{p^na_i+p^{n+m(V)}\Z}(z-p^na_i)^{k-2-j-j'}d\mu(z)\right|,
\end{multline*}
soit, en utilisant (\ref{normedistr}) :
\begin{multline*}
\left|\int_{\Z}(f_{n+1}(z)-f_n(z))d\mu(z)\right| \\ 
\leq \|\mu\|_{{\rm val}(\alpha_p)}p^{-n(2{\rm val}(\alpha_p)-k+2)}\cdot 
\sup_{0\leq j'\leq k-2-j}p^{-nj'}p^{-(n+m(V))(k-2-j-j'-{\rm val}(\alpha_p))},
\end{multline*}
soit encore :
$$\left|\int_{\Z}(f_{n+1}(z)-f_n(z))d\mu(z)\right|\leq C\|\mu\|_{{\rm val}(\alpha_p)}p^{n(j-{\rm val}(\alpha_p))},$$
o\`u $C\=\sup_{0\leq j'\leq k-2-j}p^{-m(V)(k-2-j-j'-{\rm val}(\alpha_p))}$. D'o\`u le r\'esultat puisque $j<{\rm val}(\alpha_p)$.  
\end{proof} 

On a un lemme analogue en \'echangeant $\alpha$ et $\beta$. On note $L(\alpha)$ l'adh\'erence dans $B(\alpha)$ du sous-$L$-espace vectoriel engendr\'e par les fonctions $z^j$ et $(\beta\alpha^{-1})(z-a)| z-a |^{-1} (z-a)^{k-2-j}$ pour $a\in \Q$ et $j$ entier, $0\leq j<{\rm val}(\alpha_p)$. De m\^eme, on note $L(\beta)$ l'adh\'erence dans $B(\beta)$ du sous-$L$-espace vectoriel engendr\'e par les fonctions $z^j$ et $(\alpha\beta^{-1})(z-a)|z-a |^{-1} (z-a)^{k-2-j}$ pour $a\in \Q$ et $j$ entier, $0\leq j<{\rm val}(\beta_p)$. Notons que, lorsque $\alpha = \beta\nrm$, $L(\beta)$ est de dimension finie et s'identifie aux polyn\^omes en $z$ de degr\'e au plus $k-2$.

\begin{lemm}
Le sous-espace $L(\alpha)$ (resp. $L(\beta)$) est stable par $\G$ dans $B(\alpha)$ (resp. $B(\beta)$).
\end{lemm}
\begin{proof} 
Exercice.
\end{proof} 

\begin{defi}
On pose $\breuil(V)\=B(\alpha)/L(\alpha)$.
\end{defi}

Il s'agit encore d'un $L$-espace de Banach (avec la topologie quotient) muni d'une action de $\G$ par automorphismes continus. Nous allons voir que l'application $\G\times \breuil(V)\to \breuil(V)$ est continue, que $\breuil(V)$ est unitaire et que l'on a un morphisme continu $\G$-\'equivariant $\widehat I:B(\beta)/L(\beta)\to \breuil(V)$ qui est un isomorphisme si $\alpha\ne \beta\nrm$. En particulier, lorsque $\val(\alpha)=\val(\beta)$, les $\G$-repr\'esentations $B(\alpha)/L(\alpha)$ et $B(\beta)/L(\beta)$ sont topologiquement isomorphes et il n'y a donc pas d'ambiguit\'e dans ce cas sur la d\'efinition de $\breuil(V)$.

\subsection{Une autre description de $\breuil(V)$}\label{definition2}

Le but de ce paragraphe est de donner une description plus intrins\`eque de $B(\alpha)/L(\alpha)$ pour en d\'eduire certaines propri\'et\'es de l'action de $\G$ (continuit\'e, unitarit\'e, entrelacements) peu \'evidentes sur la d\'efinition pr\'ec\'edente. On conserve les notations du \S\ref{definition1}.

Soit :
$$\pi(\alpha)\={\rm Sym}^{k-2}L^2\otimes_L{\rm Ind}_{\B}^{\G}\alpha\otimes \beta\nrm^{-1}$$ 
la repr\'esentation de $\G$ produit tensoriel de la repr\'esentation
alg\'ebrique ${\rm Sym}^{k-2}L^2$ par l'induite parabolique lisse
${\rm Ind}_{\B}^{\G}\alpha\otimes \beta\nrm^{-1}$ (c'est-\`a-dire l'espace des fonctions localement
constantes $h:\G\to L$ v\'erifiant une \'egalit\'e analogue
\`a (\ref{equfonc}) avec action \`a gauche de $\G$ par translation \`a
droite). On munit $\pi(\alpha)$ de l'unique topologie localement
convexe (au sens de \cite{Sc}) telle que les ouverts sont les
sous-$\O$-modules g\'en\'erateurs (sur $L$). La repr\'esentation
$\pi(\alpha)$ est dite localement alg\'ebrique (cf. l'appendice de
\cite{ST2}) et n'est autre que la repr\'esentation ${\rm Alg}(V)
\otimes_L {\rm Lisse}(V)$ de l'introduction.

On identifie ${\rm Sym}^{k-2}L^2$ \`a l'espace vectoriel des polyn\^omes $P(z)$ de degr\'e au plus $k-2$ \`a coefficients dans $L$ munis de l'action \`a gauche de $\G$ :
\begin{equation}\label{algebraic}
\begin{pmatrix}a&b\\c&d\end{pmatrix}\cdot P(z)
=(-cz+a)^{k-2}P\left(\frac{dz-b}{-cz+a}\right).
\end{equation}
Comme en (\ref{fonction}), on identifie ${\rm Ind}_{\B}^{\G}\alpha\otimes d^{k-2}\beta\nrm^{-1}$ \`a l'espace vectoriel des fonctions $f:\Q\to L$ localement constantes telles que $(\beta\alpha^{-1})(z)|z|^{-1} f(1/z)$ se prolonge sur $\Q$ en une fonction localement constante avec action \`a gauche de $\G$ comme en (\ref{action}) mais sans le facteur $(-cz+a)^{k-2}$. On en d\'eduit une injection $\G$-\'equivariante continue :
\begin{equation}\label{commefonction}
\pi(\alpha)\hookrightarrow B(\alpha),\ P(z)\otimes f(z)\mapsto P(z)f(z).
\end{equation} 
Par le th\'eor\`eme \ref{amice}, l'image de $\pi(\alpha)$ est dense dans $B(\alpha)$. En particulier, on a une injection continue $B(\alpha)^*\hookrightarrow \pi(\alpha)^*$. On d\'efinit de m\^eme $\pi(\beta)$ et une injection $\G$-\'equivariante continue d'image dense $\pi(\beta)\hookrightarrow B(\beta)$. Ces injections induisent des applications $\G$-\'equivariantes continues $\pi(\alpha)\to B(\alpha)/L(\alpha)$ et $\pi(\beta)\to B(\beta)/L(\beta)$.

Si $\pi^0$ est un sous-$\O$-module g\'en\'erateur d'un $L$-espace vectoriel $\pi$, rappelons qu'on appelle compl\'et\'e de $\pi$ par rapport \`a $\pi^0$ l'espace de Banach :
$$B\= (\varprojlim_n \pi^0/\pi_L^n\pi^0)\otimes_{\O}L.$$
On a un morphisme canonique d'image dense $\pi\to B$ qui n'est pas injectif en g\'en\'eral (si $\pi^0=\pi$, on a $B=0$). Le dual continu $B^*$ de $B$ s'identifie en tant qu'espace de Banach \`a ${\rm Hom}_{\O}(\pi^0,\O)\otimes_{\O}L$ (avec ${\rm Hom}_{\O}(\pi^0,\O)$ comme boule unit\'e). Si $\pi$ est un espace localement convexe tonnel\'e (cf. \cite[\S6]{Sc}) muni d'une action continue d'un groupe topologique localement compact $G$ telle que $\pi^0$ est ouvert et stable par $G$, il est facile de v\'erifier en utilisant le th\'eor\`eme de Banach-Steinhaus (cf. \cite[proposition 6.15]{Sc}) que $B$ et $B^*$ sont des $G$-Banach unitaires et que la fl\`eche canonique $\pi\to B$ est continue.

\begin{theo}\label{complete}
L'application $\pi(\alpha)\to B(\alpha)/L(\alpha)$ induit un isomorphisme topologique $\G$-\'equivariant entre $B(\alpha)/L(\alpha)$ et le compl\'et\'e de $\pi(\alpha)$ par rapport \`a un quelconque sous-$\O$-module g\'en\'erateur de $\pi(\alpha)$ stable par $\G$ et de type fini comme $\O[\G]$-module. On a le m\^eme r\'esultat en rempla\c cant $\alpha$ par $\beta$.
\end{theo}
\begin{proof}
Notons que le compl\'et\'e ne d\'epend pas du choix du sous-$\O[\G]$-module de type fini g\'en\'erateur de $\pi(\alpha)$ car ces $\O$-modules sont tous commensurables dans $\pi(\alpha)$. En utilisant $\G=\B\K$ et le fait que $\K$ est compact, on voit facilement qu'il suffit de compl\'eter par rapport \`a un sous-$\O[\B]$-module de type fini g\'en\'erateur quelconque, par exemple :
$$\sum_{j=0}^{k-2}\O[\B]({\mathbf 1}_{\Z}(z)z^j)+\sum_{j=0}^{k-2}
\O[\B]\left((\beta\alpha^{-1})(z)| z|^{-1}{\mathbf 1}_{\Q-\Z}(z)z^j\right)\subset \pi(\alpha)$$ 
o\`u ${\bf 1}_U$ est la fonction caract\'eristique de l'ouvert $U$. Le dual du compl\'et\'e cherch\'e est donc isomorphe au Banach :
\begin{multline}\label{thedual} 
\{\mu\in \pi(\alpha)^*\mid\forall\ g\in\B, \forall\ j\in
\{0,\hdots,k-2\}, |\mu(g({\mathbf 1}_{\Z}(z)z^j))|\leq 1\\   
 {\rm et}\ \left|\mu(g({\mathbf 1}_{\Q-\Z}(z)(\beta\alpha^{-1})(z)| z|^{-1}z^j))\right|\leq 1\}\otimes_{\O}L.  
\end{multline}
En utilisant l'int\'egralit\'e du caract\`ere central, il est \'equivalent de prendre $g\in \left(\begin{smallmatrix}1&\Q \\0&\Q^{\times}\end{smallmatrix} \right)$ dans (\ref{thedual}). Pour $f\in \pi(\alpha)$, vue comme fonction sur $\Q$ via (\ref{commefonction}), et $U$ ouvert de $\Q$, on \'ecrit $\int_{U}f(z)d\mu(z)$ pour $\mu({\bf 1}_U(z)f(z))$. Un calcul donne alors que les conditions sur $\mu$ dans (\ref{thedual}) sont \'equivalentes \`a l'existence d'une constante $C\in L$ ind\'ependante de $\mu$ telle que, pour tout $a\in \Q$, tout $j\in \{0,\hdots,k-2\}$ et tout $n\in \ZZ$ :  
\begin{eqnarray}\label{chaud1} 
\int_{a+p^n\Z}(z-a)^jd\mu(z)&\in & Cp^{n(j-{\rm val}(\alpha_p))}\O\\  
\label{chaud2}\int_{\Q-(a+p^n\Z)}(\beta\alpha^{-1})(z-a)| z-a|^{-1}(z-a)^{k-2-j}d\mu(z)&\in & Cp^{n({\rm val}(\alpha_p)-j)}\O
\end{eqnarray}
(si $g=\left(\begin{smallmatrix}1&\lambda \\0&\mu\end{smallmatrix} \right)$, poser $n=-{\rm val}(\mu)$ et $a=\lambda/\mu$; en fait, on peut prendre $C=1$). En \'ecrivant $(a+p^{n-1}\Z)-(a+p^n\Z)=\cup_{a_i\in S}a+p^{n-1}a_i+p^{n-1+m(V)}\Z$ et en d\'eveloppant $(z-a)^{k-2-j}=((z-a-p^{n-1}a_i)+p^{n-1}a_i)^{k-2-j}$ comme dans la preuve du lemme \ref{appartient}, on d\'eduit facilement de (\ref{chaud1}) quitte \`a modifier $C$ :
\begin{equation}\label{chaud3} 
\int_{(a+p^{n-1}\Z)-(a+p^n\Z)} (\beta\alpha^{-1})(z-a)| z-a|^{-1}(z-a)^{k-2-j}d\mu(z)\in Cp^{n({\rm val}(\alpha_p)-j)}\O.
\end{equation}
En d\'ecompo\-sant $\Q-(a+p^n\Z)=\Q-(a+p^{n+1}\Z)\setminus (a+p^{n}\Z)-(a+p^{n+1}\Z)$, puis $\Q-(a+p^{n+1}\Z)=\Q-(a+p^{n+2}\Z)\setminus (a+p^{n+1}\Z)-(a+p^{n+2}\Z)$ etc. jusqu'\`a arriver \`a $\Q-(a+p^{n+m}\Z)$ avec $n+m\geq 0$, on d\'eduit de (\ref{chaud3}) que (\ref{chaud2}) pour $j<{\rm val}(\alpha_p)$ d\'ecoule de (\ref{chaud1}) et de (\ref{chaud2}) pour $n\geq 0$ (utiliser (\ref{chaud3}) pour les morceaux compacts dans la d\'ecomposition et (\ref{chaud2}) avec $n'=n+m\geq 0$ pour le restant). Si $a\ne 0$, en d\'ecompo\-sant $\Q-(a+p^n\Z)=\Q-(a+p^{n-1}\Z)\amalg (a+p^{n-1}\Z)-(a+p^n\Z)$, puis $\Q-(a+p^{n-1}\Z)=\Q-(a+p^{n-2}\Z)\amalg (a+p^{n-2}\Z)-(a+p^{n-1}\Z)$ etc. jusqu'\`a arriver \`a $\Q-(a+p^{n-m}\Z)$ avec $n-m<{\rm val}(a)$ et $n\leq m$, on d\'eduit de (\ref{chaud3}) que (\ref{chaud2}) pour $a\ne 0$ et $j\geq {\rm val}(\alpha_p)$ d\'ecoule de (\ref{chaud1}) et de (\ref{chaud2}) pour $a=0$ et $n\leq 0$ (utiliser (\ref{chaud3}) pour les morceaux compacts dans la d\'ecomposition puis d\'evelopper $(z-a)^{k-2-j}$ et utiliser (\ref{chaud2}) avec $a=0$ et $n'=n-m\leq 0$ pour le restant). Autrement dit, $(\ref{chaud1})$ et $(\ref{chaud2})$ sont \'equivalents \`a :
\begin{itemize}
\item[(i)] (\ref{chaud1});
\item[(ii)] (\ref{chaud2}) pour $n\geq 0$;
\item[(iii)] (\ref{chaud2}) pour $a=0$ et $n\leq 0$.
\end{itemize}
Par ailleurs, tout $\mu\in \pi(\alpha)^*$ s'\'ecrit $\mu=(\mu_1,\mu_2)$ o\`u $\mu_i\in {\rm Pol}^{k-2}(\Z,L)^*$ (si $f\in \pi(\alpha)$,
$\int_{\Q}f(z)d\mu(z)=\int_{\Z}f_1(z)d\mu_1(z)+\int_{\Z}f_2(z)d\mu_2(z)$). Un calcul facile (laiss\'e au lecteur) montre que $\mu_1$ et $\mu_2$ v\'erifient (\ref{borneamice}) pour $r={\rm val}(\alpha_p)$ et $d=k-2$ avec $\|\mu_i\|_{{\rm val}(\alpha_p),k-2}\leq C$ (i.e. $\mu_1, \mu_2\in \mathcal{C}^{{\rm val}(\alpha_p)}(\Z,L)$ par le th\'eor\`eme \ref{amice2} avec leurs normes born\'ees, i.e. $\mu$ est dans une boule de $B(\alpha)^*\subset \pi(\alpha)^*$) si et seulement si $\mu$ v\'erifie (quitte \`a modifier $C$) :
\begin{equation}\label{ch0} 
\int_{a+p^{n}\Z}(z-a)^jd\mu(z)\in Cp^{n(j-{\rm val}(\alpha_p))}\O 
\end{equation} 
pour tout $a\in p\Z$, tout $j\in\{0,\cdots,k-2\}$ et tout entier $n\geq 1$, puis :
\begin{equation}\label{ch1} 
\int_{a^{-1}+p^{n-2{\rm val}(a)}\Z}(\beta\alpha^{-1})(z)| z|^{-1}z^{k-2-j}(1-az)^jd\mu(z)\in Cp^{n(j-{\rm val}(\alpha_p))}\O 
\end{equation} 
pour tout $a\in \Z-\{0\}$, tout $j\in\{0,\hdots,k-2\}$ et tout entier $n>{\rm val}(a)$, et enfin : 
\begin{equation}\label{ch2} 
\int_{\Q-p^{n}\Z}(\beta\alpha^{-1})(z)| z|^{-1}z^{k-2-j}d\mu(z)\in Cp^{n({\rm val}(\alpha_p)-j)}\O  
\end{equation} 
pour tout $j\in\{0,\hdots,k-2\}$ et tout entier $n\leq 0$. En d\'eveloppant $z^{k-2-j}=((z-a^{-1})+a^{-1})^{k-2-j}$ dans (\ref{ch1}), un calcul montre que, quitte \`a modifier $C$, (\ref{ch0}), (\ref{ch1}) et (\ref{ch2}) sont \'equivalents \`a :
\begin{itemize}
\item[(iv)] (\ref{chaud1}) pour $a\ne 0$;
\item[(v)] (\ref{chaud1}) pour $a=0$ et $n\geq 0$;
\item[(vi)] (\ref{chaud2}) pour $a=0$ et $n\leq 0$.
\end{itemize}
Si $\mu$ est comme en (\ref{thedual}), i.e. si $\mu$ v\'erifie (i) \`a
(iii), alors {\it a fortiori} $\mu$ v\'erifie (iv) \`a (vi) et donc
$\mu\in B(\alpha)^*\subset \pi(\alpha)^*$. Mais on a plus. En faisant
tendre $n$ vers $-\infty$ dans (\ref{chaud1}) lorsque $a=0$, on voit
que (\ref{chaud1}) pour $j<{\rm val}(\alpha_p)$ et $a=0$ implique que
$\mu$ annule les fonctions $z^j\in B(\alpha)$. En faisant tendre $n$
vers $+\infty$ dans (\ref{chaud2}), on voit que (\ref{chaud2}) pour
$j<{\rm val}(\alpha_p)$ implique que $\mu$ annule les fonctions
$(\beta\alpha^{-1})(z-a)| z-a|^{-1}(z-a)^{k-2-j}\in
B(\alpha)$. Un examen plus approfondi (sans difficult\'e mais que nous
omettons pour ne pas allonger la preuve) montre que les conditions (i)
\`a (iii) pr\'ec\'edentes sont en fait {\it \'equivalentes} aux
conditions (iv) \`a (vi) avec les deux conditions suppl\'ementaires
que $\mu$ annule les fonctions $z^j$ pour $j<{\rm val}(\alpha_p)$ et les
fonctions $(\beta\alpha^{-1})(z-a)| z-a|^{-1}(z-a)^{k-2-j}$ pour $a\in \Q$ et $j<{\rm val}(\alpha_p)$,
c'est-\`a-dire les fonctions de $L(\alpha)$. Autrement dit, on obtient
que le Banach dual du compl\'et\'e cherch\'e est isomorphe dans
$\pi(\alpha)^*$ au sous-espace de Banach de $B(\alpha)^*$ form\'e des
$\mu$ qui annulent $L(\alpha)$, c'est-\`a-dire \`a
$(B(\alpha)/L(\alpha))^*$. En particulier, $(B(\alpha)/L(\alpha))^*$
est un $\G$-Banach unitaire. Comme les Banach ne sont pas r\'eflexifs,
nous allons devoir faire un passage par les topologies faibles pour
d\'eduire l'isomorphisme de l'\'enonc\'e. L'injection
$B(\alpha)/L(\alpha)\hookrightarrow ((B(\alpha)/L(\alpha))^*)^*$
\'etant une immersion ferm\'ee $\G$-\'equivariante,
$B(\alpha)/L(\alpha)$ est aussi un $\G$-Banach unitaire. Cela
entra\^\i ne facilement que l'application $\pi(\alpha)\to
B(\alpha)/L(\alpha)$ induit un morphisme $\G$-\'equivariant continu du
compl\'et\'e unitaire ci-dessus de $\pi(\alpha)$ vers
$B(\alpha)/L(\alpha)$, donc un morphisme continu sur les duaux munis
de leur topologie faible (qui sont des {\og modules compacts \`a
isog\'enie pr\`es \fg} au sens de \cite{ST3}). Mais on vient de voir que ce morphisme sur les duaux \'etait bijectif (et m\^eme un isomorphisme topologique pour les topologies fortes). Par \cite[lemme 4.2.2]{Br2}, on en d\'eduit que c'est aussi un isomorphisme topologique pour les topologies faibles. Par dualit\'e (cf. \cite[th\'eor\`eme 1.2]{ST3}), on obtient l'isomorphisme topologique $\G$-\'equivariant de l'\'enonc\'e. Le cas $\beta$ se traite de m\^eme.
\end{proof}

Rappelons qu'il existe, \`a multiplication par un scalaire non nul pr\`es, un unique morphisme non nul $\G$-\'equivariant :
\begin{equation}\label{entrelisse}
I^{\rm lisse}:{\rm Ind}_{\B}^{\G}\beta\otimes \alpha\nrm^{-1}\longrightarrow {\rm Ind}_{\B}^{\G}\alpha\otimes \beta\nrm^{-1}
\end{equation}
qui est un isomorphisme non trivial lorsque $\alpha\ne\beta$ et $\alpha\ne \beta\nrm$, qui est l'identit\'e lorsque $\alpha=\beta$ et qui a un noyau et un conoyau de dimension $1$ lorsque $\alpha=\beta\nrm$ (voir \cite[\S4.5]{Bu} par exemple). En termes de fonctions localement constantes sur $\Q$, ce morphisme lorsque $\alpha\ne\beta$ est donn\'e explicitement par :
\begin{eqnarray}\label{integrale} 
I^{\rm lisse}(h)(z)&=&\int_{\Q}\!\!(\alpha\beta^{-1})(x)|x|^{-1}h(z+x^{-1})dx\\
\nonumber &=&\int_{\Q}\!\!(\beta\alpha^{-1})(x)|x|^{-1} h(z+x)dx\\
\nonumber &=&\int_{\Q}\!\!(\beta\alpha^{-1})(x-z)|x-z|^{-1} h(x)dx
\end{eqnarray} 
o\`u $dx$ est la mesure de Haar sur $\Q$ (\`a valeurs dans $\Q \subset L$). Comme la th\'eorie des repr\'esentations lisses est alg\'ebrique, il n'y a pas de probl\`e\-mes de convergence dans les int\'egrales ci-dessus car on peut toujours remplacer les sommes infinies aux voisinages de $0$ ou de $-\infty$ par des expressions alg\'ebriques en $p\alpha_p\beta_p^{-1}$ parfaitement d\'efinies. En tensorisant par l'application identit\'e sur ${\rm Sym}^{k-2}L^2$, on en d\'eduit un morphisme non nul $\G$-\'equivariant :
\begin{equation}\label{intertw}
I:\pi(\beta)\longrightarrow \pi(\alpha)
\end{equation}
qui est un isomorphisme lorsque $\alpha\ne \beta\nrm$. 

\begin{coro}\label{entrepadique}
Les repr\'esentations $B(\alpha)/L(\alpha)$ et
$B(\beta)/L(\beta)$ sont des $\G$-Banach unitaires et on a un diagramme commutatif $\G$-\'equiva\-riant :
$$\begin{matrix}B(\beta)/L(\beta)&\overset{\widehat I}{\longrightarrow} &B(\alpha)/L(\alpha)\\
\uparrow &&\uparrow \\ \pi(\beta)&\overset{I}{\longrightarrow} &\pi(\alpha)\end{matrix}$$  
o\`u $I$ est le morphisme $\G$-\'equivariant de (\ref{intertw}). Lorsque $\alpha\ne \beta\nrm$, les fl\`eches $I$ et $\widehat I$ sont des isomorphismes.
\end{coro} 
\begin{proof} 
Cela d\'ecoule du th\'eor\`eme \ref{complete} car l'image par une fl\`eche $\G$-\'equivariante d'un $\O[\G]$-module de type fini est aussi un $\O[\G]$-module de type fini.
\end{proof} 

Lorsque $\alpha=\beta\nrm$ (ce qui implique ${\rm val}(\beta_p)=(k-2)/2$), il est \'evident que $B(\beta)/L(\beta)$ est non nul puisque $L(\beta)$ est dans ce cas une repr\'esentation de dimension finie, isomorphe via (\ref{algebraic}) \`a $(\beta\circ\det)\otimes_L {\rm Sym}^{k-2}L^2$. On verra quels sont alors l'image et le noyau de $\widehat I$ au \S\ref{steinberg}. Dans les autres cas, le th\'eor\`eme \ref{complete} ne d\'emontre en rien que les espaces de Banach $B(\alpha)/L(\alpha)$ et $B(\beta)/L(\beta)$ sont non nuls. Mais on a :

\begin{prop}\label{existencereseau}
Si $\alpha\ne \beta\nrm$, le Banach $B(\alpha)/L(\alpha)$ (resp. $B(\beta)/L(\beta)$) est non nul si et seulement si $\pi(\alpha)$ (resp. $\pi(\beta)$) poss\`ede un $\O$-r\'eseau stable par $\G$. Si $\alpha=\beta\nrm$, le Banach $B(\alpha)/L(\alpha)$ est non nul si et seulement si $\pi(\alpha)$ poss\`ede un $\O$-r\'eseau stable par $\G$. 
\end{prop}
\begin{proof}
Rappelons qu'un $\O$-r\'eseau est par d\'efinition un sous-$\O$-module g\'en\'erateur qui ne contient pas de $L$-droite. Supposons d'abord $\alpha\ne \beta\nrm$, de sorte que les repr\'esentations $\pi(\alpha)$ et $\pi(\beta)$ sont (alg\'ebriquement) irr\'eductibles. Si $B(\alpha)/L(\alpha)\ne 0$, l'application canonique $\pi(\alpha)\to B(\alpha)/L(\alpha)$ est injective car non nulle (car d'image dense) et une boule unit\'e de $B(\alpha)/L(\alpha)$ stable par $\G$ induit un r\'eseau stable par $\G$ sur $\pi(\alpha)$. Inversement, supposons que $\pi(\alpha)$ poss\`ede un $\O$-r\'eseau stable par $\G$, alors pour tout $f$ non nul dans $\pi(\alpha)$, $\O[\G]f\subset \pi(\alpha)$ est un $\O$-r\'eseau de $\pi(\alpha)$ de type fini comme $\O[\G]$-module. Il est g\'en\'erateur car $\pi(\alpha)$ est irr\'eductible et il ne contient pas de $\O$-droite car, \`a multiplication pr\`es par un scalaire, il est contenu dans un $\O$-r\'eseau stable par $\G$ de $\pi(\alpha)$. L'application de $\pi(\alpha)$ dans son compl\'et\'e par rapport \`a $\O[\G]f$, qui est $B(\alpha)/L(\alpha)$ par le th\'eor\`eme \ref{complete}, est  alors injective et en particulier $B(\alpha)/L(\alpha)\ne 0$. Lorsque $\alpha=\beta\nrm$, ce qui suppose $k>2$, $\pi(\alpha)$ n'est plus irr\'eductible et a un quotient isomorphe \`a $(\beta\circ\det)\otimes_L{\rm Sym}^{k-2}L^2$. Si $B(\alpha)/L(\alpha)\ne 0$, l'application non nulle $\pi(\alpha)\to B(\alpha)/L(\alpha)$ reste injective sinon elle induirait une injection non nulle $(\beta\circ\det)\otimes_L{\rm Sym}^{k-2}L^2\hookrightarrow B(\alpha)/L(\alpha)$ ce qui est impossible car, pour $k>2$, $(\beta\circ\det)\otimes_L{\rm Sym}^{k-2}L^2$ ne poss\`ede pas de $\O$-r\'eseau stable par $\G$.
\end{proof}

Notons que, lorsque $\alpha=\beta\nrm$, $\pi(\beta)$ ne peut poss\'eder de $\O$-r\'eseau stable par $\G$ puisque sa sous-repr\'esentation irr\'eductible $(\beta\circ\det)\otimes_L{\rm Sym}^{k-2}L^2$ n'en poss\`ede pas. Dans ce cas, l'application $\pi(\beta)\to B(\beta)/L(\beta)$ est non injective (son noyau est pr\'ecis\'ement $(\beta\circ\det)\otimes_L{\rm Sym}^{k-2}L^2$). Lorsque $\alpha$ et $\beta$ sont non ramifi\'es, via la proposition \ref{existencereseau}, on peut montrer pour des petites valeurs de $k$ ou pour les valeurs de $(k,\alpha,\beta)$ provenant des formes modulaires que les Banach $B(\alpha)/L(\alpha)$ et $B(\beta)/L(\beta)$ sont non nuls, voir par exemple \cite{Br1}, \cite[\S1.3]{Br2}, \cite{Br3}, \cite{Em}. On va voir dans la suite que la non nullit\'e pour tout $k$ et tout $\alpha$, $\beta$, au moins si $\alpha\ne\beta$, d\'ecoule de la th\'eorie des $(\varphi,\Gamma)$-modules. Une autre approche possible, purement en termes de th\'eorie des repr\'esentations, est pr\'esent\'ee dans \cite[\S2]{Em} et \cite[\S\S 5-6]{Em2}.

\section{Repr\'esentations de $\G$ et $(\varphi,\Gamma)$-modules}

Le but de cette partie est de d\'emontrer l'existence d'un isomorphisme topologique canonique entre $(\varprojlim_{\psi} \dfont(V))^{\rm b}$ et le dual $\breuil(V)^*$ (muni de sa topologie faible) lorsque $\alpha\ne\beta$ et d'en d\'eduire que $\breuil(V)$ alors est toujours non nul, topologiquement irr\'eductible et admissible. Ces \'enonc\'es \'etaient conjectur\'es (et des cas particuliers d\'emontr\'es) dans \cite{Br1} et \cite{Br2}. Le fait remarquable est que ces \'enonc\'es, enti\`erement du c\^ot\'e ${\rm GL}_2$, se d\'emontrent en passant par le c\^ot\'e {\it galoisien}. On fixe une fois pour toutes une repr\'esentation apc irr\'eductible $V$ comme au \S\ref{cristallin} avec $\dcris(V)=D(\alpha,\beta)$ et on suppose jusqu'\`a la fin que $\alpha\ne\beta$.

\subsection{Deux lemmes}\label{deuxlemmes}

Le but de ce paragraphe est de d\'emontrer deux lemmes techniques mais importants utilis\'es dans les paragraphes suivants. On utilise sans commentaire certaines notations du \S\ref{cristallin}.

\begin{lemm}\label{explicit} 
Soit $m\in \N$, $m\geq m(V)$, $w_{\alpha}$, $w_{\beta}\in \calR^+$ et $\mu_{\alpha}$, $\mu_{\beta}$ les distributions localement analytiques sur $\Z$ correspondantes par (\ref{Yvette}). 
La condition :  
$$\varphi^{-m}(w_{\alpha}\otimes e_{\alpha}+w_{\beta}\otimes
e_{\beta})\in {\rm Fil}^0(L_m[[t]]\otimes_LD(\alpha,\beta))$$  
est \'equivalente aux \'egalit\'es dans $\Qpbar$ :  
$$\left(\sum_{x\in\Z^{\times}/(1+p^{m(V)}\Z)}(\beta\alpha^{-1})(x)\eta_{p^m}^{p^{m-m(V)}x}\right)\alpha_p^m\int_{\Z}z^j\eta_{p^m}^zd\mu_{\alpha}(z)
=\beta_p^m\int_{\Z}z^j\eta_{p^m}^zd\mu_{\beta}(z)$$  
pour tout $j\in \{0,\hdots,k-2\}$ et toute racine primitive $p^m$-i\`eme $\eta_{p^m}$ de $1$ dans $\Qpbar$.
\end{lemm}
 
\begin{proof}
On a :
\begin{equation*}
\varphi^{-m}(X)=\zeta_{p^m}{\rm exp}(t/p^m)-1=\zeta_{p^m}({\rm exp}(t/p^m)-1)+\zeta_{p^m}-1
\end{equation*} 
dans $F_m[[t]]$ (voir \S\ref{cristallin}). En posant $w_{\alpha}=\sum_{i=0}^{+\infty}a_iX^i$ et
$w_{\beta}=\sum_{i=0}^{+\infty}b_iX^i$ ($a_i$, $b_i\in L$), la condition sur ${\rm Fil}^0$ est \'equivalente \`a :  
\begin{multline} \label{filo}
\alpha_p^m\sum_{i=0}^{+\infty}a_i\varphi^{-m}(X)^i 
e_{\alpha}+\beta_p^m\sum_{i=0}^{+\infty}b_i\varphi^{-m}(X)^ie_{\beta}\in\\
L_m[[t]](e_{\alpha}+G(\beta\alpha^{-1})e_{\beta})\oplus   
({\rm exp}(t/p^m)-1)^{k-1}(L_m[[t]]e_{\beta}) 
\end{multline} 
en notant que ${\rm exp}(t/p^m)-1$ engendre ${\rm gr}^1(\Q[[t]])=\Q\overline
t$. On peut supposer $L$ aussi grand que l'on veut en (\ref{filo}), et en particulier contenant $F_m$. En utilisant :
$$L_m=F_m \otimes_{\Q} L=\prod_{F_m \hookrightarrow L} L,$$ 
en d\'eveloppant $\varphi^{-m}(X)^i=(\zeta_{p^m}({\rm exp}(t/p^m)-1)+\zeta_{p^m}-1)^i$ et en rempla\c cant $e_{\alpha}$ par $(e_{\alpha}+G(\beta\alpha^{-1})e_{\beta})-G(\beta\alpha^{-1})e_{\beta}$ dans le membre de gauche de (\ref{filo}), 
un calcul facile montre que la condition (\ref{filo}) est \'equivalente aux \'egalit\'es dans $L$ :
\begin{multline}\label{combin} 
\left(\sum_{x\in\Z^{\times}/(1+p^{m(V)}\Z)}(\beta\alpha^{-1})(x)\eta_{p^m}^{p^{m-m(V)}x}\right)\alpha_p^m\sum_{i=j}^{+\infty}a_i\binom{i}{j}(\eta_{p^m}-1)^{i-j} \\
= \beta_p^m\sum_{i=j}^{+\infty}b_i\binom{i}{j}(\eta_{p^m}-1)^{i-j} 
\end{multline} 
pour tout $j\in \{0,\hdots,k-2\}$ et {\it toute} racine primitive $p^m$-i\`eme $\eta_{p^m}$ de $1$. Noter que les s\'eries en (\ref{combin}) convergent bien car ${\rm val}(\eta_{p^m}-1)>0$. En se souvenant que $\alpha_i=\int_{\Z}\binom{z}{i}d\mu_{\alpha}(z)$ et en utilisant le d\'eveloppement de Mahler (\ref{mahler}) :  
$$\binom{z}{j}\eta_{p^m}^{z-j}=\sum_{i=j}^{+\infty}\binom{i}{j}(\eta_{p^m}-1)^{i-j}\binom{z}{i}$$  
on obtient : 
\begin{eqnarray*} 
\sum_{i=j}^{+\infty} a_i\binom{i}{j}(\eta_{p^m}-1)^{i-j}&=&\int_{\Z}
\left(\sum_{i=j}^{+\infty}\binom{z}{i}\binom{i}{j}(\eta_{p^m}-1)^{i-j}\right)d\mu_{\alpha}(z)\\  
&=&\eta_{p^m}^{-j}\int_{\Z}\binom{z}{j}\eta_{p^m}^{z}d\mu_{\alpha}(z) 
\end{eqnarray*} 
(la s\'erie $\sum_{i=j}^{n}\binom{i}{j}(\eta_{p^m}-1)^{i-j}\binom{z}{i}$ convergeant vers $\sum_{i=j}^{+\infty}\binom{i}{j}(\eta_{p^m}-1)^{i-j}\binom{z}{i}$ dans ${\rm An}(\Z,L)$ (cf. \cite[\S2.1.2]{Co2}), on peut inverser $\int$ et $\sum$). On a la m\^eme \'egalit\'e avec $b_i$ et $\mu_{\beta}$. Avec (\ref{combin}), on en d\'eduit le r\'esultat. 
\end{proof} 

Via l'identification $\Q/\Z=\ZZ[1/p]/\ZZ$, on peut d\'efinir le nombre complexe alg\'ebrique $e^{2i\pi z}$ pour tout $z\in \Q$ (par exemple, $e^{2i\pi z}=1$ si $z\in \Z$). En fixant des plongements $\overline{\mathbf Q}\hookrightarrow {\mathbf C}$ et $\overline{\mathbf Q}\hookrightarrow \Qpbar$, on peut voir $e^{2i\pi z}$ comme un \'el\'ement de $\Qpbar$. On obtient ainsi un caract\`ere additif localement constant $\Q\to \Qpbar^{\times}$, $z\mapsto e^{2i\pi z}$ trivial sur $\Z$. Ce que l'on fera dans la suite ne d\'ependra pas du choix de ce caract\`ere, i.e. du choix des plongements.

Notons ${\rm Pol}^d(\Q,L)$ le $L$-espace vectoriel des fonctions localement polynomiales \`a support compact $f:\Q\to L$ de degr\'e (local) au plus $d$. Si $\mu$ est une forme lin\'eaire sur ${\rm Pol}^d(\Q,L)$, $U$ un ouvert compact de $\Q$ et $f:\Q\to L$ une fonction localement polynomiale de degr\'e au plus $d$ \`a support quelconque, on note comme d'habitude $\int_Uf(z)d\mu(z)\=\mu({\bf 1}_U(z)f(z))$.

\begin{lemm}\label{crucial}
Soient $\mu_{\alpha}$ et $\mu_{\beta}$ deux formes lin\'eaires sur l'espa\-ce ${\rm Pol}^{k-2}(\Q,L)$. Les
\'enonc\'es suivants sont \'equivalents :
\begin{itemize}
\item[(i)] Pour tout $j\in \{0,\hdots,k-2\}$, tout $y\in \Q^{\times}$ et tout $N\geq {\rm val}(y)+m(V)$, on a : 
\begin{multline}\label{rajout} 
\int_{p^{-N}\Z}z^je^{2i\pi zy}d\mu_{\beta}(z)
 \\ =
\left(\sum_{x\in\Z^{\times}/(1+p^{m(V)}\Z)}(\beta\alpha^{-1})(x)
e^{\frac{2i\pi xy} {p^{\val(y)+m(V)}} } \right)\left(\frac{\beta_p}{\alpha_p}\right)^{{\rm val}(y)}\int_{p^{-N}\Z}z^je^{2i\pi zy}d\mu_{\alpha}(z). 
\end{multline} 
\item[(ii)] Pour tout $f\in \pi(\alpha)$ \`a support compact (comme fonction sur $\Q$ via (\ref{commefonction})) tel que $I(f)\in \pi(\beta)$ est aussi \`a support compact (comme fonction sur $\Q$ via (\ref{commefonction})), on a :  
\begin{equation}\label{debutinter} 
\int_{\Q}I(f)(z)d\mu_{\alpha}(z)=C(\alpha_p,\beta_p)\int_{\Q}f(z)d\mu_{\beta}(z)
\end{equation}
o\`u : 
$$C(\alpha_p,\beta_p)\=\frac{1-\frac{\beta_p}{p\alpha_p}}{1-\frac{\alpha_p}{\beta_p}}\ \text{si $\beta\alpha^{-1}$ est non ramifi\'e,}$$
et  $$C(\alpha_p,\beta_p)\=\left(\frac{\beta_p} {p\alpha_p}\right)^{m(V)}\ \text{si $\beta\alpha^{-1}$ est ramifi\'e.}$$
\end{itemize}
\end{lemm} 
\begin{proof}
Soit $h:\Q\to L$ une fonction localement constante \`a support compact et $\widehat h$ la transform\'ee de Fourier (usuelle) de $h$. Rappelons que $\widehat h$ est aussi une fonction localement constante sur $\Q$ \`a support compact telle que $\widehat h(x)=\int_{\Q}h(z)e^{-2i\pi zx}dz$ et $h(z)=\int_{\Q}\widehat h(x)e^{2i\pi zx}dx$ o\`u $dx$, $dz$ d\'esignent la mesure de Haar sur $\Q$. Pour $|z|\gg 0$, on a par (\ref{integrale}) :
$$I^{\rm lisse}(h)(z)=(\beta \alpha^{-1})(z)| z|^{-1}\int_{\Q}h(x)dx=(\beta \alpha^{-1})(z)| z|^{-1}\widehat h(0)$$
et on voit que $I^{\rm lisse}(h)$ est \`a support compact dans $\Q$ si et seulement si $\widehat h(0)=0$. Supposons donc $\widehat h(0)=0$ et soit $N\in \NN$ tel que $h$ et $I^{\rm lisse}(h)$ ont leur support dans $p^{-N}\Z$ et tel que $\widehat h_{\mid p^N\Z}=0$. Pour $j\in \{0,\hdots,k-2\}$ et $z\in p^{-N}\Z$, on a :  
\begin{eqnarray*} 
I(z^jh)(z)&=&z^jI^{\rm lisse}(h)(z)\ =\ z^j\int_{p^{-N}\Z} (\beta\alpha^{-1})(x)|x|^{-1}h(z+x)dx\\ 
&=&z^j\int_{p^{-N}\Z} (\beta\alpha^{-1})(x)|x|^{-1}\left(\int_{\Q-p^N\Z}\widehat h(y)e^{2i\pi y(z+x)}dy\right)dx\\ 
&=&z^j\int_{\Q-p^N\Z}\widehat h(y)e^{2i\pi zy}\left(\int_{p^{-N}\Z}
 (\beta\alpha^{-1})(x)|x|^{-1}e^{2i\pi xy}dx\right)dy. 
\end{eqnarray*} 
On a : 
$$\int_{p^{-N}\Z}(\beta\alpha^{-1})(x)|x|^{-1}e^{2i\pi xy}dx=\sum_{\ell=-N}^{+\infty}p^\ell \int_{p^\ell \Z^{\times}}(\beta\alpha^{-1})(x)e^{2i\pi xy}dx$$ 
o\`u, avec les notations de la preuve du lemme \ref{appartient} :
\begin{eqnarray*}
\int_{p^\ell \Z^{\times}}(\beta\alpha^{-1})(x)e^{2i\pi xy}dx&=&\sum_{a_i\in S}\int_{p^\ell a_i+p^{\ell+m(V)}\Z}(\beta\alpha^{-1})(x)e^{2i\pi xy}dx\\
&=&\left(\frac{\alpha_p}{\beta_p}\right)^\ell \sum_{a_i\in S}(\beta\alpha^{-1})(a_i)\int_{p^\ell a_i+p^{\ell+m(V)}\Z}e^{2i\pi xy}dx.
\end{eqnarray*}
Si $\ell+m(V)<-{\rm val}(y)$, on v\'erifie facilement que $\int_{p^\ell a_i+p^{\ell+m(V)}\Z}e^{2i\pi xy}dx=0$ et si $\ell+m(V)\geq -{\rm val}(y)$, on a $\int_{p^\ell a_i+p^{\ell+m(V)}\Z}e^{2i\pi xy}dx=p^{-\ell-m(V)}e^{2i\pi p^\ell a_i y}$. Supposons d'abord $\beta\alpha^{-1}$ ramifi\'e. Alors on a $\sum_{a_i\in S}(\beta\alpha^{-1})(a_i)e^{2i\pi p^\ell a_iy}=0$ si $\ell+m(V)> -{\rm val}(y)$ et, si $\ell+m(V)= -{\rm val}(y)$ :
$$\sum_{a_i\in S}(\beta\alpha^{-1})(a_i)e^{2i\pi p^\ell a_iy}=\sum_{x\in\Z^{\times}/(1+p^{m(V)}\Z)}(\beta\alpha^{-1})(x)e^{\frac{2i\pi xy}{p^{{\rm val}(y)+m(V)}}}.$$
Comme $N\geq {\rm val}(y)+m(V)$, on a donc :
\begin{multline*}
\int_{p^{-N}\Z}(\beta\alpha^{-1})(x)|x|^{-1}e^{2i\pi xy}dx \\ 
=\left(\frac{\beta_p}{p\alpha_p}\right)^{m(V)}\left(\frac{\beta_p}{\alpha_p}\right)^{{\rm val}(y)} \sum_{x\in\Z^{\times}/(1+p^{m(V)}\Z)} (\beta\alpha^{-1})(x) e^{\frac{2i\pi xy}{p^{{\rm val}(y)+m(V)}}}.
\end{multline*}
Supposons maintenant $\beta\alpha^{-1}$ non ramifi\'e. Alors on a $\sum_{a_i\in S}(\beta\alpha^{-1})(a_i)e^{2i\pi p^\ell a_iy}=p-1$ si $\ell +1> -{\rm val}(y)$ et $\sum_{a_i\in S}(\beta\alpha^{-1})(a_i)e^{2i\pi p^ \ell a_iy}=-1$ si $\ell+1= -{\rm val}(y)$ (rappelons que $m(V)=1$). Comme $N\geq {\rm val}(y)+m(V)$, on a donc :
\begin{eqnarray*}
\int_{p^{-N}\Z}(\beta\alpha^{-1})(x)|x|^{-1}e^{2i\pi xy}dx&=&-\frac{1}{p} \left(\frac{\alpha_p}{\beta_p}\right)^{-{\rm val}(y)-1}+\frac{p-1}{p}\sum_{\ell=-{\rm val}(y)}^{+\infty}\left(\frac{\alpha_p}{\beta_p} \right)^\ell \\
&=&\frac {1-\frac{\beta_p}{p \alpha_p}} {1-\frac{\alpha_p}{\beta_p}}
\left(\frac{\beta_p}{\alpha_p}\right)^{{\rm val}(y)}.
\end{eqnarray*}
Pour $y\in \Q^{\times}$, posons $H(\beta\alpha^{-1})(y)\=1$ si $\beta\alpha^{-1}$ est non ramifi\'e et : 
$$H(\beta\alpha^{-1})(y)\=\sum_{x\in\Z^{\times}/(1+p^{m(V)}\Z)} (\beta\alpha^{-1})(x)e^{\frac{2i\pi xy}{p^{{\rm val}(y)+m(V)}}}$$ 
si $\beta\alpha^{-1}$ est ramifi\'e.
On en d\'eduit : 
\begin{equation}\label{fourier1} 
I(z^jh)(z)=C(\alpha_p,\beta_p)\int_{\Q-p^N\Z}\widehat h(y)z^je^{2i\pi zy}\left(\frac{\beta_p}{\alpha_p}\right)^{{\rm val}(y)}H(\beta\alpha^{-1})(y)dy 
\end{equation} 
pour $z\in p^{-N}\Z$ et $I(z^jh)(z)=0$ sinon. De m\^eme, on a :  
\begin{equation}\label{fourier2} 
z^jh(z)=\int_{\Q-p^N\Z}\widehat h(y)z^je^{2i\pi zy}dy 
\end{equation} 
pour $z\in p^{-N}\Z$ et $z^jh(z)=0$ sinon. Notons que (\ref{fourier1}) and (\ref{fourier2}) sont en fait des sommes finies sur le m\^eme ensemble (fini) de valeurs de $y$. En rempla\c cant $I(z^jh)(z)$ et $z^jh(z)$ dans (\ref{variant}) ci-dessous par les sommes finies (\ref{fourier1}) et (\ref{fourier2}), on voit que (i) entra\^\i ne (ii). R\'eciproquement, supposons que pour tout $j\in \{0,\hdots,k-2\}$, tout $N\in \N$ et toute fonction $h$ comme ci-dessus localement constante \`a support dans $p^{-N}\Z$, on a :
\begin{equation}\label{variant} 
\int_{p^{-N}\Z}I(z^jh)(z)d\mu_{\alpha}(z)=C(\alpha_p,\beta_p)\int_{p^{-N}\Z}z^jh(z)d\mu_{\beta}(z).
\end{equation} 
Soit $y\in \Q^{\times}$ tel que $N\geq {\rm val}(y)+m(V)$, $\widehat h(z)\= {\bf 1}_{y+p^N\Z}(z)$ et :
$$h(z)\=\int_{\Q}\widehat h(x)e^{2i\pi zx}dx=\int_{p^{N}\Z}e^{2i\pi z(y+x)}dx=\frac{1}{p^N}{\bf 1}_{p^{-N}\Z}(z)e^{2i\pi zy}.$$
Alors $I(z^jh)$ est aussi \`a support compact car $\widehat h_{\mid p^N \Z}=0$ et un calcul via (\ref{fourier1}) montre que :
$$I(z^jh)(z)=C(\alpha_p,\beta_p)\frac{1}{p^N}\left(\frac{\beta_p}{\alpha_p}\right)^{{\rm val}(y)}H(\beta\alpha^{-1})(y){\bf 1}_{p^{-N}\Z}(z)z^je^{2i\pi zy}.$$
On peut donc appliquer l'\'egalit\'e (\ref{variant}) \`a $h$ qui est alors exactement l'\'egalit\'e (\ref{rajout}) multipli\'ee par $p^{-N}$. Cela montre que (ii) entra\^\i ne (i) et ach\`eve la preuve.
\end{proof} 

\subsection{D'un monde \`a l'autre}\label{laouiw}

Le but de ce paragraphe est de construire un isomorphisme topologique $(\varprojlim_{\psi}\dfont(V))^{\rm b}\simeq \breuil(V)^*$ (lorsque $\alpha\ne\beta$).

Soit $T\subset V$ un $\O$-r\'eseau stable par $\g$. On reprend les notations du \S\ref{topotr}, en particulier on dispose du $\O[[X]]$-module de type fini $\dsharp(T)$ muni de la surjection $\psi:\dsharp(T)\twoheadrightarrow \dsharp(T)$ et de l'action semi-lin\'eaire de $\Gamma$ qui commute \`a $\psi$. On dispose aussi de l'isomorphisme topologique de la proposition \ref{ddiese} qui permet de remplacer $(\varprojlim_{\psi}\dfont(V))^{\rm b}$ par $(\varprojlim_{\psi}\dsharp(T))\otimes_{\O}L$ et on sait par le th\'eor\`eme \ref{limproj} que $(\varprojlim_{\psi}\dsharp(T)) \otimes_{\O}L$ co\"\i ncide avec les suites d'\'el\'ements $w_{\alpha,n}\otimes e_{\alpha} + w_{\beta,n}\otimes e_{\beta}$ de $\calR^+\otimes_L\dcris(V)$ telles que :
\begin{itemize}
\item[(i)] $\forall\ n\geq 0$, $w_{\alpha,n}$ (resp. $w_{\beta,n}$) est d'ordre ${\rm val}(\alpha_p)$ (resp. ${\rm val}(\beta_p)$) dans $\calR^+$ et $\|w_{\alpha,n}\|_{{\rm val}(\alpha_p)}$ (resp. $\|w_{\beta,n}\|_{{\rm val}(\beta_p)}$) est born\'e;
\item[(ii)] $\forall\ n\geq 0$ et $\forall\ m\geq 1$, on a :
$$\varphi^{-m}(w_{\alpha,n}\otimes e_{\alpha}+w_{\beta,n}\otimes 
e_{\beta})\in {\rm Fil}^0(L_m [[t]]\otimes_L\dcris(V));$$ 
\item[(iii)] $\forall\ n\geq 1$, $\psi(w_{\alpha,n})=\alpha_p^{-1}w_{\alpha,n-1}$ et  
$\psi(w_{\beta,n})=\beta_p^{-1}w_{\beta,n-1}$.
\end{itemize}

Nous allons d'abord d\'efinir une application $L$-lin\'eaire $(\varprojlim_{\psi}\dsharp(T))\otimes_{\O}L\to \pi(\alpha)^*$. Soit $\mu_{\alpha,n}$ et $\mu_{\beta,n}$ les distributions sur $\Z$ correspondant \`a $\alpha_p^nw_{\alpha,n}$ et $\beta_p^nw_{\beta,n}$ par la transform\'ee d'Amice-Mahler (\ref{Yvette}). On associe \`a $(\mu_{\alpha,n})_n$ et $(\mu_{\beta,n})_n$ deux distributions localement analytiques $\mu_{\alpha}$ et $\mu_{\beta}$ sur $\Q$ \`a support compact (i.e. deux formes lin\'eaires continues sur l'espace vectoriel des fonctions localement analytiques sur $\Q$ \`a support compact) en posant :
\begin{equation}\label{prolong} 
\int_{U}f(z)d\mu_{\alpha}(z)\=\int_{\Z}{\bf 1}_U(z/p^N)f(z/p^N)d\mu_{\alpha,N}(z)
\end{equation} 
(resp. avec $\beta$ au lieu de $\alpha$) o\`u $f:\Q\to L$ est localement analytique (\`a support quelconque) et $U$ est un ouvert compact de $\Q$ contenu dans $p^{-N}\Z$.

\begin{lemm}
La valeur $\int_{\Z}{\bf 1}_U(z/p^N)f(z/p^N)d\mu_{\alpha,N}(z)$ ne d\'epend pas du choix de $N$ tel que $U$ est contenu dans $p^{-N}\Z$.
\end{lemm}
\begin{proof}
Si $\mu$ est une distribution localement analytique sur $\Z$ correspondant \`a $w\in \calR^+$ par (\ref{Yvette}), il est facile de voir que la distribution localement analytique $\psi(\mu)$ correspondant \`a $\psi(w)$ v\'erifie :
\begin{equation}\label{psidis}
\int_{\Z}f(z)d\psi(\mu)(z)=\int_{p\Z}f(z/p)d\mu(z).
\end{equation}
On a donc :
\begin{eqnarray*}
\int_{\Z}{\bf 1}_U(z/p^N)f(z/p^{N})d\mu_{\alpha,N}(z)&=&\int_{\Z}{\bf 1}_U(z/p^N)f(z/p^{N})d\psi(\mu_{\alpha,N+1})(z)\\
&\overset{(\ref{psidis})}{=}&\int_{p\Z}{\bf 1}_U(z/p^{N+1})f(z/p^{N+1})d\mu_{\alpha,N+1}(z)\\
&=&\int_{\Z}{\bf 1}_U(z/p^{N+1})f(z/p^{N+1})d\mu_{\alpha,N+1}(z),
\end{eqnarray*}
en remarquant que ${\bf 1}_U(z/p^{N+1})f(z/p^{N+1})$ est \`a support dans $p\Z$.
\end{proof}

Par le lemme \ref{explicit}, la condition (ii) pr\'ec\'edente sur $(w_{\alpha,n},w_{\beta,n})_n$ est
\'equivalente aux \'egalit\'es :  
\begin{equation}\label{filtr} 
\left(\sum_{x\in\Z^{\times}/(1+p^{m(V)}\Z)}(\beta\alpha^{-1})(x)\eta_{p^m}^{p^{m-m(V)}x}\right)\alpha_p^{m-n}\int_{\Z}z^j\eta_{p^m}^zd\mu_{\alpha,n}(z)=\beta_p^{m-n}\int_{\Z}z^j\eta_{p^m}^zd\mu_{\beta,n}(z)  
\end{equation} 
pour tout $j\in \{0,\hdots,k-2\}$, tout $n\geq 0$, tout $m\geq m(V)$ et toute racine primitive $p^{m}$-i\`eme $\eta_{p^m}$ de $1$ dans $\Qpbar$. 
  
\begin{coro}\label{equival} 
Avec les notations pr\'ec\'edentes, la condition (ii) ci-dessus sur $(w_{\alpha,n},w_{\beta,n})_n$ est \'equivalente aux \'egalit\'es dans $\Qpbar$ : 
\begin{multline*}
\int_{p^{-N}\Z}z^je^{2i\pi zy}d\mu_{\beta}(z)=\\
\left(\sum_{x\in\Z^{\times}/(1+p^{m(V)}\Z)}(\beta\alpha^{-1})(x)
e^{\frac{2i\pi xy}{p^{{\rm val}(y)+m(V)}}}\right)\left(\frac{\beta_p}{\alpha_p}\right)^{{\rm val}(y)}\int_{p^{-N}\Z}z^je^{2i\pi zy}d\mu_{\alpha}(z)
\end{multline*}
pour tout $j\in \{0,\hdots,k-2\}$, tout $y\in \Q^{\times}$ et tout $N\geq {\rm val}(y)+m(V)$. 
\end{coro} 

\begin{proof} 
Cela r\'esulte de (\ref{prolong}) et de (\ref{filtr}) en remarquant que ${\bf 1}_U(z/p^N)=1$ si $U=p^{-N}\Z$ et $z\in \Z$, et en posant $n=N$ et $m=N-{\rm val}(y)\geq m(V)$. 
\end{proof} 
 
Par le lemme \ref{crucial}, on a donc : 
$$\int_{\Q}I(f)(z)d\mu_{\alpha}(z)
=C(\alpha_p,\beta_p)\int_{\Q}f(z)d\mu_{\beta}(z)$$
pour $f\in \pi(\alpha)$ \`a support compact tel que $I(f)\in \pi(\beta)$ est aussi \`a support compact. 

\begin{lemm}\label{versfin} 
Il y a une mani\`ere unique de prolonger $\mu_{\alpha}$ et $\mu_{\beta}$ comme
\'el\'ements respectivement de $\pi(\alpha)^*$ et $\pi(\beta)^*$ telle que, pour tout $f\in \pi(\beta)$ (vue comme fonction sur $\Q$ par (\ref{commefonction})) :  
\begin{equation}\label{extension} 
\int_{\Q}I(f)(z)d\mu_{\alpha}(z)=C(\alpha_p,\beta_p)\int_{\Q}f(z)d\mu_{\beta}(z)
\end{equation} 
\end{lemm} 
\begin{proof} 
Il suffit de montrer que, si $M$ est un entier suffisamment grand et si $j\in \{0,\hdots,k-2\}$, alors les int\'egrales : \[ \int_{\Q-p^{-M}\Z}(\beta\alpha^{-1})(z)|z|^{-1}z^jd\mu_{\alpha}(z)\ \text{et}\ \int_{\Q-p^{-M}\Z}(\alpha\beta^{-1})(z)|z|^{-1}z^jd\mu_{\beta}(z) \] 
sont uniquement d\'etermin\'ees. Si $h(z)\= {\mathbf 1}_{\Z}$, on a $\widehat h(z) =h(z)$ de sorte que $\widehat h(0)\ne 0$ et $I^{\rm lisse}(h)$ n'est pas \`a support compact dans $\Q$ (cf. preuve du lemme \ref{crucial}). Ainsi, quitte \`a multiplier $h$ par un scalaire non nul, on a un entier $M$ tel que $I^{\rm lisse}(h)(z)=(\beta\alpha^{-1})(z)|z|^{-1}$ dans $\pi(\alpha)$ pour $\mid\!z|\geq p^M$. L'\'egalit\'e (\ref{extension}) entra\^\i ne :
$$\int_{\Q-p^{-M}\Z}(\beta\alpha^{-1})(z)|z|^{-1}z^jd\mu_{\alpha}(z)=C(\alpha_p,\beta_p)\int_{\Z}z^jd\mu_{\beta}(z)-\int_{p^{-M}\Z}I^{\rm lisse}(h)(z)z^jd\mu_{\alpha}(z).$$
Cela permet d\'ej\`a de prolonger $\mu_{\alpha}$ \`a tout $\pi(\alpha)^*$. Le prolongement de $\mu_{\beta}$ \`a tout $\pi(\beta)^*$ s'en d\'eduit alors par (\ref{extension}) encore puisque $C(\alpha_p,\beta_p)\ne 0$.
\end{proof} 

Lorsque $\alpha=\beta \nrm$, on peut voir que la distribution $\mu_{\beta}\in \pi(\beta)^*$ du lemme \ref{versfin} est nulle contre $\beta\circ\det\otimes_L{\rm Sym}^{k-2}L^2\subset \pi(\beta)$. 

Avec les notations pr\'ec\'edentes, on d\'eduit du lemme \ref{versfin} une application $L$-lin\'eaire :
\begin{eqnarray}\label{versbanach}
(\varprojlim_{\psi}\dsharp(T))\otimes_{\O}L&\longrightarrow &\pi(\alpha)^*\\
\nonumber (w_{\alpha,n}\otimes e_{\alpha}+w_{\beta,n}\otimes
e_{\beta})_n&\longmapsto & \mu_{\alpha}\ \text{prolong\'e}
\end{eqnarray}
et notons que la d\'efinition de cette application utilise l'existence de l'entrelacement $I$ (lemme \ref{versfin}).

\begin{lemm}\label{actiongl2}
Soit $\gamma \in \Gamma$ tel que $\eps(\gamma)=a^{-1} \in \Zp^\times$, $z\in \Z$, $(v_n)_n\in \varprojlim_{\psi}\dsharp(T)$ et $\mu_{\alpha}\in \pi(\alpha)^*$ l'image de $(v_n)_n$ par (\ref{versbanach}). Alors :
\begin{itemize}
\item[(i)] $(\psi(v_n))_n$ s'envoie sur $\left(\begin{smallmatrix}1&0\\0&p\end{smallmatrix}\right)\cdot \mu_{\alpha}$;
\item[(ii)] $(\gamma (v_n))_n$ s'envoie sur $\left(\begin{smallmatrix}1&0\\0&a\end{smallmatrix}\right)\cdot \mu_{\alpha}$;
\item[(iii)] $(\varphi^n((1+X)^z)v_n)_n$ s'envoie sur $\left(\begin{smallmatrix}1&z\\0&1\end{smallmatrix}\right) \cdot \mu_{\alpha}$.
\end{itemize}
\end{lemm}
\begin{proof}
Cela d\'ecoule de (\ref{prolong}) et de propri\'et\'es simples de la transform\'ee d'Ami\-ce-Mahler (voir par exemple \cite[\S2.2.2]{Co2}). Nous laissons les d\'etails en exercice au lecteur.
\end{proof}

En particulier, le lemme \ref{actiongl2} induit une action du
groupe $\B$ sur $\varprojlim_{\psi}\dsharp(T)$, qui co\"{\i}ncide
bien s\^ur avec celle de la d\'efinition \ref{acgt} (en faisant agir les scalaires par multiplication par le caract\`ere central de $\pi(\alpha)^*$).

\begin{lemm}\label{unsens}
L'application (\ref{versbanach}) se factorise par une injection continue $\B$-\'equivariante :
$$(\varprojlim_{\psi}\dsharp(T))\otimes_{\O}L\hookrightarrow (B(\alpha)/L(\alpha))^*$$
(continue pour la topologie faible sur $(B(\alpha)/L(\alpha))^*$).
\end{lemm}

\begin{proof}
L'injectivit\'e d\'ecoule via (\ref{prolong}) de l'injectivit\'e de l'isomorphisme $\calR^+ \overset{\sim}{\to} {\rm An}(\Z,L)^*$ (cf. (\ref{Yvette})) et la $\B$-\'equivariance du lemme \ref{actiongl2}. Montrons que l'application $\varprojlim_{\psi}\dsharp(T)\to \pi(\alpha)^*$ est continue. Notons $\pi(\alpha)_{\rm c}\subset \pi(\alpha)$ (resp. $\pi(\beta)_{\rm c}\subset \pi(\beta)$) le sous-$L$-espace vectoriel des fonctions $f\in \pi(\alpha)$ (resp. $f\in \pi(\beta)$) \`a support compact dans $\Q$. La fl\`eche $\pi(\alpha)_{\rm c}\oplus \pi(\beta)_{\rm c}\overset {{\rm incl} \oplus I}{ \longrightarrow} \pi(\alpha)$ est surjective (voir e.g. la preuve du lemme \ref{versfin}) et induit une immersion ferm\'ee entre espaces de Fr\'echet :
$$\pi(\alpha)^*\hookrightarrow \pi(\alpha)_{\rm c}^*\oplus \pi(\beta)_{\rm c}^*.$$
Il suffit donc de montrer la continuit\'e des deux applications $\varprojlim_{\psi}\dsharp(T)\to \pi(\alpha)_{\rm c}^*$ et $\varprojlim_{\psi}\dsharp(T)\to \pi(\beta)_{\rm c}^*$, ce qui d\'ecoule apr\`es passage \`a la limite projective via (\ref{prolong}) de la continuit\'e de $\calR^+\overset{\sim}{\to} {\rm An}(\Z,L)^*$ et de celle de l'injection $\dsharp(T)\hookrightarrow \calR^+\otimes_L\dcris(V)$ (cf. la proposition  \ref{comtop}). Notons $\Pi(V)$ l'espace de Banach dual du module compact $\varprojlim_{\psi}\dsharp(T)$ par l'anti-\'equivalence de cat\'egorie de \cite[\S1]{ST3}. Il est muni d'une action continue unitaire de $\B$ par la proposition \ref{actgcont} (on peut utiliser les arguments de dualit\'e de la preuve de \cite[proposition 1.6]{ST3} pour la continuit\'e de l'action) et on a par ce qui pr\'ec\`ede un morphisme $\B$-\'equivariant (continu) $\pi(\alpha)\to \Pi(V)$. Par la propri\'et\'e universelle du compl\'et\'e de $\pi(\alpha)$ par rapport \`a un sous-$\O[\B]$-module g\'en\'erateur de type fini et par le th\'eor\`eme \ref{complete}, ce morphisme s'\'etend par continuit\'e en un morphisme $\B$-\'equivariant continu $B(\alpha)/L(\alpha)\to \Pi(V)$. En redualisant, ce dernier induit un morphisme continu $(\varprojlim_{\psi}\dsharp(T))\otimes_{\O}L\to (B(\alpha)/L(\alpha))^*$ qui est le morphisme de l'\'enonc\'e.
\end{proof}

Nous construisons maintenant une application continue $(B(\alpha)/L(\alpha))^*\to (\varprojlim_{\psi} \dfont(V))^{\rm b}$ inverse de la pr\'ec\'edente. 

Soit $\mu_{\alpha}\in (B(\alpha)/L(\alpha))^*$ et $\mu_{\beta}\=C(\alpha_p,\beta_p)^{-1}\widehat I\circ \mu_{\alpha}\in (B(\beta)/L(\beta))^*$ o\`u $\widehat I$ est le morphisme $\G$-\'equivariant du corollaire \ref{entrepadique}. On d\'efinit une suite $(\mu_{\alpha,n})_n$ de distributions localement analytiques sur $\Z$ en posant :
\begin{equation}\label{prolong2} 
\int_{\Z}f(z)d\mu_{\alpha,n}(z)\=\int_{p^{-n}\Z}f(p^nz)d\mu_{\alpha}(z)
\end{equation} 
et on d\'efinit de m\^eme $(\mu_{\beta,n})_n$. Soient $w_{\alpha,n},\ w_{\beta,n}\in \calR^+$ les \'el\'ements correspondant \`a $\alpha_p^{-n}\mu_{\alpha,n},\ \beta_p^{-n}\mu_{\beta,n}$ par (\ref{Yvette}). 

\begin{lemm}\label{routine}
La suite d'\'el\'ements $w_{\alpha,n}\otimes e_{\alpha} + w_{\beta,n}\otimes e_{\beta}$ de $\calR^+\otimes_L\dcris(V)$ satisfait les conditions (i), (ii) et (iii) du th\'eor\`eme \ref{limproj}.
\end{lemm}
\begin{proof}
La condition (iii) est \'evidente \`a partir de (\ref{psidis}) et la condition (ii) d\'ecoule des d\'efinitions, du lemme \ref{crucial} et du corollaire \ref{equival}. V\'erifions la condition (i). Revenant \`a la preuve du th\'eor\`eme \ref{complete}, on a en particulier que $\mu_{\alpha}$ satisfait (\ref{chaud1}) ce qui entra\^\i ne :
\begin{eqnarray*} 
\alpha_p^{-N}\int_{a+p^n\Z}(z-a)^jd\mu_{\alpha,N}(z)&=&\alpha_p^{-N}p^{Nj}
\int_{p^{-N}a+p^{n-N}\Z}(z-p^{-N}a)^jd\mu_{\alpha}(z)\\ 
&\in & C_{\mu_{\alpha}}p^{-N{\rm val}(\alpha_p)}p^{Nj}p^{(n-N)(j-{\rm val}(\alpha_p))}\O\\ 
&\in &C_{\mu_{\alpha}}p^{n(j-{\rm val}(\alpha_p))}\O
\end{eqnarray*}
pour tout $a\in \Z$, tout $j\in \{0,\hdots, k-2\}$ et tout $n\in \N$. Avec les notations du \S\ref{analyse}, cela entra\^\i ne pour tout $N\in \N$ :  
$$\|\alpha_p^{-N}\mu_{\alpha,N}\|_{{\rm val}(\alpha_p),k-2}\leq
c|C_{\mu_{\alpha}}|$$  
pour une constante $c\in \R_{\geq 0}$. On a une borne analogue pour les $\mu_{\beta,N}$. On en d\'eduit (i). 
\end{proof}

Par le lemme \ref{routine} et le th\'eor\`eme \ref{limproj}, on a une application $L$-lin\'eaire :
$$(B(\alpha)/L(\alpha))^*\longrightarrow (\varprojlim_{\psi}\dsharp(T))\otimes_{\O}L$$
et il est imm\'ediat \`a partir des d\'efinitions et du lemme \ref{versfin} de v\'erifier qu'elle est inverse de celle du lemme \ref{unsens}.

\begin{theo}\label{phigammabanach}
Il y a un unique isomorphisme topologique (\`a multiplication pr\`es par un scalaire non nul) entre les $L$-espaces vectoriels localement convexes (pour la topologie faible des deux c\^ot\'es): 
$$(\varprojlim_{\psi}\dfont(V))^{\rm b}\overset{\sim} {\longrightarrow} \breuil(V)^*$$ 
tel que l'action de $\left(\begin{smallmatrix}1&0\\0&p^\ZZ\end{smallmatrix}\right)$ sur $\breuil(V)^*$ correspond \`a $(v_n)_n\mapsto (\psi^\ZZ(v_n))_n$, l'action de $\left(\begin{smallmatrix}1&0\\0&\Z^{\times}\end{smallmatrix}\right)$ \`a celle de $\Gamma$ et l'action de $\left(\begin{smallmatrix}1&\Z\\0&1\end{smallmatrix}\right)$ \`a $(v_n)_n\mapsto ((1+X)^{p^n\Z}v_n))_n$. 
\end{theo}

\begin{proof}
L'existence d'un tel isomorphisme d\'ecoule des r\'esultats
pr\'ec\'e\-dents, sachant qu'une application bijective continue entre
deux {\og modules compacts \`a isog\'enie pr\`es \fg} est un isomorphisme
topologique (c'est la version duale par \cite{ST3} du th\'eor\`eme de
l'image ouverte entre espaces de Banach). Il reste \`a d\'emontrer
l'unicit\'e (\`a scalaire pr\`es) mais cela r\'esulte de la
proposition \ref{schur}.
\end{proof}

Rappelons que $H^i_{\rm Iw}(\Q,V)\= L \otimes_{\O} 
\varprojlim_n H^i({\rm Gal}(\Qpbar/F_n),T)$ 
o\`u $T$ est un $\O$-r\'eseau quelconque
de $V$ stable par $\g$ (voir le paragraphe \S\ref{pgmod}).

\begin{coro}\label{iwasawa}
On a un isomorphisme de $\O[[\Z^{\times}]]$-modules :
$$H^1_{\rm Iw}(\Q,V)\simeq {\breuil(V)^*}^{\left(\begin{smallmatrix}
1&0\\0&p^\ZZ\end{smallmatrix}\right)}$$
o\`u $\Z^{\times}$ agit via l'action de $\Gamma$ \`a gauche et via l'action de $\left(\begin{smallmatrix}1&0\\0&\Z^{\times}\end{smallmatrix}\right)$ \`a droite.
\end{coro}

\begin{proof}
Cela d\'ecoule du th\'eor\`eme \ref{phigammabanach} et de la proposition \ref{cciw}.
\end{proof}

\begin{rema}\label{iwasawa2}
Le $\O[[\Z^{\times}]]$-module $H^2_{\rm Iw}(\Q,V)$ s'identifie aussi
aux coinvariants de $\breuil(V)^*$ sous l'action de
$\left(\begin{smallmatrix}1&0\\0&p^\ZZ\end{smallmatrix}\right)$ 
et en fait, ces deux espaces sont nuls par le (ii) de la proposition
\ref{bpriw} parce que $V$ est irr\'eductible. 
En effet, le th\'eor\`eme \ref{phigammabanach} et le 
corollaire \ref{h2ddiese} nous disent que les
coinvariants d'un r\'eseau de $\breuil(V)^*$ sous l'action de
$\left(\begin{smallmatrix}1&0\\0&p^\ZZ\end{smallmatrix}\right)$ s'identifient \`a $H^2_{\rm Iw}(\Q,T)$. 
\end{rema}

\subsection{Irr\'eductibilit\'e et admissibilit\'e}\label{resul}

Le but de ce paragraphe est de d\'eduire de tous les r\'esultats pr\'ec\'edents la non nullit\'e, l'irr\'eductibilit\'e (topologique) et l'admissibilit\'e de $\breuil(V)$ (pour $\alpha\ne\beta$).

\begin{coro} \label{cestnonnul}
L'espace de Banach $\breuil(V)$ est non nul.
\end{coro} 

\begin{proof} 
Cela r\'esulte du th\'eor\`eme \ref{phigammabanach} et du corollaire
\ref{dpsinozero} qui implique que $\projlim_\psi \dsharp(T) \neq 0$. 
\end{proof} 
 
Pour $\beta\alpha^{-1}$ non ramifi\'e, le corollaire \ref{cestnonnul} \'etait conjectur\'e (via la proposition \ref{existencereseau}) et d\'emontr\'e pour $k\leq 2p$ si $p\ne 2$ et $k<4$ si $p=2$ dans \cite[\S3.3]{Br1} par un calcul explicite de r\'eseaux.
 
\begin{coro}\label{irreductible} 
Le $\G$-Banach unitaire $\breuil(V)$ est topologiquement irr\'educti\-ble. 
\end{coro} 
\begin{proof} 
Cela r\'esulte du th\'eor\`eme \ref{phigammabanach} et de la proposition
\ref{actgirred}.
\end{proof} 
 
La proposition \ref{actgirred} montre que $\breuil(V)$ est en fait topologiquement irr\'eductible comme $\B$-repr\'esentation.
 
\begin{coro}\label{admissible} 
Le $\G$-Banach unitaire $\breuil(V)$ est admissible. 
\end{coro} 

\begin{proof} 
On ignore si le $\O$-module compact $\varprojlim_{\psi}\dsharp(T)$ est stable par $\K$ dans $\breuil(V)^*$ (via le th\'eor\`eme \ref{phigammabanach}) mais on peut le remplacer par le $\O$-r\'eseau de $\breuil(V)^*$ :
$${\mathcal M}\=\cap_{g\in \K}g(\varprojlim_{\psi}\dsharp(T))\subset \varprojlim_{\psi}\dsharp(T)$$
qui est un sous-$\O[[X]]$-module compact stable par $\G$ dans $\breuil(V)^*$ (on v\'erifie qu'il est stable par $B(\Q)$ en utilisant la d\'ecomposition d'Iwasawa de $\G$). Le $\O$-module ${\mathcal M}$ poss\`ede alors {\it deux} structures naturelles de $\O[[X]]$-modules : l'une est celle d\'ej\`a d\'efinie et l'autre est :
$$(\lambda,v)\in \O[[X]]\times {\mathcal M}\mapsto
\begin{pmatrix}0&1\\1&0\end{pmatrix}\lambda\begin{pmatrix}0&1\\1&0\end{pmatrix}v.$$
La premi\`ere structure est telle que la multiplication par $(1+X)^{\Z}$ correspond \`a l'action de $\left(\begin{smallmatrix}1&\Z\\0&1\end{smallmatrix}\right)$ et la deuxi\`eme est telle que la multiplication par $(1+X)^{\Z}$ correspond \`a l'action de $\left(\begin{smallmatrix}1&0\\\Z&1\end{smallmatrix}\right)$. Soit ${\rm pr}:{\mathcal M}\to \dsharp(T)$ la projection sur la premi\`ere composante et $M\={\rm pr}({\mathcal M})$ : $M$ est un sous-$\O[[X]]$-module (de type fini) de $\dsharp(T)$. Posons ${\mathcal N}\={\rm Ker}({\rm pr})\subsetneq {\mathcal M}$. L'application :
\begin{equation}\label{tordu}
{\mathcal N}\to M,\ v\mapsto {\rm pr}\left(\begin{pmatrix}0&1\\1&0\end{pmatrix}v\right)
\end{equation}
est injective : si $v$ a pour image $0$, sa distribution associ\'ee $\mu_{\alpha}\in B(\alpha)^*\simeq \mathcal{C}^{{\rm val}(\alpha_p)}(\Z,L)^*\oplus \mathcal{C}^{{\rm val}(\alpha_p)}(\Z,L)^*$ par (\ref{versbanach}) et (\ref{f1f2}) est nulle sur les deux copies de $\mathcal{C}^{{\rm val}(\alpha_p)}(\Z,L)$, donc est nulle dans $\breuil(V)^*$. En pensant encore en termes de distributions, on voit que ${\mathcal N}$ est un $\O[[X]]$-module pour la premi\`ere structure mais seulement un $\varphi(\O[[X]])$-module pour la deuxi\`eme structure. De plus, pour cette deuxi\`eme structure, l'injection (\ref{tordu}) est $\varphi(\O[[X]])$-lin\'eaire. Comme $M$ est de type fini sur $\O[[X]]$, donc sur $\varphi(\O[[X]])$, on obtient que le $\varphi(\O[[X]])$-module ${\mathcal N}$ pour la deuxi\`eme action de $\varphi(\O[[X]])$ est de type fini. Fixons maintenant des \'el\'ements $(e_1,\hdots,e_m)\in {\mathcal M}$ (resp. $(f_1,\hdots,f_n)\in {\mathcal N}$) tels que les ${\rm pr}(e_i)$ (resp. les $f_i$) engendrent $M$ sur $\O[[X]]$ (resp. $\mathcal N$ sur $\varphi(\O[[X]])$). Soit $v\in {\mathcal M}$. Il existe $\lambda_1,\cdots,\lambda_m$ dans $\O[[X]]$ tels que $v-\sum \lambda_ie_i\in {\mathcal N}$ et il existe $\mu_1,\cdots,\mu_n$ dans $\varphi(\O[[X]])$ tels que $v-\sum \lambda_ie_i=\sum \left(\begin{smallmatrix}0&1\\1&0\end{smallmatrix}\right)\mu_i\left(\begin{smallmatrix}0&1\\1&0\end{smallmatrix}\right)f_i$. Comme les $\lambda_i$ correspondent \`a l'action d'\'el\'ements de l'alg\`ebre de groupe de $\left(\begin{smallmatrix}1&\Z\\0&1\end{smallmatrix}\right)$ et les $\left(\begin{smallmatrix}0&1\\1&0\end{smallmatrix}\right)
\mu_i\left(\begin{smallmatrix}0&1\\1&0\end{smallmatrix}\right)$ \`a l'action d'\'el\'ements de l'alg\`ebre de groupe de $\left(\begin{smallmatrix}1&0\\p\Z&1\end{smallmatrix}\right)$, on voit que $\mathcal M$ est {\it a fortiori} de type fini sur l'alg\`ebre de groupe de $\K$, d'o\`u l'admissibilit\'e. 
\end{proof} 

Pour $\beta\alpha^{-1}$ non ramifi\'e, les corollaires \ref{irreductible} et \ref{admissible} \'etaient conjectur\'es et d\'emontr\'es par un argument de r\'eduction modulo $p$ pour $k\leq 2p$ (et $k<4$ si $p=2$) dans \cite[\S1.3]{Br2} avec l'hypoth\`ese suppl\'ementaire ${\rm val}(\alpha_p + \beta_p)\ne 1$ pour le premier.

On peut d\'eduire des r\'esultats pr\'ec\'edents deux autres corollaires, l'un sur les r\'eseaux dans $\pi(\alpha)$ et $\pi(\beta)$, l'autre sur les vecteurs localement analytiques dans $\breuil(V)$.

\begin{coro}\label{reseaurigolo}
Supposons $\alpha\ne \beta\nrm$, alors $\pi(\alpha)$ (resp. $\pi(\beta)$) poss\`ede des $\O$-r\'eseaux stables par $\G$ et tous les $\O$-r\'eseaux stables par $\G$ dans $\pi(\alpha)$ (resp. $\pi(\beta)$) sont commensurables entre eux. Supposons $\alpha= \beta\nrm$, alors on a le m\^eme r\'esultat pour $\pi(\alpha)$.
\end{coro}
\begin{proof}
L'existence de tels $\O$-r\'eseaux r\'esulte du corollaire \ref{cestnonnul} et de la proposition \ref{existencereseau}. Pour montrer qu'ils sont tous commensurables entre eux, il est \'equivalent de montrer qu'ils sont tous commensurables aux $\O$-r\'eseaux de type fini sur $\O[\G]$. Le $\O$-dual d'un $\O$-r\'eseau stable par $\G$ est toujours contenu dans le $\O$-dual d'un $\O$-r\'eseau de type fini sur $\O[\G]$. Par le th\'eor\`eme \ref{complete} et le corollaire \ref{admissible}, ce dernier dual est de type fini sur l'alg\`ebre de groupe compl\'et\'ee de $\K$. Comme c'est une alg\`ebre noeth\'erienne, il en est de m\^eme du premier dual. Cela entra\^\i ne que le compl\'et\'e de $\pi(\alpha)$ (ou $\pi(\beta)$ si $\alpha\ne \beta\nrm$) par rapport \`a un $\O$-r\'eseau stable par $\G$ quelconque est aussi admissible, et donc topologiquement isomorphe \`a $\breuil(V)$ par le corollaire \ref{irreductible} et le fait que la cat\'egorie des $\K$-Banach admissibles est ab\'elienne (\cite[\S3]{ST3}). Tous les $\O$-r\'eseaux stables par $\G$ induisent donc des normes \'equivalentes sur $\pi(\alpha)$ (ou $\pi(\beta)$ si $\alpha\ne \beta\nrm$) ce qui ach\`eve la preuve.
\end{proof}

\begin{rema}
Lorsque $\alpha=\beta$, on s'attend \`a ce que les corollaires \ref{cestnonnul} \`a \ref{reseaurigolo} restent vrais (cela se d\'eduit de \cite{Br1} pour $k\leq 2p$ et $k\ne 4$ si $p= 2$ par un calcul direct, cf. \cite[Th.1.3.3]{Br2}), mais on ignore si l'on a encore un isomorphisme $(\varprojlim_{\psi}\dfont(V))^{\rm b}\overset{\sim}{\longrightarrow} \breuil(V)^*$ comme au th\'eor\`eme \ref{phigammabanach}.
\end{rema}

Comme dans \cite[\S7]{ST4}, on note $\breuil(V)_{\rm an}$ le sous-$L$-espace vectoriel de $\breuil(V)$ des vecteurs localement analytiques, i.e. des vecteurs $v\in \breuil(V)$ tels que l'application orbite $\G\to \breuil(V)$, $g\mapsto g\cdot v$ est localement analytique. Il est muni d'une topologie naturelle d'espace localement convexe de type compact (cf. \cite[\S7]{ST4}).

Soit :
$${A}(\alpha)\=\left({\rm Ind}_{\B}^{\G}\alpha\otimes d^{k-2}\beta\nrm^{-1}\right)^{\rm an}$$ 
l'induite parabolique localement analytique au sens de \cite{ST1}. On d\'efinit de m\^eme ${A}(\beta)$ en \'echangeant $\alpha$ et $\beta$. On a des injections naturelles continues $\G$-\'equivariantes ${A}(\alpha)\hookrightarrow B(\alpha)$ et ${A}(\beta)\hookrightarrow B(\beta)$. 

\begin{coro}\label{anal}
Supposons $\alpha\ne \beta\nrm$. On a une injection continue $\G$-\'equiva\-riante :
$${\rm A}(\beta)\oplus_{\pi(\beta)}{\rm A}(\alpha)\hookrightarrow \breuil(V)_{\rm an}$$ 
o\`u $\pi(\beta)$ s'envoie dans $A(\alpha)$ via l'entrelacement (\ref{intertw}).
\end{coro} 
\begin{proof} 
Par \cite[\S4]{ST2}, $\pi(\alpha)$ (resp. $\pi(\beta)$) est le seul sous-objet topologiquement irr\'eductible non nul dans $A(\alpha)$ (resp. $A(\beta)$). Par le th\'eor\`eme \ref{complete}, on d\'eduit que les injections ci-dessus induisent encore des injections ${A}(\alpha)\hookrightarrow B(\alpha)/L(\alpha)$ et ${A}(\beta)\hookrightarrow B(\beta)/L(\beta)$. Le r\'esultat d\'ecoule alors du corollaire \ref{entrepadique}.
\end{proof}

Il est naturel de conjecturer :

\begin{conj}
Supposons $\alpha\ne \beta\nrm$. L'application ${\rm A}(\beta)\oplus_{\pi(\beta)}{\rm A}(\alpha)\hookrightarrow \breuil(V)_{\rm an}$ du corollaire \ref{anal} est un isomorphisme topologique.
\end{conj}

\subsection{Le cas non g\'en\'erique}\label{steinberg}

On ach\`eve ici l'examen complet du cas $\alpha= \beta\nrm$ (relations entre les Banach $B(\beta)/L(\beta)$ et $B(\alpha)/L(\alpha)$, vecteurs localement analytiques).

Rappelons que $L(\beta)\simeq (\beta\circ\det)\otimes_L {\rm Sym}^{k-2}L^2\subset B(\beta)$ s'identifie au sous-espace des polyn\^omes de degr\'e $\leq k-2$ \`a coefficients dans $L$ (voir \S\ref{definition1}). Notons $K(\beta)\subseteq B(\beta)$ l'adh\'erence du sous-$L$-espace vectoriel engendr\'e par les fonctions de $L(\beta)$ et les fonctions $f:\Q\to L$ de la forme :
\begin{eqnarray}\label{valuation}
f(z)=\sum_{j \in J}\lambda_j(z-z_j)^{n_j}\val(z-z_j)
\end{eqnarray}
o\`u $J$ est un ensemble fini, $\lambda_j\in L$, $z_j\in \Q$, $n_j\in \{\lfloor\frac{k-2}{2}\rfloor +1,\ldots,k-2\}$ et ${\rm deg}(\sum_{j\in J}\lambda_j(z-z_j)^{n_j})<(k-2)/2$. Pour que $K(\beta)$ soit bien contenu dans $B(\beta)$, il faut v\'erifier le lemme suivant, dont on laisse les d\'etails au lecteur (voir par exemple \cite[lemmes 3.3.1 et 3.3.2]{Br3}) :

\begin{lemm}
Les fonctions $f$ comme en (\ref{valuation}) appartiennent \`a $B(\beta)$.
\end{lemm}

La proposition suivante donne pr\'ecis\'ement le d\'efaut pour l'entrelacement $\widehat I$ du corollaire \ref{entrepadique} d'\^etre un isomorphisme dans ce cas.

\begin{prop}\label{noyaust}
On a une suite exacte $\G$-\'equivariante d'espaces de Banach :
$$0\longrightarrow K(\beta)/L(\beta)\longrightarrow B(\beta)/L(\beta)\overset{\widehat{I}}{\longrightarrow} B(\alpha)/L(\alpha)\longrightarrow 0$$
o\`u $\widehat I$ est le morphisme du corollaire \ref{entrepadique}.
\end{prop}

\begin{proof}
Notons $\st$ la repr\'esentation de Steinberg de $\G$, c'est-\`a-dire $\left({\rm Ind}_{\B}^{\G}1\right)/1$. On a des extensions de repr\'esentations localement alg\'ebriques de $\G$ :
$$0\to L(\beta)\otimes_L\st\to \pi(\alpha)\to L(\beta)\to 0$$
et :
$$0\to L(\beta)\to \pi(\beta)\to L(\beta)\otimes_L\st\to 0.$$
L'entrelacement $\pi(\beta)/L(\beta)\to \pi(\alpha)$ induit par $I$ (cf. \S\ref{definition1}) n'est autre dans ce cas que l'injection $L(\beta)\otimes_L\st\hookrightarrow \pi(\alpha)$. En proc\'edant comme dans \cite[\S\S 2.1-2.2]{Br2}, on v\'erifie que $\pi(\alpha)$ s'identifie aux fonctions $H:\Q\to L$ localement polynomiales de degr\'e au plus $k-2$ telles que, pour $\mid\!z| \gg 0$, on a $H(z)=Q(z)-2P(z){\rm val}(z)$ o\`u $P$ et $Q$ sont des polyn\^omes en $z$ de degr\'e au plus $k-2$ et o\`u l'action de $\G$ est donn\'ee par :
\begin{multline*}
\left[\begin{pmatrix}a & b \\ c & d\end{pmatrix}\cdot H\right](z) \\ 
\=\beta(ad-bc)(-cz+a)^{k-2}\left[H\left(\frac{dz-b}{-cz+a}\right)+
P\left(\frac{dz-b}{-cz+a}\right){\rm val}\left(\frac{ad-bc}{(-cz+a)^2}\right)\right]
\end{multline*}
(prolong\'e par continuit\'e en $z$ tel que $-cz+a=0$). Dans cette identification, la sous-repr\'esentation $L(\beta)\otimes_L\st$ correspond au sous-espace des fonctions $H$ telles que $H(z)=Q(z)$ pour $\mid\!z| \gg 0$ (i.e. $P=0$). Le compl\'et\'e de $\pi(\alpha)$ par rapport \`a un quelconque $\O$-r\'eseau invariant de type fini sur $\O[\G]$ (s'il en existe) se calcule alors par dualit\'e exactement comme dans la preuve de \cite[th\'eor\`eme 3.3.3]{Br3}, en rempla\c cant partout les $\log_{\mathcal L}(z)$ par des ${\rm val}(z)$. En particulier, on obtient que l'injection $L(\beta)\otimes_L\st\hookrightarrow \pi(\alpha)$ induit une surjection sur les compl\'et\'es par rapport \`a des r\'eseaux invariants de type fini, i.e. l'application $\widehat I:B(\beta)/L(\beta)\to B(\alpha)/L(\alpha)$ est surjective, et que le noyau de cette surjection est exactement $K(\beta)/L(\beta)$.
\end{proof}

Lorsque $\alpha=\beta\nrm$, le Banach $\breuil(V)$ admet donc {\it trois} descriptions diff\'erentes. La premi\`ere comme $B(\alpha)/L(\alpha)$, la deuxi\`eme comme compl\'et\'e de $\pi(\alpha)$ et la troisi\`eme comme $B(\beta)/K(\beta)$. En fait, dans ce cas, l'isomorphisme $B(\beta)/K(\beta)\overset{\sim}{\longrightarrow} B(\alpha)/L(\alpha)$ de la proposition \ref{noyaust} doit \^etre vu comme rempla\c cant l'isomorphisme $B(\beta)/L(\beta)\overset{\sim}{\longrightarrow}  B(\alpha)/L(\alpha)$ du cas $\alpha\ne\beta\nrm$.

Concernant les vecteurs localement analytiques dans $\breuil(V)$, on a le r\'esultat suivant dont la preuve est analogue \`a celle du corollaire \ref{anal} en rempla\c cant l'isomorphisme $B(\beta)/L(\beta)\simeq B(\alpha)/L(\alpha)$ par l'isomorphisme $B(\beta)/K(\beta)\simeq B(\alpha)/L(\alpha)$ :

\begin{coro}\label{anal2}
Supposons $\alpha=\beta\nrm$. On a une injection continue $\G$-\'equivariante :
$$A(\beta)/L(\beta) \oplus_{L(\beta)\otimes_L\st}A(\alpha)\hookrightarrow \breuil(V)_{\rm an}.$$
\end{coro}

Comme en \S\ref{resul}, on termine avec la :

\begin{conj}
Supposons $\alpha= \beta\nrm$. L'application $A(\beta)/L(\beta) \oplus_{L(\beta)\otimes_L\st}A(\alpha)\hookrightarrow \breuil(V)_{\rm an}$ du corollaire \ref{anal2} est un isomorphisme topologique.
\end{conj}


\begin{thebibliography}{Ber04b}
 
\bibitem[Ber02]{Ber1}
\textsc{L. Berger} --
\textit{Repr\'esentations $p$-adiques et \'equations diff\'erentielles.}
Invent. Math. 148 (2002), no. 2, 219--284. 

\bibitem[Ber04a]{Be1} 
\textsc{L. Berger} --
\textit{Limites de repr\'esen\-tations cristallines}.
Compos. Math. 140 (2004), no. 6, 1473--1498.

\bibitem[Ber04b]{Be2} 
\textsc{L. Berger} --
\textit{\'Equations diff\'erentielles $p$-adiques et $(\varphi,N)$-modules filtr\'es}. 
Pr\'e\-publication 2004.

\bibitem[BB04]{BB}
\textsc{L. Berger, C. Breuil} --
\textit{Towards a $p$-adic Langlands programme}. 
Notes d'un cours donn\'e \`a l'\'Ecole d'\'et\'e de
Hangzhou (ao\^ut 2004), disponibles \`a l'adresse : 
\texttt{www.ihes.fr/${}^{\sim}$breuil/publications.html}

\bibitem[Bre03a]{Br1} 
\textsc{C. Breuil} -- 
\textit{Sur quelques repr\'esen\-tations modulaires 
et $p$-adiques de $\mathrm{GL}_2(\Qp)$ II}.
J. Institut Math. Jussieu 2, 2003, 23--58.
 
\bibitem[Bre03b]{Br2} 
\textsc{C. Breuil} -- 
\textit{Invariant $\mathcal{L}$ et s\'erie sp\'eciale $p$-adique}.
Ann. Sci. Ecole Norm. Sup. (4) 37 (2004), no. 4, 559--610.
 
\bibitem[Bre03c]{Br3} 
\textsc{C. Breuil} -- 
\textit{S\'erie sp\'eciale $p$-adique et cohomologie \'etale compl\'et\'ee}. Pr\'epublication 2003, disponible \`a l'adresse : \texttt{www.ihes.fr/${}^{\sim}$breuil/publications.html}

\bibitem[Bum98]{Bu} 
\textsc{D. Bump} --
\textit{Automorphic forms and representations}. 
Cambridge Studies in Advanced Math. 55, Cambridge University Press, 1998.
  
\bibitem[CC98]{CC98} 
\textsc{F. Cherbonnier, P. Colmez} --
\textit{Repr{\'e}sentations $p$-adiques surconvergentes}.
Invent. Math. 133 (1998), no. 3, 581--611.

\bibitem[CC99]{CC99} 
\textsc{F. Cherbonnier, P. Colmez} --
\textit{Th{\'e}orie d'Iwasawa des repr{\'e}sentations $p$-adiques d'un corps local}. 
J. Amer. Math. Soc. 12, 1999, 241--268.

\bibitem[Col99]{Co4}
\textsc{P. Colmez} --
\textit{Repr\'esentations cristallines et repr\'esen\-tations 
de hauteur finie.}
J. Reine Angew. Math. 514 (1999), 119--143.

\bibitem[Col04a]{Co2} 
\textsc{P. Colmez} --
\textit{Une correspondance de Langlands locale $p$-adique pour les repr\'esen\-tations semi-stables de dimension $2$}. 
Pr\'epublication 2004.

\bibitem[Col04b]{Co3} 
\textsc{P. Colmez} --
\textit{S\'erie principale unitaire pour $\G$ et repr\'esen\-tations triangulines de dimension $2$}.
Pr\'epublication 2004.

\bibitem[Col05]{Co5} 
\textsc{P. Colmez} --
\textit{Fonctions d'une variable $p$-adique}.
Pr\'epublication 2005.

\bibitem[CF00]{CF} 
\textsc{P. Colmez, J-M. Fontaine} -- 
\textit{Construction des repr\'esen\-tations $p$-adiques semi-stables}.
Inv. Math. 140, 2000, 1--43.

\bibitem[Eme04]{Em} 
\textsc{M. Emerton} --
\textit{$p$-adic $L$-functions and unitary completions of representations of $p$-adic reductive groups}. 
Pr\'epublication 2004.

\bibitem[Eme05]{Em2} 
\textsc{M. Emerton} --
\textit{A local-global compatibility conjecture in the $p$-adic Langlands programme for $\mathrm{GL}_{2/\mathbf{Q}}$}.
Pr\'epublication 2005. 

\bibitem[Fon90]{F90} 
\textsc{J-M. Fontaine} --
\textit{Repr{\'e}sentations $p$-adiques des corps locaux I}. 
The Grothendieck Festschrift, Vol. II, 249--309, Progr. Math. 87, Birkh{\"a}user Boston, Boston, MA 1990.
 
\bibitem[Fon94a]{F3}
\textsc{J.-M. Fontaine} --
\textit{Le corps des p\'eriodes $p$-adiques}.
Ast\'erisque No. 223 (1994), 59--111.

\bibitem[Fon94b]{F4}
\textsc{J.-M. Fontaine} --
\textit{Repr\'esentations $p$-adiques semi-stables}.
Ast\'erisque No. 223 (1994), 113--184.

\bibitem[Fon94c]{F5}
\textsc{J.-M. Fontaine} --
\textit{Repr\'esentations $\ell$-adiques potentiellement semi-stables}.
Ast\'e\-risque No. 223 (1994), 321--347.

\bibitem[Ked04]{KK04} 
\textsc{K. Kedlaya} --
\textit{A $p$-adic local monodromy theorem}. 
Ann. of Math. (2) 160 (2004), no. 1, 93--184.

\bibitem[Per94]{BP94} 
\textsc{B. Perrin-Riou} --
\textit{Th\'eorie d'Iwasawa des repr\'esen\-tations $p$-adiques sur un corps local}. 
Inv. Math. 115, 1994, 81--161.

\bibitem[Sch01]{Sc} 
\textsc{P. Schneider} -- 
\textit{Nonarchimedean Functional Analysis}. 
Springer-Verlag, 2001.

\bibitem[ST01]{ST2} 
\textsc{P. Schneider, J. Teitelbaum} --
\textit{$U(\mathfrak g)$-finite locally analytic representations} (with an appendix by D. Prasad). 
Representation Theory 5, 2001, 111--128.
 
\bibitem[ST02a]{ST1} 
\textsc{P. Schneider, J. Teitelbaum} --
\textit{Locally analytic distributions and $p$-adic representation theory, with applications to ${\rm GL}_{2}$}. 
J. Amer. Math. Soc. 15, 2002, 443--468.

\bibitem[ST02b]{ST3} 
\textsc{P. Schneider, J. Teitelbaum} --
\textit{Banach space representations and Iwasawa theory}. 
Israel J. Math. 127, 2002, 359--380.
 
\bibitem[ST03]{ST4} 
\textsc{P. Schneider, J. Teitelbaum} --
\textit{Algebras of $p$-adic distributions and admissible representations}. 
Inv. Math. 153, 2003, 145--196.

\bibitem[Sen80]{Sn80}
\textsc{S. Sen} --
\textit{Continuous cohomology and $p$-adic Galois representations}.
Inv. Math. 62 (1980/81) 89--116.

\bibitem[Wa96]{W96} 
\textsc{N. Wach} -- 
\textit{Repr{\'e}sentations $p$-adiques potentiellement cristallines}. 
Bull. Soc. Math. France 124, 1996, 375--400.

\end{thebibliography}
\end{document}